\documentclass[12pt]{amsbook}

\usepackage[latin1]{inputenc}
\usepackage[active]{srcltx}

\usepackage[all]{xy}

\usepackage{cases}

\usepackage{enumerate} \usepackage{float} \textwidth=6truein
\usepackage{amssymb}
\newcommand{\Z}{{\mathbb Z}} \newcommand{\Q}{{\mathbb Q}}

 \newcommand{\pcom}{_{p}^{\wedge}}

\newcommand{\liminv}{\operatornamewithlimits{\hbox{$\varprojlim$}}}
\newcommand{\limdir}{\operatornamewithlimits{\hbox{$\varinjlim$}}}

\newcommand{\hocolim}{\operatornamewithlimits{hocolim}}

\newcommand{\Iso}{\operatorname{Iso}\nolimits}
\newcommand{\Hom}{\operatorname{Hom}\nolimits}

\newcommand{\Aut}{\operatorname{Aut}\nolimits}

\renewcommand{\ker}{\operatorname{Ker}\nolimits}
\newcommand{\im}{\operatorname{Im}\nolimits}
\newcommand{\coker}{\operatorname{Coker}\nolimits}
\newcommand{\coim}{\operatorname{Coim}\nolimits}

\renewcommand{\dim}{\operatorname{rk}\nolimits}

\newcommand{\Set}{\operatorname{Set}\nolimits}
\newcommand{\Grp}{\operatorname{Grp}\nolimits}
\newcommand{\Ab}{\operatorname{Ab}\nolimits}
\newcommand{\Cat}{\operatorname{Cat}\nolimits}
\newcommand{\SSet}{\operatorname{SSet}\nolimits}
\newcommand{\Top}{\operatorname{Top}\nolimits}

\newcommand{\Ob}{\operatorname{Ob}\nolimits}

\newcommand{\A}{\ifmmode{\mathcal{A}}\else${\mathcal{A}}$\fi}
\newcommand{\B}{\ifmmode{\mathcal{B}}\else${\mathcal{B}}$\fi}
\newcommand{\C}{\ifmmode{\mathcal{C}}\else${\mathcal{C}}$\fi}
\newcommand{\D}{\ifmmode{\mathcal{D}}\else${\mathcal{D}}$\fi}
\newcommand{\G}{\ifmmode{\mathcal{G}}\else${\mathcal{G}}$\fi}
\newcommand{\I}{\ifmmode{\mathcal{I}}\else${\mathcal{I}}$\fi}
\newcommand{\J}{\ifmmode{\mathcal{J}}\else${\mathcal{J}}$\fi}
\newcommand{\K}{\ifmmode{\mathcal{K}}\else${\mathcal{K}}$\fi}
\renewcommand{\O}{\ifmmode{\mathcal{O}}\else${\mathcal{O}}$\fi}
\renewcommand{\P}{\ifmmode{\mathcal{P}}\else${\mathcal{P}}$\fi}
\newcommand{\U}{\ifmmode{\mathcal{U}}\else${\mathcal{U}}$\fi}
\newcommand{\M}{\ifmmode{\mathcal{M}}\else${\mathcal{M}}$\fi}
\newcommand{\N}{\ifmmode{\mathcal{N}}\else${\mathcal{N}}$\fi}
\newcommand{\Ss}{\ifmmode{\mathcal{S}}\else${\mathcal{S}}$\fi}
\newcommand{\T}{\ifmmode{\mathcal{T}}\else${\mathcal{T}}$\fi}
\newcommand{\Ff}{\ifmmode{\mathcal{F}}\else${\mathcal{F}}$\fi}
\newcommand{\Ll}{\ifmmode{\mathcal{L}}\else${\mathcal{L}}$\fi}
\newtheorem{Thm}{Theorem}[section]
\newtheorem{Prop}[Thm]{Proposition}
\newtheorem{Cor}[Thm]{Corollary}
\newtheorem{Lem}[Thm]{Lemma}
\newtheorem{Claim}{Claim}[Thm]

\theoremstyle{definition}
\newtheorem{Defi}[Thm]{Definition}
\newtheorem{Rmk}[Thm]{Remark}
\newtheorem{Ex}[Thm]{Example}

\theoremstyle{remark}

\numberwithin{section}{chapter}

\newtheorem{Intro_Thm}{Theorem}
\newtheorem{Intro_Defi}[Intro_Thm]{Definition}
\newtheorem{Intro_Prop}[Intro_Thm]{Proposition}
\newtheorem{Intro_Ex}[Intro_Thm]{Example}
\newtheorem{Intro_Lem}[Intro_Thm]{Lemma}
\newtheorem{Intro_Prob}[Intro_Thm]{Problem}

\theoremstyle{plain}

\newcommand{\definicio}{\stackrel{\text{def}}{=}}

\title{Homological Algebra on Graded Posets}

\author{Antonio D{\'\i}az Ramos}

\begin{document}

\maketitle


\tableofcontents

\newpage

\chapter*{Introduction}

This work has its origin in a problem related with $p$-local
finite groups \cite{blo2}. These algebraic objects are triples
$(S,\Ff,\Ll)$ where $S$ is a finite $p$-group and $\Ff$ and $\Ll$
are categories which encode fusion information among the subgroups
of $S$. In fact, they contain the data needed to build topological
spaces which behave as $p$-completions of classifying spaces of
finite groups ($BG\pcom$). While a finite group give rise to a
$p$-local finite group in a natural way, there are $p$-local
finite groups which do not arise from finite groups. These are
called \emph{exotic} $p$-local finite groups, and examples have
already been found in \cite{blo2}, \cite{2-solomon}, \cite{rv} and
\cite{drv}.

Prof. A. Viruel wondered a few years ago if these exotic $p$-local
finite group arise or not from an infinite group. The work of
Aschbacher and Chermak shows that this is the case for some of the
Solomon $2$-local finite groups \cite{2-solomon}.

An approach to find a candidate to this infinite group is to use
the normalizer decomposition for a $p$-local finite group
$(S,\Ff,\Ll)$ \cite{libman}. It gives a description of the nerve
$|\Ll|$ as a homotopy colimit (in the sense of Bousfiled and Kan
\cite{BK}) of spaces which has the homotopy type of classifying
spaces of finite groups:
$$
\hocolim_{\overline{s}dC} \delta_C \simeq |\Ll|.
$$
Here, $\overline{s}dC$ is a partially ordered set (poset for
short) built from the $p$-local finite group $(S,\Ff,\Ll)$ and
$\delta_C:\overline{s}dC\rightarrow \Top$ is a functor which
satisfies
$$
\delta_C(\sigma)\simeq B\Aut_\Ll(\sigma)
$$
for every $\sigma\in \overline{s}dC$ (where $B:\Grp\rightarrow
\Top$ is a functor which sends a discrete small group to a
classifying space for it). If this homotopy colimit has the
homotopy type of an Eilenberg-MacLane space
$$
\hocolim_{\overline{s}dC} \delta_C \simeq B\pi,
$$
then $\pi$ is a candidate for infinite group from which the
original $p$-local finite group arises. Thus, the problem we have
focused on is:
\begin{Intro_Prob}\label{Intro_Prob_hocolim}
 Given a diagram of groups $G:\C\rightarrow \Grp$ and a cone
$\tau:G\Rightarrow \pi$, when is the natural map
$$\hocolim_\C BG\rightarrow B\pi$$
a homotopy equivalence? And, in case it is not, can we compute the
fiber of this map?
\end{Intro_Prob}
A classical result of J.H.C. Whitehead states that if $\C$ is the
pushout category
$$
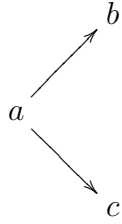
\xymatrix{&b\\
a\ar[ru]\ar[rd]&\\
&c }
$$
and $\pi$ is the amalgamated product $G(b)*_{G(a)}G(c)$, i.e., the
direct limit $\limdir G$, then $\hocolim_\C BG\simeq B\pi$. Thus,
in this case we have an affirmative answer to our problem. It is
worthwhile noticing that, also in this particular case, we have
that the category $\C$ is contractible and that the arrows
$\tau_a$, $\tau_b$ and $\tau_c$ are monomorphisms:
$$
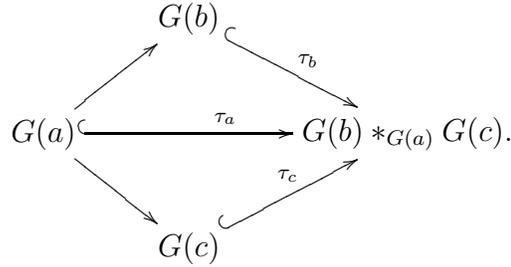
\xymatrix{&G(b)\ar@{^(->}[rd]^{\tau_b}\\
G(a)\ar@{^(->}[rr]^{\tau_a}\ar[ru]\ar[rd]&&G(b)*_{G(a)}G(c).\\
&G(c)\ar@{^(->}[ru]^{\tau_c}}
$$

Studying further the fiber in the general case we find
that:
\begin{Intro_Thm}\label{Intro_Thm_hocolim} Let $G:\C\rightarrow \Grp$ be a functor and
let $\tau: G\Rightarrow \limdir G$ be the limit cone to the direct
limit $\limdir G$. Assume $\C$ is contractible and $\tau$ is a
monomorphism on each object. Then, if $F$ is the fiber of the map
$$\hocolim_\C BG\rightarrow B\limdir G,$$
we have that:
\begin{itemize}
\item $F$ is simply connected, and \item  $H_j(F)=\limdir_{j-1} H$
for each $j\geq 2$, where $H:\C\rightarrow \Ab$ is a
functor.\end{itemize}
\end{Intro_Thm}
Thus, when we assume the same hypothesis as
in Whitehead's Theorem, it turns out that the homology groups are
given by the higher direct limits $\limdir_j H$ for a functor
$H:\C\rightarrow \Ab$ from $\C$ to abelian groups. Before continue
put
\begin{Intro_Defi}\label{Intro_Defi_direct acyclic} Let $H:\C\rightarrow \Ab$ be a functor. Then
we say that it is \emph{$\limdir$-acyclic} if $\limdir_j H=0$ for
$j\geq 1$.
\end{Intro_Defi}
Then it is clear by G.W. Whitehead's Theorems \cite{whitehead},
that $F$ is homotopic to a point, i.e., $\hocolim_\C BG\simeq
B\limdir G$, if and only if $H$ is $\limdir$-acyclic. Thus, we
have reduced the Homotopy Theory Problem \ref{Intro_Prob_hocolim}
to the following problem of Homological Algebra:
\begin{Intro_Prob}\label{Intro_Prob_homological algebra_direct}
Given a functor $H:\C\rightarrow \Ab$, when does it hold that $H$
is $\limdir$-acyclic? And, in case it does not hold, can we
compute $\limdir_j H$?
\end{Intro_Prob}
This is a very general problem, and we shall assume one additional
hypothesis to attack it. It consist in restricting the category
$\C$: mainly because the category $\overline{s}dC$ considered above
is a graded partially ordered set (a graded poset for short) we
assume that $\C$ is a graded poset. These are special posets in
which we can assign an integer to each object (called the degree
of the object) in such a way that preceding elements are assigned
integers which differs just in $1$. Thus a graded poset can be
divided into a set of ``layers" (the objects of a fixed degree),
and these layers are linearly ordered. This restriction is no so
hard as simplicial complexes and subdivision categories are graded
posets. Moreover, each $CW$-complex is (strong) homotopy
equivalent to a graded poset.

Coming back to Problem \ref{Intro_Prob_homological algebra_direct}
consider the functor $H:\C\rightarrow \Ab$ as an object of the
abelian category $\Ab^\C$ of functors from $\C$ to abelian groups.
Because the higher limits $\limdir_j H$ are the left derived
functors of the direct limit functor $\limdir:\Ab^\C\rightarrow \C$,
then $H$ projective is enough for $H$ $\limdir$-acyclic (as left
derived functors vanishes on projective objects). This is a
partial answer to Problem \ref{Intro_Prob_homological
algebra_direct}, inasmuch that we know the projective objects of
$\Ab^\C$. In order to study these projective objects we introduce
\begin{Intro_Defi}\label{Intro_Defi_coker}
Let $H:\C\rightarrow \Ab$ be a functor over the graded poset $\C$,
and let $i_0$ be an object of $\C$. Then define
$$
\coker(i_0)=H(i_0)/\langle \{H(\alpha)|\alpha:i\rightarrow i_0,\alpha\neq
1_{i_0}\}\rangle.$$
\end{Intro_Defi}
That is, $\coker(i_0)$ is the quotient of the value of $H$ on
$i_0$ by the images of the non-trivial morphisms arriving to
$i_0$. For example, if $i_0$ is a minimal object of the graded
poset $\C$, i.e., it has no non-trivial arriving arrows, then
$\coker(i_0)=H(i_0)$. The next theorem relates the projectiveness
of $H$ in $\Ab^\C$ with that of $\coker(i_0)$ in $\Ab$:
\begin{Intro_Thm}\label{Intro_Thm_projectives of Ab^C}
Let $\C$ be a graded poset and suppose $H\in \Ab^\C$. Then $H$ is
projective in $\Ab^\C$ if and only if
\begin{itemize}
\item $\coker(i_0)$ is projective in $\Ab$ for each object $i_0$
of $\C$. \item $H$ is pseudo-projective.
\end{itemize}
\end{Intro_Thm}

This theorem characterizes the projective objects of $\Ab^\C$.
Although the second condition, pseudo-projectiveness, is
technical, and thus I do not state it here for simplicity, both
conditions in the theorem are easy to check for a given functor
$H\in \Ab^\C$. As we commented above, $H$ projective implies that
$H$ is $\limdir$-acyclic, but this last condition is clearly
weaker than $H$ being projective. Can we weak the condition of
projectiveness and still have $\limdir$-acyclicity? Then answer is
yes and next theorem provide us the weaker condition:
\begin{Intro_Thm}\label{Intro_Thm_pseudo implica acyclic direct}
Let $\C$ be a graded poset and suppose $H\in \Ab^\C$. If $H$ is
pseudo-projective then $H$ is $\limdir$-acyclic.
\end{Intro_Thm}
Thus pseudo-projectiveness gives also an answer to Problem
\ref{Intro_Prob_homological algebra_direct}. This last theorem is
proven by constructing a spectral sequence for a given functor
$H\in \Ab^\C$ over a graded poset $\C$. It is the grading of $\C$
which allows us to define certain filtered differential graded
$\Z$-module from which we build the spectral sequence:
\begin{Intro_Prop}\label{Intro_Prop_ss_direct}
Let $\C$ be a graded poset and suppose $H\in \Ab^\C$. Then there is
a homological type spectral sequence $E^*_{*,*}$ with target
$\limdir_* H$.
\end{Intro_Prop}

To answer the second question in Problem
\ref{Intro_Prob_homological algebra_direct} we have made a
``dimension shifting argument" in the sense it is made in group
cohomology \cite{brown}. This is done by constructing a short
exact sequence
$$
0\Rightarrow K_1\Rightarrow H'\Rightarrow H\Rightarrow 0,
$$
where the functor $H'$ is pseudo-projective. Iterating this
process we obtain
\begin{Intro_Lem}
Let $\C$ be a graded poset and $H:\C\rightarrow \Ab$ a functor.
Then, for each $j\geq 1$,
$$\limdir_j H=\limdir_1 K_{j-1},$$
where $H=K_0$, $K_1$, $K_2$,... are functors in $\Ab^\C$.
\end{Intro_Lem}
We have also given an interpretation of the higher limit of order
$1$, $\limdir_1 H$, as a flow problem in the directed graph
associated to $\C$. All these tools can be applied to an example
told to us by A. Libman. This example shows that the conditions
$\C$ contractible and $\tau$ a monomorphism on each object are not
enough to have $\hocolim_\C BG\simeq B\limdir G$ in general. We
compute the fiber of the map $\hocolim_\C BG\rightarrow B\limdir G$
by the methods above:
\begin{Intro_Ex}\label{Intro_Ex_Libman}
For any group $\pi$ consider the graded poset $\C$ and the functor
\linebreak $G:\C\rightarrow \Grp$ with values
$$
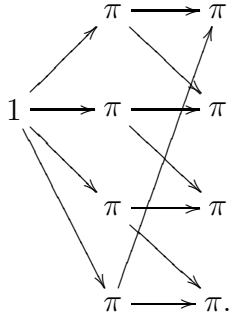
\xymatrix{
      & \pi\ar[r]\ar[rd] & \pi \\
1\ar[r]\ar[ru]\ar[rd]\ar[rdd] & \pi\ar[r]\ar[rd] & \pi\\
      & \pi\ar[r]\ar[rd] & \pi\\
      & \pi\ar[r]\ar[ruuu] & \pi. }
$$
Then there is a fibration
$$\bigvee_{\alpha\in \pi\setminus\{1\}} (S^2)_\alpha\rightarrow \hocolim_\C BG\rightarrow B\pi.$$
\end{Intro_Ex}

All the results concerning higher direct limits are contained in
Chapter \ref{section_on direct limit}, while the spectral sequence
is introduced in Chapter \ref{spectral}. Results about the
homotopy colimit, as Theorem \ref{Intro_Thm_hocolim} and Example
\ref{Intro_Ex_Libman}, are presented in Chapter \ref{hocolim}.
This chapter also contains a proof of the theorem of J.H.C.
Whitehead about the pushout, and a proof of that the classifying
space of a locally finite group is the homotopy colimit of the
classifying spaces of its subgroups. In order to wide our view of
the developments till here we decided to dualize all the
Homological Algebra results. This is made in Chapter
\ref{section_on inverse limit}, and the results are summarized
below:
\begin{Intro_Defi}\label{Intro_Defi_inverse acyclic}
 Let $H:\C\rightarrow \Ab$ be a functor. Then
we say that it is \emph{$\liminv$-acyclic} if $\liminv^j H=0$ for
$j\geq 1$.
\end{Intro_Defi}
\begin{Intro_Defi}\label{Intro_Defi_ker}
Let $H:\C\rightarrow \Ab$ be a functor over the graded poset $\C$,
and let $i_0$ be an object of $\C$. Then define
$$
\ker(i_0)=\bigcap_{\alpha:i_0\rightarrow i,\alpha\neq 1_{i_0}}
H(\alpha).$$
\end{Intro_Defi}
\begin{Intro_Thm}\label{Intro_Thm_injectives of Ab^C}
Let $\C$ be a graded poset and suppose $H\in \Ab^\C$. Then $H$ is
injective in $\Ab^\C$ if and only if
\begin{itemize}
\item $\ker(i_0)$ is injective in $\Ab$ for each object $i_0$ of
$\C$. \item $H$ is pseudo-injective.
\end{itemize}
\end{Intro_Thm}
\begin{Intro_Thm}\label{Intro_Thm_pseudo implica acyclic inverse}
Let $\C$ be a graded poset and suppose $H\in \Ab^\C$. If $H$ is
pseudo-injective then $H$ is $\liminv$-acyclic.
\end{Intro_Thm}
\begin{Intro_Lem}\label{Intro_Lem_C0C1C2...}
Let $\C$ be a graded poset and $H:\C\rightarrow \Ab$ a functor.
Then, for each $j\geq 1$,
$$\liminv_j H=\liminv_1 C_{j-1},$$
where $H=C_0$, $C_1$, $C_2$,... are functors in $\Ab^\C$.
\end{Intro_Lem}
If we take $H=c_\Z:\C\rightarrow \Ab$ in this lemma, i.e., $H$
equals the functor of constant value $\Z$, we obtain an approach
to compute the integer cohomology groups $H^*(|\C|;\Z)=\limdir^*
c_\Z$ of the realization $|\C|$ of the graded poset $\C$. From now
onwards we assume that $H=c_\Z$, and we explain how a sharpener
version of Lemma \ref{Intro_Lem_C0C1C2...} is obtained for
$H=c_\Z$ (this is made in detail in Chapter
\ref{section_cohomology}).

First step is imposing some extra structure in the graded poset
$\C$. In fact, this extra structure is not very restrictive as each
$CW$-complex is still homotopy equivalent to a graded poset which
carries this extra information. We denote this structure by $\J$
and we call it a \emph{covering family} for $\C$. It consist of a
collection of subsets of objects of the category $\C$. We explain
this briefly.

This collection of subsets, $\J$, is of ``local" nature in the
sense that for each object $i_0$ of $\C$ we have several subsets of
objects of $\C$. More precisely, for a fixed $i_0$, we must choose
subsets of objects of the subcategory $(i_0\downarrow \C)$. This
under category corresponds to the full subcategory of $\C$ with
objects $\{i|\exists i_0\rightarrow i\}$, i.e., it is exactly
composed by the objects greater or equal to $i_0$. What we need on
each category $(i_0\downarrow \C)$ amounts to a subset $J^{i_0}_p$
of objects from $(i_0\downarrow \C)$ for each appropriate degree
$p$ (recall that objects of $\C$ are graded). An example shall make
this clearer:
\begin{Intro_Ex}\label{Intro_Ex_cf_1}
Suppose $(i_0\downarrow \C)$ has the following shape
$$
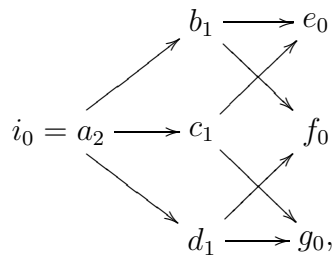
\xymatrix{
& b_1\ar[r]\ar[rd] & e_0\\
i_0=a_2\ar[r]\ar[ru]\ar[rd] & c_1\ar[ru]\ar[rd] & f_0\\
 & d_1\ar[ru]\ar[r] & g_0,  }
$$
where subindexes point out the degree of each object. We can
define $J^{i_0}_2=\{a\}$, $J^{i_0}_1=\{c,d\}$, and
$J^{i_0}_0=\{g\}$.
\end{Intro_Ex}

The collection of subsets $\J$ must fulfill a ``covering"
condition and a ``inheritance" condition. The ``covering"
condition states that each object of $(i_0\downarrow \C)$ of a
given degree $(p-1)$ must be preceded by an object of $J^{i_0}_p$.
The ``inheritance condition" states that if $i\in J^{i_0}_{p+1}$
then $J^i_p\subseteq J^{i_0}_p$. We come back to the earlier
example to show what this means:
\begin{Intro_Ex}\label{Intro_Ex_cf_2}
The ``covering" condition for the object $e$ of degree $0$ if
fulfilled as it is preceded by the object $c\in J^{i_0}_1$:
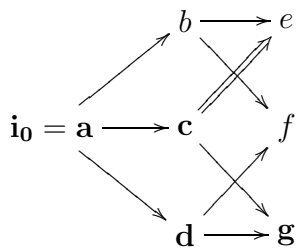
$$
\xymatrix{& b\ar[r]\ar[rd] & e\\
\mathbf{i_0=a}\ar[r]\ar[ru]\ar[rd] & \mathbf{c}\ar@{=>}[ru]\ar[rd] & f\\
 & \mathbf{d}\ar[ru]\ar[r] & \mathbf{g}}
$$
(elements in $\J$ are boldfaced). It is also fulfilled for $f$ and
for $g$ as $d\in J^{i_0}_1$ precedes $f$ and $c$ and $d$ precede
$g$. If we define $J^{d}_1=\{d\}$ and $J^{d}_0=\{g\}$:
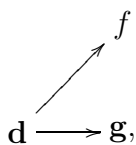
$$
\xymatrix{ &f\\
\mathbf{d}\ar[ru]\ar[r]&\mathbf{g},}
$$
then the object $d\in J^{i_0}_1$ fulfills the ``inheritance"
condition because $J^i_0=\{g\}\subseteq J^{i_0}_0=\{g\}$.
\end{Intro_Ex}

The existence of this family $\J$ of subsets for each object $i_0$
of $\C$, plus a numerical condition (\emph{adequateness}) involving
the number of elements in the subsets of $\J$, gives:
\begin{Intro_Prop}\label{Intro_prop cohomology with adequate covering family}
Let $\C$ be a graded poset category and let $\J$ be an adequate
covering family for $\C$. Then, for each $j\geq 1$,
$$\liminv_j c_\Z=\liminv_1 F_{j-1},$$
where $c_\Z=F_0$, $F_1$, $F_2$,... are functors in $\Ab^\C$ which
take free abelian groups as values.
\end{Intro_Prop}
The sequence of functors $F_0$, $F_1$, $F_2$,... of Proposition
\ref{Intro_prop cohomology with adequate covering family} have
some advantages over the sequence of functors $C_0$, $C_1$,
$C_2$,... of Lemma \ref{Intro_Lem_C0C1C2...} (applied to the
functor $H=c_\Z$):
\begin{itemize}
\item $F_j(i_0)$ is a free abelian group for each object $i_0$ of
$\C$,\item $F_j(i_0)$ has far less generators than $C_j(i_0)$ for
each object $i_0$ of $\C$, and \item $\liminv F_j\simeq \liminv
F_j|{\C_{j+1,j}}$, where $\C_{j+1,j}$ is the full subcategory of $\C$
with objects of degrees $j+1$ and $j$.
\end{itemize}
The third statement above is the main feature of $F_j$. It is a
generalization of the fact that the connected components of a
simplicial complex are computed by just looking to the vertices
and to the edges, i.e., to the objects of degrees $0$ and $1$. The
nice properties of $F_j$ causes our main theorem about integer
cohomology:
\begin{Intro_Thm}\label{Intro_Thm_adequate_local y global characterization acyclicity}
Let $\C$ be a bounded graded poset for which there exist an
adequate covering family $\J$ and an adequate global covering
family $\K$. Then $|\C|$ is acyclic if and only if $|K_0|$ equals
the number of connected components of $\C$. Moreover, in this case
$H^0(|\C|;\Z)=\Z^{|K_0|}$.
\end{Intro_Thm}
Here, we are using the concept of adequate \emph{global} covering
family, which plays the role of covering family for the whole
category $\C$ instead of for each subcategory $(i_0\downarrow \C)$.
The family $\K$ is composed of subsets $K_p$, one for each
appropriate degree $p$. This theorem reduces the
$\liminv$-acyclicity of the functor $c_\Z\in \Ab^\C$, i.e., the
acyclicity of the space $|\C|$, to the integral equation
$|K_0|=|\pi_0(\C)|=|\pi_0(|\C|)|$.

Although the hypotheses in Theorem \ref{Intro_Thm_adequate_local y
global characterization acyclicity} look very hard, the existence
of the family $\J$ is automatic for $\C^{op}$ (of course $\C^{op}$
and $\C$ have the same homotopy type) when $\C$ is a simplicial
complex or a subdivision category. In fact, it applies to $\C^{op}$
whenever the graded poset $\C$ is locally as a simplicial complex,
i.e., such that $(\C\downarrow i_0)$ is isomorphic to the poset of
all non-empty sets of a finite set (with inclusion as order
relation). We call these posets \emph{simplex-like} posets. Thus,
the difficult point resides in finding a global family $\K$.

In Chapter \ref{section_applications Webb_conjecture} we describe
a particular case where this global family $\K$ exists, yielding a
proof of part of a conjecture of Webb. The conjecture is related
with the Brown's complex \cite{brown} (denoted $\Ss_p(G)$) of a
finite group $G$ and a prime $p$. Webb conjectured that the orbit
space $\Ss_p(G)/G$ is contractible. This orbit space has as
objects the $G$-conjugation classes of chains of subgroups of
$\Ss_p(G)$. The conjecture was first proven by Symonds in
\cite{symonds}, generalized for blocks by Barker
\cite{barker1,barker2} and extended to arbitrary (saturated)
fusion system by Linckelmann \cite{markus}.

The works of Symonds and Linckelmann prove the contractibility of
the orbit space by showing that it is simply connected and
acyclic, and invoking G.W. Whitehead's theorem. Both proofs of
acyclicity work on the subposet of normal chains (introduced by
Knörr and Robinson \cite{robinson} for groups). Symonds uses the
results from Thevenaz and Webb \cite{G-homot} that the subposet of
normal chains is $G$-equivalent to Brown's complex. Linckelmann
proves on his own that, also for fusion systems, the orbit space
and the orbit space on the normal chains has the same cohomology
\cite[Theorem 4.7]{markus}.

We apply Theorem \ref{Intro_Thm_adequate_local y global
characterization acyclicity} to prove that the orbit space on the
normal chains is acyclic. The definition of the global covering
family $\K$ is related with the pairing defined by Linckelmann in
\cite[Definition 4.7]{markus}. We describe $\K$ briefly: let
$(S,\Ff)$ be a saturated fusion system where $S$ is a $p$-group
and consider its subdivision category $\Ss(\Ff)$ and the poset
$[\Ss(\Ff)]$. An object of degree $n$ in the orbit space of the
normal chains $[\Ss_\lhd(\Ff)]$ is an $\Ff$-isomorphism class of
chains
$$
[Q_0<Q_1<...<Q_n]
$$
where the $Q_i$ are subgroups of $S$ normal in $Q_n$. The
subcategory \linebreak$([\Ss_\lhd(\Ff)]\downarrow [Q_0<...<Q_n])$
has objects $[Q_{i_0}<...<Q_{i_m}]$ with $0\leq m\leq n$ and
$0\leq i_1<i_2<...<i_m\leq n$ (see \cite[2]{markus}). For example,
\linebreak $([\Ss_\lhd(\Ff)]\downarrow [Q_0<Q_1<Q_2])$ is
$$
\xymatrix{ [Q_0]\ar[r]\ar[rd] & [Q_0<Q_1]\ar[rd]\\
[Q_1]\ar[ru]\ar[rd]&[Q_0<Q_2]\ar[r]&[Q_0<Q_1<Q_2].\\
[Q_2]\ar[r]\ar[ru]&[Q_1<Q_2]\ar[ru] }
$$
Then it is clear that $[\Ss_\lhd(\Ff)]$ is a simplex-like poset
and thus Theorem \ref{Intro_Thm_adequate_local y global
characterization acyclicity} applies to $[\Ss_\lhd(\Ff)]^{op}$.
The definition of the global family $\K$ follows
$$
K_n=\big\{[Q_0<...<Q_n]|[Q_0<...<Q_n]=[Q'_0<...<Q'_n]\Rightarrow
\cap_{i=0}^n N_S(Q'_i)=Q'_n\big\}.
$$
It is easy to see that $[\Ss_\lhd(\Ff)]^{op}$ has just one
connected component. Thus, by Theorem
\ref{Intro_Thm_adequate_local y global characterization
acyclicity}, we can conclude that $[\Ss_\lhd(\Ff)]^{op}$ is
acyclic if it is the case that $|K_0|=1$. By definition
$$
K_0=\big\{[Q_0]|[Q_0]=[Q'_0]\Rightarrow N_S(Q'_0)=Q'_0\big\}.
$$
Because the unique subgroup of the $p$-group $S$ whose normalizer
equals $S$ is $S$ itself we have that $|K_0|=1$.

This finishes the exposition of the main results of the present
work. By chapter, the contents are the following:
\begin{itemize}
\item Chapter 1: Notation and preliminaries. The notation used
throughout the work is introduced, as well as the definitions of
derived functors and the normalization theorems for simplicial
abelian groups. Likewise, graded posets are defined, which are the
categories over which the main results of the work are developed.
This kind of category have been chosen because it is the prototype
of the category  $\overline{s}dC$ used in the normalizer
decomposition mentioned above \cite{libman}. Moreover, they
include simplicial complexes (Section 1.4) and subdivision
categories \cite{markus}, and they contain all the homotopic
information of $CW$-complexes (Section 1.5).

\item Chapter 2: A spectral sequence. A spectral sequence is built
from differential graded modules associated to a functor
$H:\C\rightarrow \Ab$, where $\C$ is a graded poset. The target of
this spectral sequence is the higher limits $\limdir_* H$. In a
similar way another spectral sequence is built whose target is the
higher inverse limits $\liminv^* H$.

\item Chapter 3: Higher direct limits. Projective objects in the
abelian category $\Ab^\C$ are characterized, where $\C$ is a graded
poset. Moreover, due to the spectral sequence from Chapter 2,
another family of objects of $\Ab^\C$ with vanishing higher direct
limits is found, they are the pseudo-projective functors. Finally,
a way of reducing the higher limit $\lim_j H$ of order $j$ to
limits of lower order is given, and some applications of these
results are worked out.

\item Chapter 4: Higher inverse limits. This chapter is the
dualization of Chapter 3, in which the injective objects of
$\Ab^\C$ are described and another family of functors with
vanishing higher inverse limits is found (the pseudo-injective
functors). As application some tools for the computation of
integer cohomology of categories are developed in Chapter 5.

\item Chapter 5: Cohomology. Through some additional structure
(called covering families) on the graded poset $\C$ it is proven
that such a category is acyclic if and only if some integral
equation involving geometrical elements of $\C$ is satisfied.

\item Chapter 6: Application: Webb's conjecture. Results from the
earlier chapter are used to prove the cohomological part of Webb's
conjecture, which has been proven in its maximal generality in
\cite{markus}. This conjecture states that the orbit space of the
$p$-subgroups poset of a finite group $G$ is contractible.

\item Chapter 7: Application: homotopy colimit. Although the study
of homotopy colimits is the origin of this work, it appears here
as an application of earlier chapters. It is proven Whitehead's
Theorem, which states that the pushout of Eilenberg-MacLane spaces
with injective morphisms is an Eilenberg-MacLane space, and an
example is given of explicit computation of the fiber $F$. Also it
is proven the well known fact that the classifying space of a
locally finite group is the homotopy colimit of the classifying
spaces of its finite subgroups.
\end{itemize}

Finally, we comment other applications where we hope the theory
can contribute with some results. First of them is to translate to
the initial diagram of groups $G:\C\rightarrow \Grp$ the conditions
of pseudo-projectiveness over $H$ (the functor $H$ comes from Theorem
\ref{Intro_Thm_hocolim}). In this way it should be obtained a
generalized Whitehead's Theorem for diagrams larger than the
pushout. Moreover, it would be interesting to study its relation
with Mather's cube Theorem.

The second application is take up again the original problem about
$p$-local finite groups. We have already found conditions which
imply that the functor
\linebreak$\delta_C:\overline{s}dC\rightarrow \Top$ factorizes as
$$
\xymatrix{ \overline{s}dC\ar[rd]\ar[rr]^{\delta_C}&& \Top.\\
&\Grp\ar[ru]^{B}}
$$
This makes possible to apply the theory to a ``honest" diagram of
groups. Moreover, these conditions imply that $\overline{s}dC$ is
contractible, fulfilling one of the hypotheses of Theorem
\ref{Intro_Thm_hocolim}. Thus, the next step is determine when is
$|\Ll|$ an Eilenberg-MacLane space or compute the fiber $F$
$$
F\rightarrow |\Ll|\rightarrow B\pi
$$
for appropriate group $\pi$ in favorable cases. For example, it
can be shown that all the $p$-rank $2$ $p$-local finite groups for
$p$ odd (which are described in \cite{drv}) have an
Eilenberg-MacLane space as non-completed classifying space
$|\Ll|$.

The third possible application is related to Quillen's conjecture
\cite{quillen-p}, which states that the poset of the non-trivial
$p$-subgroups of a finite group $G$, $\Ss_p(G)$ (Brown's complex
\cite{brown}), is contractible if and only if $\O_p(G)\neq 1$
(where $\O_p(G)$ is the largest normal $p$-subgroup of
$G$). A stronger formulation of this conjecture is that the
integer reduced cohomology of this poset is trivial if and only if
$\O_p(G)\neq 1$ (this formulation is used in \cite{asch-smith}).
Due to Theorem \ref{Intro_Thm_adequate_local y global
characterization acyclicity} this stronger formulation is
equivalent to find an adequate covering family $\K(G)$ of
$\Ss_p(G)$ for each finite group $G$ which satisfies the condition
$$|K(G)_0|=1\Leftrightarrow \O_p(G)\neq 1.$$
By results of Bouc \cite{bouc} and of Thévenaz and Webb
\cite{G-homot} we can also work with the poset of non-trivial
$p$-radical subgroups or with the poset of non-trivial elementary
abelian $p$-subgroups.

Final note: For reference it follows a list of the main points of
the introduction and their corresponding statements in the manuscript.
For simplicity some statements over the boundedness of the graded
poset $\C$ were omitted in the introduction. These conditions are
automatically satisfied in the applications and for finite graded
posets. We add them bracketed in the list (as well as other
comments).
\begin{itemize}
\item Definition \ref{Intro_Defi_direct acyclic} = Definition
\ref{Defi_direct_acyclic}. \item Definition
\ref{Intro_Defi_coker} = Definition \ref{defi_im}. \item Definition
\ref{Intro_Defi_inverse acyclic} = Definition
\ref{Defi_inverse_acyclic}. \item Definition
\ref{Intro_Defi_ker} = Definition \ref{defi_ker}. \item Definition
of pseudo-projectiveness = Definition
\ref{property_pseudo-projective}. \item Definition of
pseudo-injectiveness = Definition \ref{property_pseudo-injective}.
\item Definition of covering family = Definition
\ref{defi_covering_family}. \item Definition of adequate covering
family = Definition \ref{defi_adequate_covering_family}. \item
Definition of global covering family = Definition \ref{defi_global
covering family}. \item Definition of adequate global covering
family = Definition \ref{defi global_adequate_covering_family}. \item
Theorem \ref{Intro_Thm_hocolim} = Theorem \ref{hocolim_technic} with
$G_0=\limdir G$. \item Theorem \ref{Intro_Thm_projectives of
Ab^C} = Lemma \ref{lem_proj_pseudo_graded_coker_projective} + Lemma
\ref{lem_proj_pseudo_graded} + Proposition
\ref{pro_projective_bounded_graded} ($\C$ must be bounded below).
\item Theorem \ref{Intro_Thm_pseudo implica acyclic
direct} = Theorem \ref{pro_acyclic_graded} ($\C$ must be bounded
below). \item Theorem \ref{Intro_Thm_injectives of Ab^C} = Lemma
\ref{lem_proj_pseudo_graded_ker_injective} + Lemma
\ref{lem_inj_pseudo_graded} + Proposition
\ref{pro_injective_bounded_graded} ($\C$ must be bounded above).
\item Theorem \ref{Intro_Thm_pseudo implica acyclic
inverse} = Theorem \ref{pro_inv_acyclic_graded} ($\C$ must be bounded
above). \item Proposition \ref{Intro_Prop_ss_direct} = Proposition
\ref{spectral_sequence_1} + Proposition \ref{spectral_sequence_2} (also in Chapter \ref{spectral}
is built the dual spectral sequence with target $\liminv^* H$).
\item Proposition \ref{Intro_prop cohomology with adequate
covering family} = Proposition \ref{cohomology with adequate
covering family} ($\C$ must be bounded above). \item Example
\ref{Intro_Ex_Libman} = Example \ref{Ex_Libman}. \item Examples
\ref{Intro_Ex_cf_1} and \ref{Intro_Ex_cf_2} = come from Section
\ref{section_locally delta}. \item Simplex-like posets = Definition
\ref{simplex-like poset} + Section \ref{section_locally
delta} (relation with covering families). \item Graded
posets = Section \ref{section graded_posets} + Section
\ref{realization} (realize all the $CW$-homotopy types). \item
Locally finite groups = Example \ref{ex_locally finite group}. \item
Properties of functors $F_0$, $F_1$, $F_2$... = Proposition
\ref{cohomology with adequate covering family} + Section
\ref{section_coho_local2} + Lemma \ref{lemma_liminvFp
p+2_esqueleto}. \item Whitehead's Theorem = Example
\ref{ex_hocolim_pushout}.
\end{itemize}

\chapter*{Resumen en castellano}

\section*{Introducción.}

Este trabajo tiene su origen en un problema relacionado con grupos
$p$-locales finitos \cite{blo2}. Estos objetos algebraicos son
triples $(S,\Ff,\Ll)$ donde $S$ es un $p$-grupo finito y $\Ff$ y
$\Ll$ son categorías que codifican la información de fusión entre
los subgrupos de $S$. De hecho, contienen los datos necesarios
para construir espacios topológicos que se comportan como
$p$-completaciones de espacios clasificadores de grupos finitos
$BG\pcom$. Mientras que un grupo finito da lugar a un grupo
$p$-local finito de una forma natural, hay grupos $p$-locales
finitos que no surgen de este modo. Son los llamados grupos
$p$-locales finitos \emph{exóticos}, de los cuales se han
encontrado ya varias familias  \cite{blo2}, \cite{2-solomon},
\cite{rv} y \cite{drv}.

El profesor A. Viruel formuló hace varios años la cuestión de si
estos grupos $p$-locales finitos exóticos provienen o no de un
grupo infinito. El trabajo de Aschbacher y Chermak muestra que
esto ocurre en algunos de los grupos finitos $2$-locales de
Solomon \cite{2-solomon}.

Un enfoque para encontrar candidatos a este grupo infinito es usar
la descomposición por normalizadores \cite{libman} para un grupo
$p$-local finito $(S,\Ff,\Ll)$. Esta descomposición da una
descripción del nervio $|\Ll|$ como un colímite homotópico de
espacios del tipo de homotopía de espacios clasificadores de
grupos finitos:
$$
\hocolim_{\overline{s}dC} \delta_C \simeq |\Ll|.
$$
Aquí, $\overline{s}dC$ es cierto poset construido a partir del
grupo $p$-local finito $(S,\Ff,\Ll)$ y \linebreak
$\delta_C:\overline{s}dC\rightarrow \Top$ es un funtor que
satisface
$$
\delta_C(\sigma)\simeq B\Aut_\Ll(\sigma).
$$
para todo $\sigma\in \overline{s}dC$ (donde $B:\Grp\rightarrow
\Top$ asocia funtorialmente a cada grupo discreto un espacio
clasificador). Si este colímite homotópico tiene el tipo de
homotopía de un espacio de Eilenberg-MacLane
$$
\hocolim_{\overline{s}dC} \delta_C \simeq B\pi,
$$
entonces $\pi$ es un candidato a grupo infinito del que proviene
el grupo $p$-local finito original.

Así que el problema del que partió este trabajo es si dado un
diagrama de grupos $G:\C\rightarrow \Grp$ y un cono
$\tau:G\Rightarrow \pi$ la aplicación natural
$$\hocolim_\C BG\rightarrow B\pi$$
es o no una equivalencia homotópica. Con este propósito estudiamos
la fibra $F$ de esta aplicación (Capítulo 7). Bajo hipótesis
débiles $F$ es simplemente conexa y su homología está dada por los
límites directos superiores $\limdir_j H$, donde $H:\C\rightarrow
\Ab$. Hemos estudiado entonces condiciones para que estos límites
derivados se anulen y herramientas para su cómputo (Capítulo 3),
así como los resultados duales para límites inversos superiores
$\liminv^j H$ (Capítulo 4). Estos resultados están fundamentados
en una sucesión espectral (Capítulo 2), y tienen aplicaciones al
cálculo de cohomología entera de categorías (Capítulo 5), como el
caso de la conjetura de Webb (Capítulo 6).

Una versión más detallada del contenido es la siguiente:
\begin{itemize}
\item \textbf{Capítulo 1: Notación y preliminares.} Se introducen
las notaciones que se usarán así como las definiciones adecuadas
de funtores derivados y los teoremas de normalización para grupos
abelianos simpliciales. Así mismo, se introducen los posets
graduados, que son las categorías sobre las que se desarrollan los
resultados principales de este trabajo. Este tipo de categorías ha
sido elegido ya que es el prototipo de la categoría
$\overline{s}dC$ usada en la descomposición antes mencionada
\cite{libman}. Además incluye a los complejos simpliciales
(Sección 1.4) y a las categorías de subdivisión \cite{markus}, y
contiene toda la información homotópica de los $CW$-complejos
(Sección 1.5).

\item \textbf{Capítulo 2: Un sucesión espectral.} Se construye a
partir de módulos diferenciales graduados una sucesión espectral
asociada a un funtor $H:\P\rightarrow \Ab$, donde $\P$ es un poset
graduado. El límite de esta sucesión espectral es los límites
superiores $\limdir_j H$. Análogamente se construye otra sucesión
espectral que converge a los límites inversos superiores
$\liminv^j H$.

\item \textbf{Capítulo 3: Límites directos superiores.} Se
caracterizan los objetos proyectivos en la categoría abeliana
$\Ab^\P$ donde $\P$ es un poset graduado. Además, gracias a la
sucesión espectral del Capítulo 2, se encuentra otra familia de
objetos de $\Ab^\P$ cuyos límites directos superiores se anulan.
Finalmente se muestran aplicaciones de estos resultados.

\item \textbf{Capítulo 4: Límites inversos superiores.} Es la
dualización del Capítulo 3, en el cual se tratan los objetos
inyectivos de $\Ab^\P$ y otra familia de funtores cuyos límites
inversos superiores se anulan. Como aplicación se desarrollan
herramientas para el cálculo de cohomología de categorías en el
Capítulo 5.

\item \textbf{Capítulo 5: Cohomología.} A través de una estructura
adicional sobre el poset graduado $\P$ se establece que dicha
categoría es acíclica si y sólo si se verifica cierta ecuación
entera en la que están involucrados elementos geométricos de $\P$.

\item \textbf{Capítulo 6: Aplicación: la conjetura de Webb.} Se
usan los resultados del capítulo anterior para probar la parte
cohomológica de la conjetura de Webb, que ha sido demostrada en su
máxima generalidad en \cite{markus}. Esta conjetura afirma que el
espacio de órbitas del poset de los $p$-subgrupos de un grupo
finito $G$ es contráctil.

\item \textbf{Capítulo 7: Aplicación: colímite homotópico.} Aunque
el estudio de colímites homotópicos es el origen de este trabajo,
aparece aquí como aplicación de los capítulos anteriores. Se
prueba el Teorema de Whitehead, el cual afirma que el pushout de
espacios de Eilenberg-MacLane con aplicaciones inyectivas es un
espacio de Eilenberg-MacLane, y se da un ejemplo de cálculo
explícito de la fibra $F$. También se demuestra el hecho conocido
de que el espacio clasificador en un grupo localmente finito es el
colímite homotópico de los espacios clasificadores de sus
subgrupos finitos.
\end{itemize}

Finalmente, comentamos otras aplicaciones donde esperamos que la
teoría pueda contribuir con algunos resultados:
\begin{itemize}
\item Trasladar al diagrama inicial de grupos $G:\C\rightarrow
\Grp$ las condiciones sobre $H$ que implican que $\limdir_j H=0$
para $j>0$. De esta forma se obtendría una versión generalizada
del Teorema de Whitehead para diagramas más grandes que el
pushout. Estudiar la relación  con el Teorema del cubo de
Mather.\item Habiendo encontrado ya condiciones que implican que
el funtor $\delta_C:\overline{s}dC\rightarrow \Top$ factoriza como
$$
\xymatrix{ \overline{s}dC\ar[rd]\ar[rr]^{\delta_C}&& \Top,\\
&\Grp\ar[ru]^{B}}
$$
determinar cuando $|\Ll|$ es un espacio de Eilenberg-MacLane o
calcular la fibra $F$
$$
F\rightarrow |\Ll|\rightarrow B\pi
$$
para un grupo $\pi$ adecuado en casos favorables. Por ejemplo, se
puede demostrar que el espacio clasificador sin completar $|\Ll|$
de todos los grupos $p$-locales finitos de $p$-rango $2$ con $p$
impar (descritos en \cite{drv}) es un espacio de
Eilenberg-MacLane. \item Un posible enfoque de la conjetura de
Quillen \cite{quillen-p}, la cual afirma que el poset de los
$p$-subgrupos no triviales de un grupo finito $G$ (complejo de
Brown \cite{brown}) es contráctil si y sólo si $\O_p(G)\neq 1$
(donde $\O_p(G)$ es el mayor $p$-subgrupo normal de $G$). Una
formulación más fuerte de esta conjetura es que la cohomología
entera reducida de este poset es trivial si y sólo si $O_p(G)\neq
1$ (esta formulación es usada, por ejemplo, en \cite{asch-smith}).

Gracias a los resultados del Capítulo 5 esta formulación más
fuerte es consecuencia de encontrar para cada grupo finito $G$ un
subconjunto adecuado $K(G)$ del complejo de Brown's que satisface
la ecuación entera
$$|K(G)|=1$$
exactamente cuando $\O_p(G)\neq 1$. Por resultados de Bouc
\cite{bouc} y de Thévenaz y Webb \cite{G-homot} también podemos
trabajar con el poset de los subgrupos \linebreak $p$-radicales no
triviales o con el poset de los $p$-subgrupos elementales
abelianos no triviales.
\end{itemize}

\section*{Capítulo 1. Notación y preliminares.}
\subsection*{1.1. Categorías.}
En esta sección se enumeran las categorías que se usarán a lo
largo del trabajo, y que son las siguientes: $\Set$ (conjuntos),
$\Grp$ (grupos), $\Ab$ (grupos abelianos), $\Top$ (espacios
topológicos), $\SSet$ (conjuntos simpliciales) y $\Cat$ (categoría
cuyos objetos son las categorías pequeñas y cuyas flechas son los
funtores entre ellas). Los objetos y flechas de una categoría $\C$
se denotan por $\Ob(\C)$ y $\Hom(\C)$ respectivamente.

Después se introducen los conceptos de categoría \emph{conexa} así
como las \emph{sobre-categorías} $(S\downarrow c_0)$ y
\emph{bajo-categorías} $(c_0\downarrow S)$ dado un funtor
$S:\C\rightarrow \D$ y un objeto $c$ de $\C$. Particularizando a
$S=1_\C:\C\rightarrow \C$ se definen las categorías de objetos
sobre $(\C\downarrow c_0)$ y bajo $(c_0\downarrow \C)$ un objeto
dado $c$ de la categoría $\C$.

La principal referencia para esta sección es \cite{MCL}.

\subsection*{1.2. Grupos y grupos abelianos.}
Tras introducir las notaciones multiplicativas y aditivas para
$\Grp$ y $\Ab$ respectivamente, se fijan las notaciones para
\emph{subgrupo generado}, \emph{producto directo} y \emph{suma
directa} (coproducto) en la categoría de grupos abelianos $\Ab$.
Así mismo se introducen ciertos grupos abelianos (infinito cíclico
$\Z$, finito cíclico $\Z_n$, racionales $\Q$, $\Z[p^\infty]$ con
$p$ primo) y algunas definiciones relativas a homomorfismos en
$\Ab$.

Seguidamente se introducen tanto los funtores $F:\C\rightarrow
\Grp$ como las transformaciones naturales entre ellos que se
comportan de manera inyectiva o sobreyectiva sobre cada flecha o
sobre cada objeto respectivamente. De la referencia
\cite{griffith} se extrae algunos hechos básicos de $\Ab$,
incluyendo la descripción exacta de sus objetos proyectivos e
inyectivos.

Finalmente se introducen las aplicaciones $f:A\rightarrow B$ de
$\Ab$ donde $A$ y $B$ son libres, y tales que el conúcleo $B/f(A)$
es libre. Una aplicación inyectiva de este tipo es un isomorfismo
cuando $A$ y $B$ tienen el mismo rango finito ($\dim A=\dim B$).

\subsection*{1.3. Posets graduados.}
En esta sección se introduce el tipo de categoría sobre el que se
establecerán los principales resultados de este trabajo. No son
sino ciertos conjuntos parcialmente ordenados (\emph{posets}) para
los que existe una \emph{función de grado} $deg:\Ob(\P)\rightarrow
\Z$ cuyos valores sobre objetos que se preceden se diferencian en
una unidad.

Tras varios ejemplos y algunos hechos básicos se prueba que la
función de grado se puede extender a los morfismos del poset
graduado, y se definen los subconjuntos $\Ob_S(\C)$ y $\Hom_S(\C)$
como los objetos y flechas, respectivamente, cuyo \emph{grado}
pertenece a un subconjunto fijado $S\subseteq \Z$.

\subsubsection*{Acotación en posets graduados.}
En esta subsección se definen los posets gra\-duados $\P$ cuya
función de grado $deg:\Ob(\P)\rightarrow \Z$ tiene imagen acotada
superior o inferiormente o ambas. Varios ejemplos ilustran los
diferentes casos, así como la no acotación.

\subsubsection*{Grafos y posets graduados} Se introduce el grafo
(dirigido) asociado a un poset graduado $\P$ como aquél que tiene
los mismos objetos y tiene por aristas las flechas de grado $1$ de
$\P$, es decir, las flechas que corresponden a dos elementos que
se preceden en $\P$. A partir de aquí se introducen los conceptos
de árbol y árbol maximal para posets graduados.

\subsection*{1.4. Complejos simpliciales.} Después de definir los objetos
geométricos que constituyen los complejos simpliciales, se los
caracteriza como aquellos posets que localmente son
\emph{simplex-like}, y tales que cada par de elementos que tienen
una cota inferior poseen un ínfimo \cite{garrett}. Además, se
muestra que todo conjunto simplicial puede ser visto como un poset
graduado.

\subsection*{1.5. Realización y tipo de homotopía.}
En este punto se ha visto que tenemos las inclusiones
$$
\textit{Conjuntos simpliciales }\subseteq\textit{Posets graduados
}\subseteq \textit{Categorías}.
$$
Usando que todo espacio topológico tiene el tipo de homotopía
(débil) de un $CW$-complejo \cite{whitehead}, los cuales a su vez
tienen el tipo de homotopía de un complejo simplicial
\cite{hatcher}, y usando también la \emph{subdivisión
baricéntrica} \cite{spanier} se llega a\\

\textbf{Proposición.} \textit{Para cualquier espacio $X\in \Top$
existe un poset graduado $\P$ y una equivalencia de homotopía
débil:
$$
X\simeq |N\P|.
$$
}

Esto significa que la restricción a posets graduados no implica
pérdida alguna desde el punto de vista homotópico.

\subsection*{1.6. Functores derivados del límite directo e inverso.}
En esta sección se dan definiciones   para los funtores derivados
izquierdos $\limdir_i$ del funtor límite directo
$\limdir:\Ab^\C\rightarrow \C$ y para los funtores derivados
derechos $\liminv^i$ del límite inverso $\liminv:\Ab^\C\rightarrow
\C$, donde $\C$ es una categoría pequeña.

Sabemos \cite[IX, Proposition 3.1]{maclanehom} que $\Ab^\C$ es una
categoría abeliana, y que $\limdir$ y $\liminv$ son exactos por la
derecha y por la izquierda respectivamente porque son los adjuntos
por la izquierda y por la derecha del funtor $\Ab\rightarrow
\Ab^\C$ definido por $A\mapsto c_A$, donde $c_A:\C\rightarrow \Ab$
es el funtor de valor constante $A$.

En vez de definir los funtores derivados a través de resoluciones
proyectivas e inyectivas se definen como la homología
(cohomología) de ciertos complejos de cadena (cocadena), basándose
en \cite[Appendix II.3]{calculus}, \cite[XII.5.5]{BK} y
\cite[XI.6.2]{BK}.

\subsection*{1.7. Teorema de Normalización}
En esta sección se enuncian los Teoremas de Normalización para
grupos abelianos simpliciales y cosimpliciales, los cuales se
pueden encontrar en \cite[22.1, 22.3]{may} y \cite[VIII.1]{goerss}
respectivamente. Estos teoremas establecen que la homología
(cohomología) de un grupo abeliano simplicial (cosimplicial) se
puede calcular omitiendo los símplices (cosímplices) degenerados.

\section*{Capítulo 2. Una sucesión espectral.}
En este capítulo se introducen un par de sucesiones espectrales
que tienen por límites $\limdir_i F$ y $\liminv^i F$, donde
$F:\C\rightarrow \Ab$ es un funtor y $\C$ es un poset graduado, y
se dan condiciones de convergencia (débil). La construcción se
hace a partir de módulos diferenciales graduados \cite{mccleary}.

Estas sucesiones espectrales son el pilar fundamental en los
desarrollos posteriores que se realizan en este trabajo. Las
páginas iniciales se calculan como homología (cohomología) de
subgrupos abelianos simpliciales (cosimpliciales) cuyos símplices
(cosímplices) son cadenas que empiezan y acaban con objetos de
$\C$ en los cuales la función de grado toma valores fijos.

\subsection*{2.1. Ejemplos.} En esta sección se describen numerosos
ejemplos para ilustrar la construcción de las sucesiones
espectrales en casos concretos. En particular se trabaja sobre las
siguientes categorías:
\begin{itemize}
\item pushout $ \xymatrix{a  & b \ar[l]_{f} \ar[r]^{g} &  c,}$
\item pullback $ \xymatrix{ & c \ar[d]^{g} \\ a \ar[r]_{f} & b,}$
\item telescopio $ \xymatrix{a_0 \ar[r]^{f_1} & a_1 \ar[r]^{f_2} &
a_2 \ar[r]^{f_3} &  a_3 \ar[r]^{f_4} & a_4 ...}$ y
\item ciclo $ \xymatrix{
a \ar[rd]^<<<<{g} \ar[r]^{f} & c \\
b \ar[ru]^<<{h} \ar[r]_{i} & d. \\
}
$
\end{itemize}
Algunos de los comportamientos particulares de las sucesiones
espectrales que se observan en estos ejemplos se generalizan en
capítulos posteriores a teoremas de ámbito general.

\section*{Capítulo 3. Límites directos superiores.}
En este capítulo se encuentran condiciones suficientes para que
los límites derivados superiores de un funtor dado
$F:\C\rightarrow \Ab$ se anulen, es decir, para que $\limdir_i
F=0$ si $i\geq 1$. A estos funtores se les llamará funtores
\emph{$\limdir$-acíclicos}. Así mismo se \linebreak desarrollan
ciertas herramientas para el cálculo de estos límites superiores
en el caso general, y se encuentran aplicaciones basándose en una
analogía con el concepto de flujo en un grafo.

El desarrollo es como sigue: en la Sección 3.1 se caracterizan los
objetos proyectivos de la categoría abeliana $\Ab^\P$, donde $\P$
es un poset graduado. Recordemos que los funtores proyectivos son
$\limdir$-acíclicos. A partir de esta caracterización se deduce
una condición más débil, \emph{pseudo-proyectividad}, la cual
también implica $\limdir$-aciclicidad. Este hecho se prueba en la
Sección 3.2 mediante el uso de las sucesiones espectrales del
Capítulo 2.

En la Sección 3.3 se construye un funtor pseudo-proyectivo
asociado a un funtor arbitrario $F:\P\rightarrow \Ab$. Mediante su
uso y un argumento clásico de ``cambio de dimensiones'' se
demuestra como reducir el cálculo de un límite superior $\limdir_i
F$ al cálculo de un límite superior de orden $1$ $\limdir_1 G$,
donde $G$ es un funtor que se construye a partir de $F$.

Por último, en la Sección 3.4, se aplican los resultados de la
sección anterior al caso en el que cierta subcategoría, denotada
\emph{core(\P)}, de la categoría $\P$ es un árbol. Todo el
capítulo está ilustrado con ejemplos en los se concretan las ideas
teóricas desarrolladas.

\subsection*{3.1. Objetos proyectivos de $\Ab^\P$.}
Dado un funtor $F:\P\rightarrow \Ab$, donde $\P$ es un poset
graduado, podemos considerar, para cada objeto $i\in \Ob(\P)$, el
cociente de $F(i)$ por las imágenes $F(\alpha)$ donde $\alpha\in
\Hom(\P)$ es un morfismo no trivial que termina en $i$. A este
cociente lo denotamos $\coker(i)$.

Pues bien, en esta sección se demuestra que si $F$ es proyectivo
en $\Ab^\P$ entonces:
\begin{itemize}
\item $\coker(i)$ es proyectivo en $\Ab$ para todo $i\in \Ob(\P)$.
\item $F$ verifica una condición técnica adicional que llamamos
\emph{pseudo-proyectividad} y la cual, en casos extremos, se
reduce a la inyectividad de los morfismos $F(\alpha)$, $\alpha\in
\Hom(\P)$.
\end{itemize}

Posteriormente se demuestra que si el poset graduado $\P$ es
acotado inferiormente entonces estas dos condiciones mencionadas
arriba implican la proyectividad de $F$, con lo cual se obtiene
una caracterización completa de los proyectivos de $\Ab^\P$.

\subsection*{3.2. Pseudo-proyectividad.}
Recuérdese que la proyectividad de un objeto implica que los
funtores derivados izquierdos se anulan sobre él, en particular
los objetos proyectivos de $\Ab^\P$ son $\limdir$-acíclicos.
Teniendo en cuenta que la pseudo-proyectividad es más débil que la
proyectividad en $\Ab^\P$ (ver Sección 3.1), es natural
preguntarse si esa condición implica por sí sola
$\limdir$-aciclicidad.

Esta sección se dedica exclusivamente a probar este hecho, que
la\linebreak pseudo-proyectividad implica $\limdir$-aciclicidad.
Para ello se hace un uso intensivo de las sucesiones espectrales
del Capítulo 2. Se asume la hipótesis adicional de que $\P$ es
acotado inferiormente, lo cual permite hallar una página de la
sucesión espectral $E_r$, para $r$ suficientemente grande, en la
cual probar que las contribuciones al límite son cero.

\subsection*{3.3. Calculando límites superiores.}
En esta sección se asocia a cada funtor $F:\P\rightarrow \Ab$ un
funtor pseudo-proyectivo $F'$ y una transformación natural
$F'\Rightarrow F$. Este funtor tiene varias propiedades
interesantes, entre las que destaca ser pseudo-proyectivo y por
tanto, gracias a la Sección 3.2, $\limdir$-acíclico.

La sucesión exacta corta
$$
0\Rightarrow K_F\Rightarrow F'\Rightarrow F\Rightarrow 0
$$
nos da $\limdir_i F=\limdir_{i-1} K_F$ para $i>1$. La iteración de
este argumento reduce el cálculo de todos los límites superiores
de $F$ al cálculo de límites de orden $1$ sobre sucesivos
funtores. Como aplicación se obtiene que un funtor (con morfismos
inyectivos) sobre un árbol es $\limdir$-acíclico.

\subsection*{3.4. $Lim_1$ como un problema de flujo.}
Los problemas de flujo son un tópico en la Teoría de Grafos
\cite{bela}. Un flujo no es más que la asignación de valores a las
aristas de un grafo de forma que en cada vértice el flujo entrante
iguala al saliente. Mientras que clásicamente estos valores son
números naturales, un sencillo argumento permite describir los
$1$-ciclos (del complejo de cadena cuya homología en grado $1$ es
$\limdir_1 F$) como ``flujos'' cuyo valor en la arista
$\alpha:i_0\rightarrow i_1$ es un elemento de $F(i_0)$.

Esto da inmediatamente que $\limdir_1 F=0$, es decir, todo
$1$-ciclo es $1$-borde, si y sólo si cada ``flujo'' se escribe
como suma de ciertos flujos triviales minimales. \linebreak
Posteriormente, la relativa rigidez de la estructura de un poset
graduado permite re\-escribir cada ``flujo''  como un ``flujo''
sobre cierta subcategoría de $\P$ denotada $core(\P)$.

Finalmente, se obtiene que la aciclicidad de $core(\P)$, es decir,
que esta subcategoría sea un árbol, junto con la inyectividad de
$F$ implican la $\limdir$-aciclicidad de este último.

\section*{Capítulo 4. Límites inversos superiores.}
Este capítulo es la dualización del Capítulo 3 y su desarrollo es
como sigue:
\begin{itemize}
\item \textbf{4.1} Se caracterizan los objetos inyectivos de
$\Ab^\P$. \item \textbf{4.2} Se prueba, mediante el uso de las
sucesiones espectrales del Capítulo 2, que la pseudo-inyectividad
implica $\liminv$-aciclicidad. \item \textbf{4.3} Mediante
sucesiones exactas cortas se reducen el cálculo de $\liminv^i F$ a
límites $\liminv_1 G$, donde $G$ es un funtor que se construye a
partir de $F$.
\end{itemize}

Como aplicación se han desarrollado herramientas para el cálculo
de cohomología de categorías y, debido a su extensión, esto
constituye el Capítulo 5

\section*{Capítulo 5. Cohomología}
En esta sección, y como aplicación de la teoría desarrollada en el
capítulo anterior, de desarrollan herramientas para el cálculo de
la cohomología $H^*(\P;\Z)$ de un poset graduado $\P$, bajo
ciertas hipótesis estructurales adicionales sobre $\P$.

Los conceptos introducidos a lo largo de las distintas secciones
son de naturaleza local, en el sentido de que dependen de las
subcategorías $(i\downarrow \P)$ y la restricciones
$F_{(i\downarrow \P)}$ donde $i$ es un objeto arbitrario de $\P$.

Inicialmente se asume que $F$ es un funtor que toma por valores
grupos abelianos libres, para especializarse posteriormente al
caso en que $F=c_\Z$ y de este modo calcular la cohomología de
$\P$. La última sección describe una familia de posets \linebreak
graduados, los \emph{simplex-like} posets, a la cual pertenecen
los complejos simpliciales, y cuyos miembros verifican las
hipótesis estructurales adicionales.

La idea que guía todo el capítulo es la siguiente: en el capítulo
anterior se \linebreak construyó una sucesión exacta corta
$$
0\Rightarrow F\Rightarrow F'\Rightarrow C_F\Rightarrow 0,
$$
donde $F'$ es un funtor $\liminv$-acíclico, con objeto de reducir
el cálculo de los limítes inversos superiores de $F$ a limítes de
orden inferior sobre otro funtor. Los funtores $F'$ y $C_F$ puede
ser ``muy grandes'' (en el sentido de cantidad de generadores), y
con objeto de disminuir este tamaño se considera otro funtor
$\ker_F$ ``más pequeño'', y se buscan condiciones bajo las cuales
podamos construir sucesiones exacta cortas del tipo
$$
0\Rightarrow F\Rightarrow \ker'_F\Rightarrow G\Rightarrow 0
$$
iterativamente.
\subsection*{5.1. Functores $p$-condensados.}
En esta sección se considera el funtor $\ker_F$ asociado al funtor
$F:\P\rightarrow \Ab$, y cuyo valor en $i\in\Ob(\P)$ es la
intersección de los morfismos no triviales que salen de $i$.
Fijado el entero $p$, diremos que $F$ es \emph{$p$-condensado} si
\begin{itemize}
\item $F(i)=0$ si $deg(i)<p$, y \item $\ker_F(i)=0$ si $deg(i)>p$.
\end{itemize}

Como veremos, esta condición significa que la información del
funtor $F$ está contenida de alguna manera en sus objetos de grado
$p$. Técnicamente nos permite \linebreak construir la sucesión
exacta corta deseada:
$$
0\Rightarrow F\Rightarrow \ker'_F\Rightarrow G\Rightarrow 0.
$$

La sección termina dando una caracterización de cuándo el funtor
$G$ de la sucesión exacta corta de arriba es $(p+1)$-condensado.
Esto es importante ya que el objetivo final del capítulo es hacer
esta construcción iterativamente. Además, esta caracterización se
usa en la siguiente sección para reducir el problema de cuándo es
$G$ $(p+1)$-condensado a un conjunto de ecuaciones enteras con
coeficientes que dependen de la estructura local de $\P$.

\subsection*{5.2. Familias recubridoras.}
Dado un poset graduado $\P$, una \emph{familia recubridora} $\J$
del mismo es una familia de subconjuntos $J^i_p\subseteq
(i\downarrow \P)_p$, donde $i$ recorre los objetos de $\P$ y $p$
recorre los grados $p\leq deg(i)$, sometida a ciertas
restricciones combinatorias. Como se ve la existencia o no de una
familia recubridora sólo depende de la estructura local de $\P$.
Esto permite (ver Sección 5.4) construir familias recubridoras
para posets que tiene una estructura local homogénea y adecuada,
como son los simplex-like posets.

Un funtor $F:\P\rightarrow \Ab$ que toma por valores grupos
abelianos libres, se dirá que es \emph{$\J$-determinado}, donde
$\J$ es una familia recubridora para $\P$, si la información de
$F$ está contenida en los objetos de los subconjuntos de $\J$. Así
se llega a la\\

\textbf{Proposición
\ref{lem_covering_familiy_determination_hereditary}. }\textit{Sea
$\P$ un poset graduado y $\J$ una familia recubridora para $\P$.
Si el funtor $F$ es $p$-condensado y $\J$-determinado entonces el
funtor $G$ de la sucesión exacta corta
$$
0\Rightarrow F\Rightarrow \ker'_F\Rightarrow G\Rightarrow 0
$$
es $(p+1)$-condensado y $\J$-determinado si se verifican un
conjunto de ecuaciones enteras.}

Este teorema se puede aplicar iterativamente sin más que comprobar
si ciertas ecuaciones enteras se verifican.

\subsubsection*{Familias recubridoras adecuadas.}
Esta subsección se centra en el caso $F=c_\Z$, con objeto de
calcular la cohomología $H^*(\P;\Z)$. Si $\J$ es una familia
recubridora para el poset graduado $\P$, diremos que es
\emph{adecuada} si podemos aplicar el Teorema 5.2.3
iterativamente, es decir, si las ecuaciones enteras involucradas
se verifican.

La condición de \emph{adecuación} de $\J$ se define, de nuevo,
mediante ciertas ecuaciones enteras que dependen de la estructura
local de $\P$. En el caso de que $\J$ sea adecuada se obtiene una
sucesión de funtores $F_0=c_\Z$, $F_1$, $F_2$,... cada uno de los
cuales encaja en una sucesión exacta corta
$$
0\Rightarrow F_p\Rightarrow \ker'_{F_p}\Rightarrow
F_{p+1}\Rightarrow 0.
$$

\subsubsection*{Bases y morfismos locales.}
En esta subsección se dan descripciones explícitas de bases para
los grupos abelianos libres que toman como valores los funtores
$F_0,F_1,F_2,...$ construidos en la subsección anterior. Así mismo
se dan descripciones explícitas de los morfismos $F_p(\alpha)$,
donde $\alpha\in \Hom(\P)$. Una vez más, estas descripciones
dependen de la estructura local de $\P$ (y de la familia $\J$
considerada).

\subsection*{5.3. Comportamiento global.}
Centrándose aún en el caso $F=c_\Z$, y considerando la sucesión de
funtores $F_0$, $F_1$, $F_2$,... descrita en la Subsección 5.2.1,
se describen en esta sección propiedades de los límites superiores
de $F$ que no dependen de la familia recubridora $\J$ elegida para
construir la mencionada sucesión de funtores (en la subsección
5.2.2 se describieron propiedades que sí dependen de $\J$).

En particular, se prueba un hecho análogo a que la cohomología en
grado $n$ de un $CW$-complejo sólo depende del $(n+1)$-esqueleto,
y una fórmula para la característica de Euler, que se reduce a la
conocida suma alternada del número de objetos de cada grado en
caso de que $\P$ sea un simplex-like poset (ver Sección 5.4).

\subsubsection*{Familias recubridoras globales.}
En esta subsección se encuentra el teorema principal de este
capítulo. Se comienza definiendo el concepto de \emph{familia
recubridora global $\K$}, que a diferencia de su versión local, es
decir, de familia recubridora, es una familia de subconjuntos
$K_p\subseteq \Ob_p(\P)$ donde $p$ recorre los grados posibles del
poset graduado $\P$. Obtenemos entonces el\\

\textbf{Teorema \ref{Thm_adequate_local y global characterization
acyclicity}. }\textit{Sea $\P$ un poset graduado para el que
existe una familia recubridora adecuada $\J$ y una familia
recubridora global adecuada $\K$. Entonces $\P$ es acíclico, es
decir, $H^i(\P;\Z)=0$ para $i>0$, si y sólo si $|K_0|$ iguala el
número de componentes conexas de $\P$.}

Como se ve, se ha reducido un problema de Álgebra Homológica, la
aciclicidad de $\P$, a una ecuación entera que involucra elementos
geométricos de $\P$.

\subsection*{5.4. Simplex-like posets.}
Un simplex-like poset no es más que un poset gra\-duado $\P$ cuyos
subcategorías locales $(i\downarrow \P^{op})$ son isomorfas al
poset de subdivisión del poset lineal $0<1<...<deg(i)$. Por
ejemplo, para $deg(i)=2$ la categoría $(i\downarrow \P^{op})$
tendría el aspecto:
$$
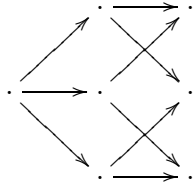
\xymatrix{ & .\ar[r]\ar[rd] & .\\
. \ar[ru]\ar[r]\ar[rd] & . \ar[ru]\ar[rd]&.\\
& .\ar[r]\ar[ru]&.}
$$

Si $\P$ es simplex-like entonces, como se prueba en esta sección,
existe una familia recubridora adecuada para $\P^{op}$, lo cual se
usará en el Capítulo 6. La sección acaba reenunciando los
resultados de la Sección 5.3 en el caso de simplex-like posets.

\section*{Capítulo 6. Aplicación: la conjetura de Webb.}
Denotemos por $\Ss_p(G)$ el complejo de Brown para el primo $p$,
que fue introducido en \cite{brown}. Webb conjeturó que el espacio
de órbitas $\Ss_p(G)/G$ es contráctil (como espacio topológico),
lo cual fue probado por Symonds en \cite{symonds}, extendido a
bloques por Barker \cite{barker1,barker2} y extendido a sistemas
de fusión (saturados) arbitrarios por Linckelmann \cite{markus}.

Los trabajos de Symonds y Linckelmann prueban la contractibilidad
del espacio de órbitas mostrando que es simplemente conexo y
acíclico, e invocando el Teorema de Whitehead. Ambas pruebas de
aciclicidad trabajan con el subposet de las cadenas normales.
Symonds usa los resultados de Thévenaz y Webb \cite{G-homot} sobre
que el subposet de las cadenas normales es $G$-equivalente al
complejo de Quillen. Linckelmann prueba por su cuenta que el
espacio de órbitas y el espacio de órbitas sobre las cadenas
normales tienen la misma cohomología entera \cite[Theorem
4.7]{markus}.

En este capítulo se aplican los resultados del Capítulo 5 para dar
una prueba alternativa de que el espacio de órbitas sobre las
cadenas normales es acíclico. Para ello, se considera un sistema
de fusión saturado $(S,\Ff)$ y el espacio de órbitas sobre las
cadenas normales correspondiente, denotado $[S_\lhd(\Ff)]$.
Fácilmente se comprueba que este poset es un simplex-like poset,
por lo cual (ver Sección 5.4) existe una familia recubridora
adecuada para $[S_\lhd(\Ff)]^{op}$.

Para construir una familia global adecuada $\K$ para
$[S_\lhd(\Ff)]^{op}$ se usa el mismo emparejamiento que usa
Linckelmann \cite[Definition 4.7]{markus}. Del hecho de que el
único subgrupo del $p$-grupo $S$ que iguala a su normalizador es
el propio $S$ se deduce que $|K_0|=1$. Esto, junto con que el
poset $[S_\lhd(\Ff)]^{op}$ es conexo, nos da, gracias al Teorema
$5.3.6$, la aciclicidad de este poset.

La prueba de que la familia $\K$ definida a través del
emparejamiento es una familia recubridora global adecuada es
bastante técnica, y se postpone a la Subsección 6.1.

\subsection*{6.1. $\K$ es una familia recubridora global adecuada.} El
subconjunto $K_n$ se define como aquellas clases de isomorfismo de
cadenas normales $[Q_0<...<Q_n]$ que tienen como normalizador al
propio $Q_n$ para cualquier representante $Q_0<...<Q_n$. Los
detalles técnicos usan profusamente los resultados de
\cite[Appendix]{blo2}, junto al hecho de que el sistema de fusión
bajo consideración $(S,\Ff)$ es saturado.

\section*{Capítulo 7. Aplicación: colímite homotópico.}
En este capítulo se estudia el problema de cuando el colímite
homotópico de un diagrama de espacios clasificadores de grupos
coincide con el espacio clasificador del colímite de los grupos.
Para ello, dado un diagrama de grupos $G:\P\rightarrow \Grp$ y un
cono $\tau:G\Rightarrow G_0$ se estudia la fibra $F$ de la
aplicación
$$
\hocolim BG\rightarrow BG_0.
$$
La relación con los desarrollos previos consiste en que la
homología de $F$ se calcula como los límites directos derivados
para cierto funtor $H:\P\rightarrow \Ab$. Este hecho es el teorema
principal de este capítulo, y su prueba se postpone a la Sección
7.1.

El resto del capítulo se dedica a ejemplos de aplicación del
teorema, como son el Teorema de Whitehead sobre el pushout y los
siguientes:
\begin{itemize}
\item Si $G_0$ es un grupo localmente finito entonces se tiene que
$$
\hocolim_{G\subseteq G_0, G finito} BG\simeq BG_0.
$$
\item Para cualquier grupo $G_0$ tenemos
$$
\hocolim_{G\subseteq G_0, \textit{G p-subgrupo finito normal}}
BG\simeq B(\lim_{G\subseteq G_0,\textit{G p-subgrupo finito
normal}} G).
$$
\end{itemize}

\subsection*{7.1. Demostración del teorema.} La demostración del teorema está basada
en la descripción de la fibra $F$ como el colímite homotópico de
las fibras sobre cada objeto \cite{cha}. También se usa la
sucesión espectral de Van Kampen \cite{stover}, la sucesión
espectral en homología de Bousfield-Kan \cite{BK} y la sucesión
larga de homotopía para fibraciones.

\subsection*{7.2. Otro ejemplo.} En esta subsección final se aplican el teorema central de este
capítulo junto con las herramientas desarrolladas para límites
directos derivados a un ejemplo concreto propuesto por A. Libman.
Se obtiene una fibración
$$\bigvee_{\alpha\in G_0\setminus\{1\}} (S^2)_\alpha\rightarrow \hocolim BG\rightarrow BG_0$$
donde el colímite homotópico se toma sobre cierto poset graduado
de dimensión $2$.

\chapter{Notation and Preliminaries}\label{chapter_preliminaries}

\section{Categories}\label{section categories} Throughout this work we use the following familiar categories:
\begin{itemize}
\item$\Set$, the category of sets, \item$\Grp$, the category of
groups, \item$\Ab$, the category of abelian groups, \item$\Top$,
the category of topological spaces with arrows continuous maps,
\item$\SSet$, the category of simplicial sets, \item$\Cat$, the
category with objects the small categories and with arrows the
functors between them,
\end{itemize}
and the pointed versions $\SSet_*$ and $\Top_*$. Any other
category is assumed to be a small category without explicit
mention. The objects of a category $\C$ are denoted by $\Ob(\C)$
and the arrows by $\Hom(\C)$. If $c,c'\in \Ob(\C)$ then
$\Hom_\C(c,c')$ denotes the arrows in $\C$ from $c$ to $c'$. Any
functor $F:\C\rightarrow \D$ is covariant if not stated otherwise.
If $s, s'\in \Ob(\Set)$ and $k\in s'$ then the constant function
from $s$ to $s'$ of value $k$ is denoted by $c_k$. Now we recall
some concepts in Category Theory (see \cite{MCL}):

\begin{Defi}\label{defi_connected}
A category $\C$ is \emph{connected} if for any two objects
$c,c'\in \Ob(\C)$ exists a chain of morphisms in $\C$ between $c$
and $c'$:
$$c\rightarrow c_1\leftarrow c_2\rightarrow ...c_{n-1}\rightarrow c_n\leftarrow c'$$
\end{Defi}

\begin{Defi}\label{defi_under_category}
Let $S:\D\rightarrow \C$ be a functor and $c_0\in \Ob(\C)$. The
category of \emph{objects $S$-under $c_0$}, \emph{$(c_0\downarrow
S$)}, has objects all the pairs $(f,d)$, where $d\in \Ob(\D)$ and
$f\in \Hom_\C(c_0,S(d))$, and arrows $h:(f,d)\rightarrow (f',d')$
those arrows $h:d\rightarrow d'$ in $\D$ for which $S(h)\circ
f=f'$.
\end{Defi}

\begin{Defi}\label{defi_over_category}
Let $S:\D\rightarrow \C$ be a category and $c_0\in \Ob(\C)$. The
category of \emph{objects $S$-over $c_0$}, \emph{$(S\downarrow
c_0)$}, has objects all the pairs $(f,d)$, where $d\in \Ob(\D)$
and $f\in \Hom_\C(S(d),c_0)$, and arrows $h:(f,d)\rightarrow
(f',d')$ those arrows $h:d\rightarrow d'$ such that $f'\circ
S(h)=f$.
\end{Defi}

In particular, for the identity functor $1_\C:\C\rightarrow \C$,
we  have the categories of \emph{objects under} and \emph{over
$c_0$}$\in \Ob(\C)$ respectively:
\begin{itemize}
\item $(c_0\downarrow \C)\definicio (c_0\downarrow 1_\C)$ and\item
$(\C \downarrow c_0)\definicio (1_\C\downarrow c_0)$.
\end{itemize}

\section{Groups and abelian groups}\label{section groups & graphs}
For the product of two elements $g,g'\in G\in\Grp$ the
multiplicative notation $gg'$ is used, while for $a,a'\in A\in\Ab$
the additive notation $a+a'$ is used. Another notations for the
category $\Ab$ are:
\begin{itemize}
\item $\sum_{i\in \I} A_i$ is the subgroup of $A\in \Ab $
generated by the subgroups $A_i\subseteq A$. \item
$\bigoplus_{i\in\I} A_i$ denotes the direct sum (coproduct) of the
abelian groups $\{A_i,\textit{ $i\in \I$}\}$, and $\prod_{i\in \I}
A_i$ denotes its direct product. \item $\Z$: infinite cyclic
group, $\Z_n$: order $n$ cyclic group, $\Q$: rational numbers as
additive group, $\Z[p^\infty]$, $p$ prime: the subgroup of $\Q/\Z$
generated by the cosets $1/p^n+\Z$ for $n\geq 0$. Recall that
$$\Q/\Z=\bigoplus_{\textit{$p$ prime}}\Z[p^\infty].$$ \item
$\Z\stackrel{\times n}\rightarrow \Z$ denotes the homomorphism
$m\mapsto nm$ and $\Z\stackrel{red_n}\rightarrow \Z/n$ is the
canonical projection, where $n\geq 0$. \item If $f:A\rightarrow B$
and $g:C\rightarrow D$ are homomorphisms then $f\oplus g:A\oplus C
\rightarrow B\oplus D$ is given by $a\oplus c\mapsto f(a)\oplus
g(c)$ and, in case $A=C$, $f\times g:A\rightarrow B\oplus
D=B\times D$ is given by $a\mapsto f(a)\oplus g(a)=(f(a),g(a))$.
Finally, in case $B=D$, $f+g:A\oplus C\rightarrow B$ is defined by
$a\oplus c\mapsto (f+g)(a\oplus c)=f(a)+g(c)$.
\end{itemize}

Consider a functor $F:\C\rightarrow \Grp$ with the category of
groups as target category. The functor $F$ is called
\label{defi_monic} \emph{monic} if $F(f)$ is a monomorphism for
each $f\in \Hom(\C)$. A natural transformation $\tau:F\Rightarrow
F'$ between functors $F,F':\C\rightarrow Grp$ is called
\emph{monic} if $\tau_i$ is a monomorphism for each $i\in
\Ob(\C)$.

Similarly we call a functor $F:\C\rightarrow \Grp$
\label{defi_epic} \emph{epic} if $F(f)$ is an epimorphism for each
$f\in \Hom(\C)$. A natural transformation $\tau:F\Rightarrow F'$
between functors $F,F':\C\rightarrow \Grp$ is called \emph{epic}
if $\tau_i$ is an epimorphism for each $i\in \Ob(\P)$. Notice
that, by the natural inclusion $\Ab\subseteq \Grp$, we have
defined the terms \emph{monic} and \emph{epic} also for functors
$\C\rightarrow \Ab$ and natural transformations between them.

We collect a few basic facts about the category $\Ab$ of abelian
groups:
\begin{enumerate}
\item The projective objects of $\Ab$ are the free abelian groups
\cite[III, Theorem 18]{griffith}. \item The injective objects of
$\Ab$ are direct sums of $\Q$ and $\Z[p^{\infty}]$ for various
primes $p$ \cite[III, Theorem 21]{griffith}.\item Subgroups of
free groups are free \cite[II, Theorem 15]{griffith}. \item
Finitely generated torsion free groups are free \cite[II, Theorem
16]{griffith}.
\end{enumerate}

We are also interested in maps $A\stackrel{f}\rightarrow B$
between free abelian groups which have free cokernel. 

\begin{Defi}\label{defi_pure_map}
Let $A\stackrel{f}\rightarrow B$ be a map between free abelian
groups. We say that $f$ is \emph{pure} if $\coker(f)$ is a free
abelian group.
\end{Defi}

If $A\cong \Z^n$ is a finitely generated free abelian group we
call $\dim(A)\definicio n$. We have the following property of pure
maps, which will be used repeatedly in successive sections,

\begin{Lem}\label{lema_fpuremono_imply_iso}
Let $A\stackrel{f}\rightarrow B$ be a map in $\Ab$ between free
abelian groups of the same rank. If $f$ is pure and injective then
it is an isomorphism.
\end{Lem}
\begin{proof}
The short exact sequence of free abelian groups
$$
0\rightarrow A\stackrel{f}\rightarrow B\rightarrow
\coker(f)\rightarrow 0
$$
implies that
$$
\dim A-\dim B+\dim(\coker(f))=0
$$
and thus, $$\dim(\coker(f))=0$$ and $\coker(f)=0$.
\end{proof}

\begin{Defi}\label{defi_free_functor}
Let $F:\C\rightarrow \Ab$ be a functor. We say that $F$ is
\emph{free} if $F(i)$ is a free abelian group for each $i\in
\Ob(\C)$.
\end{Defi}

Notice that this does not imply that $F$ is a free object in
$\Ab^{\C}$. For example, consider the category
$\C=\cdot\rightarrow \cdot$ and the functor $F\in \Ab^{\C}$ with
values
$$
\xymatrix{\Z\ar[r]^{0}&\Z.
 }
$$
The functor $F$ is free (Definition \ref{defi_free_functor}).
However, it is not projective by Corollary
\ref{cor_projective_bounded_graded}, and thus $F$ is not a free
object in $\Ab^{\C}$.

\section{Graded posets}\label{section graded_posets} In this section we
define a special kind of categories, graded posets, which are the
main ingredient in most of the results of this work.

We begin defining posets:

\begin{Defi}
A \emph{poset} is a category $\P$ in which, given objects $p$ and
$p'$,
\begin{itemize}
\item there is at most one arrow $p\rightarrow p'$, and \item if
there are arrows $p\rightarrow p'$ and $p'\rightarrow p$ then
$p=p'$.
\end{itemize}
\end{Defi}

In any poset $\P$ define a binary relation on its objects
$\Ob(\P)$ with $p\leq p'$ if and only if there is an arrow
$p\rightarrow p'$. Then $\leq$ is reflexive, symmetric and
transitive, i.e., $(\Ob(\P),\leq)$ is a partial order. Conversely
any partial order determines a poset in which the arrows are
exactly those ordered pairs $(p,p')$ for which $p\leq p'$.

It is worthwhile noticing that if $\P$ is a poset and $p_0\in \P$
then the categories $(p_0\downarrow\P)$ and $(\P\downarrow p_0)$
defined in Section \ref{section categories} are exactly the full
subcategories of $\P$ with objects $\{p| \exists p_0\rightarrow
p\}$ and $\{p| \exists p\rightarrow p_0\}$ respectively. We define
also the categories $(p_0\downarrow\P)_*$ and $(\P\downarrow
p_0)_*$ as the full subcategories of $\P$ with objects $\{p|
\exists p_0\rightarrow p,p\neq p_0\}$ and $\{p| \exists
p\rightarrow p_0,p\neq p_0\}$ respectively.

\begin{Defi}
If $\P$ is a poset and $p<p'$ then $p$ \emph{precedes} $p'$ if
$p\leq p''\leq p'$ implies that $p=p''$ or $p'=p''$.
\end{Defi}
Most of the results of this work are about the following types of
posets:
\begin{Defi}
Let $\P$ be a poset. $\P$ is called \emph{graded} if there is a
function $deg:\Ob(\P)\rightarrow \Z$, called the \emph{degree
function} of $\P$, which is order preserving and that satisfies
that if $p$ precedes $p'$ then $deg(p')=deg(p)+1$. If $p\in
\Ob(\P)$ then $deg(p)$ is called the \emph{degree} of $p$.
\end{Defi}

Notice that the degree function associated to a graded poset is
not unique (consider the translations $deg'=deg+c_k$ for $k\in
\Z$). According to the definition the degree function increases in
the direction of the arrows: we say that this degree function is
\emph{increasing}. If the degree function is order reversing and
satisfies the alternative condition that $p$ precedes $p'$ implies
$deg(p')=deg(p)-1$, i.e., $deg$ decreases in the direction of the
arrows, then we say that $deg$ is a \emph{decreasing} degree
function. Clearly both definitions are equivalent (by taking
$deg'=-deg$). A poset which is graded satisfies some structural
conditions:

\begin{Lem}\label{lem_graded_poset_structural}
If $\P$ is a graded poset and $p<p'$ then there is an integer $n$
and a finite chain
$$
p=p_0< p_1< p_2<...< p_{n-1}< p_n=p'
$$
where $p_i$ precedes $p_{i+1}$ for $i=0,..,n-1$. Moreover, if
$$
p=q_0< q_1< q_2<...< q_{m-1}< q_m=p'
$$
is another finite chain where $q_i$ precedes $q_{i+1}$ for
$i=0,..,m-1$ then $m=n$.
\end{Lem}

The proof of this lemma is straightforward.

\begin{Ex}
The ``pushout category" $b\leftarrow a\rightarrow c$, the
``telescope category" $a\rightarrow b\rightarrow c\rightarrow...$,
and the opposite ``telescope category" $...\rightarrow
c\rightarrow b\rightarrow a$ are graded posets. The integers $\Z$
is a graded poset. The rationals $\Q$ with the usual order is a
poset but it is not a graded poset by the first condition of Lemma
\ref{lem_graded_poset_structural}.
\end{Ex}

As in the next example, when drawing a poset, we picture just the
arrows $p\rightarrow p'$ where $p$ precedes $p'$.

\begin{Ex}
The poset
$$
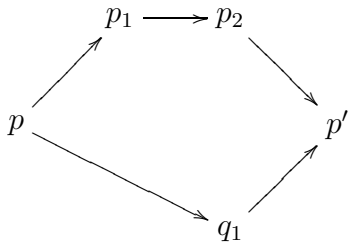
\xymatrix{ &p_1\ar[r] & p_2 \ar[rd]\\
p\ar[ru]\ar[rrd]&&& p'\\
&&q_1\ar[ru]}
$$
is not graded by the second condition of Lemma
\ref{lem_graded_poset_structural}.
\end{Ex}

If $\P$ is a graded poset we can ``extend" the degree function
$deg$ to the morphisms set $\Hom(\P)$ by $deg(p\rightarrow
p')=|deg(p')-deg(p)|$. By the preceding lemma this number does not
depend on the degree function. Whenever $\P$ is a graded poset we
denote by $\Ob_n(\P)$ the objects of degree $n$ and by
$\Hom_n(\P)$ the arrows of degree $n$ of the graded poset $\P$.
More generally:
\begin{Defi}\label{graded_posets_P_S}
Let $S\subset \Z$ and let $\P$ be a graded poset with degree
function $deg$. Then \emph{$\P_S$} is the full subcategory of $\P$
with objects $p$ such that $deg(p)\in S$, $\Ob_S(\P)=\{p\in
\Ob(\P)|deg(p)\in S\}$ and \emph{$\Hom_S(\P)$} is the set $\{f\in
\Hom(\P)|deg(f)\in S\}$.
\end{Defi}

\subsection{Boundedness on graded posets.}\label{graded_poset_boundedness}
For some results we restrict to:

\begin{Defi}
A graded poset $\P$ with increasing degree function $deg$ is
\emph{bounded below} (\emph{bounded above}) if the set
$deg(\P)\subset \Z$ has a lower bound (an upper bound).
\end{Defi}

If the degree function  $deg$ of $\P$ is decreasing then $\P$ is
\emph{bounded below} (\emph{bounded above}) if and only if
$deg(\P)\subset \Z$ has an upper bound (a lower bound). Notice
that upper boundedness means that the category $\P$ ends if you
move in the direction of the arrows, and lower boundedness means
that $\P$ finishes moving in the opposite direction.

If $\P$ is bounded below and over then $N\definicio
max(deg(\P))-min(deg(\P))$ exists and is finite. By Lemma
\ref{lem_graded_poset_structural} this number does not depend on
the degree function. We call it the dimension of $\P$, and we say
that $\P$ is $N$ dimensional.

\begin{Ex}
The ``pushout category" $b\leftarrow a\rightarrow c$ is
$1$-dimensional, the ``telescope category" $a\rightarrow
b\rightarrow c\rightarrow ..$ is bounded below but it is not
bounded over. The opposite ``telescope category" $..\rightarrow
c\rightarrow b\rightarrow a$ is bounded over but it is not bounded
below.
\end{Ex}

Notice that in a bounded above (below) graded poset there are
maximal (minimal) elements, but that the existence of maximal
(minimal) objects does not guarantee boundedness.

\begin{Rmk}\label{graded_posets_finiteness assumption on starting arrows}
When dealing with cohomology some assumptions on a bounded above
graded poset $\P$ shall be done. In particular, we shall assume
that $(i\downarrow \P)_n=(i\downarrow \P)_{\{n\}}$ (Definition
\ref{graded_posets_P_S}) is finite for each $i\in \Ob(\P)$, and
that all the maximal elements of $\P$ have the same degree.
\end{Rmk}

\subsection{Graphs and graded posets}\label{graded_poset_graphs}

Let $\P$ be a graded poset. The \emph{undirected graph associated}
to $\P$ has vertices the objects $\Ob(\P)$ and edges the arrows of
degree $1$, $\Hom_1(\P)$. The \emph{directed graph associated} to
$\P$ a has vertices the objects $\Ob(\P)$ and oriented edges the
oriented arrows of degree $1$, $\Hom_1(\P)$. \label{defi_tree}

A graded poset $\P$ is a \emph{tree} if it is associated
undirected graph is a tree (it contains no cycle). A tree or
maximal tree of $\P$ is a subcategory $\P'$ such that the
undirected subgraph associated to the graded poset $\P'$ is a tree
or maximal tree respectively of the undirected graph associated to
$\P$.

\section{Simplicial complexes}\label{section simplicial complexes}
A simplicial complex \cite{dwyer} is a pair $K=(V,S)$ where $V$ is
a set and $S$ is a collection of finite subsets of $V$ satisfying
\begin{enumerate}
\item $\sigma\in S$, $\sigma'\subseteq \sigma\Rightarrow
\sigma'\in S$. \item $v\in V\Rightarrow \{v\}\in S$.
\end{enumerate}

Elements of $V$ are called \emph{vertices} and elements of $S$ are
called \emph{simplices}. If $K=(V,S)$ is a simplicial complex we
can associate to it the poset with objects $S$ and inclusion as
order relation. This poset verify the following property:

\begin{Defi}\label{simplex-like poset}
Let $\P$ be a poset. Then $\P$ is \emph{simplex-like} if for all
$p\in \Ob(\P)$ the category $(\P\downarrow p)$ is isomorphic to
the poset of all non-empty subsets of a finite set (with inclusion
as order relation).
\end{Defi}

In fact, \cite[3.1]{garrett}, we have that a poset $\P$ arises
from a simplicial complex as above if and only if it is a
simplex-like poset and any two elements of $\P$ which have a lower
bound have an infimum, i.e., a greatest lower bound. If the poset
$\P$ arises from the simplicial complex $K=(V,S)$ then there is a
map $dim:\P\rightarrow \Z$ which assigns to each simplex $s\in S$
its dimension $dim(s)\in \Z$. If $s$ is a simplex of $S$ then all
the subsets of $s$ are in $S$ too. This implies that preceding
simplices of $\P$ differ just in one dimension, and thus the
function $dim$ is an increasing degree function and $\P$ is a
graded poset.


\section{Realization and homotopy type}\label{realization}
According to Sections \ref{section graded_posets} and \ref{section
simplicial complexes} we have inclusions
$$
\textit{Simplicial complexes $\subseteq$ Graded posets $\subseteq$
Categories.}
$$
On the one hand, we can realize in $\Top$ a simplicial complex $K$
as a space of formal linear combinations appropriately topologized
\cite[3]{dwyer}. On the other hand, for a category $\C$ we can
consider its \emph{nerve} $N\C$, which is a simplicial set $N\C\in
\SSet$, and the realization in $\Top$ of this simplicial set
\cite[3]{dwyer}. We denote by $|K|$ and $|N\C|$ the realizations
of simplicial complexes and categories respectively.

For a given simplicial complex $K$ consider the graded poset $\P$
associated to it (Section \ref{section simplicial complexes}). The
simplicial complex whose simplices consist of all of the totally
ordered subsets of $\P$ is exactly the barycentric subdivision
\cite[3.3]{spanier} of $K$, $sd K$. By \cite[3,4]{dwyer} there are
homeomorphisms $|K|\cong |sd K|\cong |N\P|$. Thus, as the
realization of a simplicial complex $K$ we can consider either of
$|K|$, $|sd K|$ or $|N\P|$.

\subsection{Homotopy type.}

From the homotopy viewpoint restricting to graded posets means no
loss:
\begin{Prop}
For any space $X\in \Top$ there is a graded poset $\P$ and a weak
homotopy equivalence:
$$
X\simeq |N\P|.
$$
\end{Prop}
\begin{proof}
It is a well known fact \cite[Theorem V.3.2]{whitehead} that $X$
has the weak homotopy type of a $CW$-complex. In fact, by
\cite[Theorem 2C.5]{hatcher}, $X$ also has the weak homotopy type
of a simplicial complex $K$, $X\simeq |K|$. By the comments above
we have $|K|\cong |sd K|\cong |N\P|$ for certain graded poset
$\P$, and thus $X\simeq|N\P|$.
\end{proof}
\begin{Rmk}
By the proof of the Theorem if $X$ is a $CW$-complex then there is
a (strong) homotopy equivalence $X\simeq |N\P|$.
\end{Rmk}

\section{Derived functors of direct and inverse
limit}\label{derivedfunctors} In this section we give definitions
for the left derived functors $\limdir_i$ of the direct limit
$\limdir:\Ab^\C\rightarrow \C$ and for the right derived functors
$\liminv^i$ of the inverse limit $\liminv:\Ab^\C\rightarrow \C$
for any small category $\C$.

Notice \cite[IX, Proposition 3.1]{maclanehom} that $\Ab^\C$ is an
abelian category in which the short exact sequences are the
object-wise ones, and that $\limdir$ and $\liminv$ are right exact
and left exact respectively because they are left adjoint and
right adjoint respectively to the functor $\Ab\rightarrow \Ab^\C$
which maps $A\mapsto c_A$.

It is well known (\cite[XI.6.1]{BK} and its dual) that in $\Ab^\C$
there are enough projectives and injectives so we can define the
derived functors of $\limdir$ and $\liminv$. Instead of
considering projective and injective resolutions for objects of
$\Ab^\C$, the definitions of $\limdir_i$ and $\liminv^i$ we
expound here have computational purposes and are based on
\cite[Appendix II.3]{calculus}. They also appear in
\cite[XII.5.5]{BK}, \cite[XI.6.2]{BK}. In \cite[p.409ff.]{goerss}
there is a summary.

Denote by $N\C$ the nerve of the small category $\C$ and by
$\sigma\in N\C_n$ an $n$-simplex of the nerve,
 that is, a chain of morphisms
$\sigma=\xymatrix{ \sigma_0 \ar[r]^{\alpha_1} & \sigma_1
\ar[r]^{\alpha_2}&...\ar[r]^{\alpha_{n-1}}&\sigma_{n-1}\ar[r]^{\alpha_n}
& \sigma_n}$. Given a covariant functor $F:\C \rightarrow \Ab$
consider the simplicial abelian group with simplices

$$C_n(\C,F)=\bigoplus_{\sigma \in {N\C}_n} F_\sigma,$$
where $F_{\sigma}=F(\sigma_o)$. The face map $d_i$ for $0\leq
i\leq n$ is the unique homomorphism which makes commute the
diagram
$$ \xymatrix{
C_n(\C,F) \ar[d]^{\pi_{\sigma}} \ar[r]^{d_i} & C_{n-1}(\C,F) \ar[d]^{\pi_{d_i(\sigma)}} \\
F_\sigma \ar[r]^{id^*} & F_{d_i(\sigma)} }$$ for every $\sigma \in
{N\C}_n$, where \begin{numcases}{id^*=}
F(\alpha_1):F(\sigma_0)\rightarrow F(\sigma_1) & for $i=0$ \nonumber \\
id_{F(\sigma_0)}:F(\sigma_0)\rightarrow F(\sigma_0) & for $0<i\leq
n.$ \nonumber
\end{numcases}

The degeneracy map $s_i$ for $0\leq
i\leq n$ is the unique homomorphism which makes commute the
diagram
$$ \xymatrix{
C_n(\C,F) \ar[d]^{\pi_{\sigma}} \ar[r]^{s_i} & C_{n+1}(\C,F) \ar[d]^{\pi_{s_i(\sigma)}} \\
F_\sigma \ar[r]^{id_{F(\sigma_0)}} & F_{s_i(\sigma)} }$$ for every
$\sigma \in {N\C}_n$.

This simplicial object gives rise to a chain complex, the Moore
complex, $(C_*(\C,F),d)$ with differential of degree $-1$,
$d:C_n(\C,F)\rightarrow C_{n-1}(\C,F)$, $d=\sum_{i=0}^n (-1)^i
d_i$.
\begin{Defi}\label{Defi:limi}
Let $\C$ be a small category and let $F$ be a covariant functor
\linebreak $F:\C\rightarrow \Ab$, then the $i$-left derived
functor of $\limdir:\Ab^\C\rightarrow \C$ is
$${\limdir}_i(F):= H_i(C_*(\C,F),d).$$
\end{Defi}

For the inverse limit $\liminv:\Ab^\C\rightarrow \C$ consider the
cosimplicial abelian group with simplices:
$$    C^n(\C,F)=\prod_{\sigma \in {N\C}_n} F^\sigma,$$
where $F^{\sigma}=F(\sigma_n)$. The coface map $d^i$ for $0\leq
i\leq n+1$ is the unique homomorphism which makes commute the
diagram
$$ \xymatrix{
C^n(\C,F) \ar[d]^{\pi_{d_i(\sigma)}} \ar[r]^{d^i} & C^{n+1}(\C,F) \ar[d]^{\pi_{\sigma}} \\
F^{d_i(\sigma)} \ar[r]^{id^*} & F^\sigma }$$ for every $\sigma \in
{N\C}_{n+1}$, where \begin{numcases}{id^*=}
F(\alpha_{n+1}):F(\sigma_n)\rightarrow F(\sigma_{n+1}) & for $i=n+1$ \nonumber \\
id_{F(\sigma_{n+1})}:F(\sigma_{n+1})\rightarrow F(\sigma_{n+1}) & for $0\leq i\leq n$ \nonumber
\end{numcases}

The codegeneracy map $s^i$ for $0\leq i\leq n$ is the unique homomorphism which makes commute the
diagram
$$ \xymatrix{
C^{n+1}(\C,F) \ar[d]^{\pi_{s_i(\sigma)}} \ar[r]^{s^i} & C^n(\C,F) \ar[d]^{\pi_\sigma} \\
F^{s_i(\sigma)} \ar[r]^{id_{F(\sigma_n)}} & F^\sigma }$$ for every
$\sigma \in {N\C}_n$.

This cosimplicial object gives rise to a cochain complex
$(C^*(\C,F),d)$ with differential of degree $1$,
$d:C^n(\C,F)\rightarrow C^{n+1}(\C,F)$, $d=\sum_{i=0}^{n+1} (-1)^i
d_i$.
\begin{Defi}\label{Defi:liminvi}
Let $\C$ be a small category and let $F$ be a covariant functor
\linebreak $F:\C\rightarrow \Ab$, then the $i$-right derived
functor of $\liminv:\Ab^\C\rightarrow \C$ is
$${\liminv}^i(F):= H^i(C^*(\C,F),d).$$
\end{Defi}

For every short exact sequence of natural trasformations
$$0\Rightarrow F\Rightarrow G\Rightarrow H\Rightarrow 0$$ in
$\Ab^\C$ there exists a pair of long exact sequences of derived
functors
$$.. \rightarrow {\limdir}_1 F \rightarrow {\limdir}_1 G \rightarrow
{\limdir}_1 H \rightarrow {\limdir} F \rightarrow {\limdir} G
\rightarrow {\limdir} H\rightarrow 0,$$
$$.. \leftarrow {\liminv}^1 H \leftarrow {\liminv}^1 G \leftarrow
{\liminv}^1  F \leftarrow {\liminv} H \leftarrow {\liminv} G
\leftarrow {\liminv} F\leftarrow 0.$$

We use the following notation for the obvious inclusions and projections:

$$F_\sigma \stackrel{i_\sigma}\hookrightarrow C_n(\C,F) \stackrel{\pi_\sigma}\twoheadrightarrow  F_\sigma,$$
$$F^\sigma \stackrel{i^\sigma}\hookrightarrow C^n(\C,F) \stackrel{\pi^\sigma}\twoheadrightarrow  F^\sigma.$$

\section{Normalization Theorem}\label{normalization}
We shall use the \emph{Normalization Theorem} for simplicial
abelian groups in order to compute the homology of these. It
states roughly that the homology of a simplicial abelian group can
be computed removing the degenerate simplices. The theorem we
state below is contained in \cite[22.1, 22.3]{may}. It can be also
found in \cite[III.2.1, III.2.4]{goerss}.

Let $S$ be a simplicial abelian group and let $(S,\sum (-1)^i
d_i)$ be the Moore chain complex associated to $S$. Consider the
chain complex $NS$ with $n$-chains
$$NS_n=\bigcap_{i=0}^{n-1} \ker(d_i:S_n\rightarrow S_{n-1})$$ and
with differential $(-1)^n d_n: NS_n\rightarrow NS_{n-1}$. Define
$DS$ as the chain subcomplex $DS$ of $(S,\sum (-1)^i d_i)$ which
$n$-chains are generated by the degenerate elements of $S$, that
is, by the elements in the image of some $s_i$.

\begin{Thm}\label{simp_normalization_thm}
Let $S$ be a simplicial abelian group, then:
$$H_*(S)=H_*(NS)=H_*(S/DS).$$
\end{Thm}

The dual version for cosimplicial abelian groups appears in
\cite[VIII.1]{goerss} and \linebreak \cite[X.7.1]{BK}. Let $C$ be
a cosimplicial abelian group and consider the cochain complex
$(C,\sum (-1)^i d^i)$. There is a cochain complex $NC$ with
$n$-cochains
$$NC^n=\bigcap_{i=0}^{n-1} \ker(s^i: C^n\rightarrow C^{n-1})$$ and
with differential $\sum (-1)^i d^i$. Define $DC$ as the subcomplex
of $(C,d=\sum (-1)^i d^i)$ which $n$-cochains are generated by the
elements in the image of some $d^i$.

\begin{Thm}\label{cosimp_normalization_thm}
Let $C$ be a cosimplicial abelian group, then:
$$H^*(C)=H^*(NC)=H^*(C/DC).$$
\end{Thm}

\begin{Rmk}\label{simp_cosimp_normalization}
We shall apply these theorems to simplicial and cosimplicial
abelian groups coming from a diagram as in Section
\ref{derivedfunctors}. In fact, we shall use the quotient chain
complex $S/DS$ when dealing with $C_n(\C,F)$ and the cochain
subcomplex $NC$ when dealing with $C^n(\C,F)$ without explicit
mention.
\end{Rmk}


\chapter{A spectral sequence}\label{spectral} In this
section we shall construct spectral sequences with targets
$\limdir_i F$ and $\liminv^i F$ for $F:\C\rightarrow \Ab$ with
$\C$ a graded poset. Some conditions for (weak) convergence shall
be given. We build the spectral sequences starting from filtered
differential modules (see \cite{mccleary}, where the notion of
weak convergence we use is also given).

First consider the complex $(C_*(\C,F),d)$ defined in Section
\ref{derivedfunctors} and choose a decreasing degree function
$deg$ over the objects of $\C$. There is a decreasing filtration
of this complex given by

$$L^pC_n(\C,F)=\bigoplus_{\sigma \in {N\C}_n, deg(\sigma_n)\geq
p} F_{\sigma}.$$

It is straightforward that the triple $(C_*(\C,F),d,L^*)$ is a
filtered differential graded $\Z$-module, so it yields a spectral
sequence $(E^{*,*}_r,d_r)$ of cohomological type whose
differential $d_r$ has bidegree $(r,1-r)$. The $E^{*,*}_1$ page is
given by
$$
    E^{p,q}_1\simeq H^{p+q}(L^pC/L^{p-1}C)
\textit{.}$$ The differential graded $\Z$-module $L^pC/L^{p-1}C$
is in fact a simplicial abelian group because the face operators
$d_i$ and the degeneracy operators $s_i$ respect the filtration
$L^*$. The $n$-simplices are
$$
(L^pC/L^{p-1}C)_n=\bigoplus_{\sigma \in {N\C}_n, deg(\sigma_n)=p}
F_{\sigma}.
$$
Moreover, for each $p$, $L^pC/L^{p-1}C$ can be filtered again by
the condition $deg(\sigma_0)\leq p'$ to obtain a homological type
spectral sequence. Then arguing as above we obtain:
\begin{Prop}\label{spectral_sequence_1}
For a (decreasing) graded poset and a functor $F:\C\rightarrow
\Ab$:
\begin{itemize}
\item There exists a cohomological type spectral sequence
$E^{*,*}_*$ with target $\limdir_*(F)$. \item There exists a
homological type spectral sequence $(E^p)_{*,*}^*$ with target the
column $E^{p,*}_1$ for each $p$.
\end{itemize}
\end{Prop}

Notice that the column $E^{p,*}_1$ is given by the cohomology of
the simplicial abelian group formed by the simplices that end on
objects of degree $p$, the column $(E^p)_{p',*}^1$ is given by the
homology of the simplicial abelian group formed by the simplices
that end on degree $p$ and begin on degree $p'$, and all the
differentials in the spectral sequences above are induced by the
completely described differential of $(C_*(\C,F),d)$. An advantage
of handling simplical abelian groups instead of chain complexes is
the chance to use the Normalization Theorem
\ref{simp_normalization_thm}.

As $\bigcup_p L^pC_n = C_n$ and $\bigcap_p L^pC_n = {0}$ for each
$n$ the spectral sequence $E^{*,*}_*$ converges weakly to its
target. In case the map $deg$ has a bounded image, i.e., when $\C$
is $N$ dimensional, the filtration $L^*$ is bounded below and
over, and so $E^{*,*}_*$ collapses after a finite number of pages.
The same assertions on weak converge and boundedness hold for the
spectral sequences $(E^p)_{*,*}^*$.

If we proceed in reverse order, i.e., filtrating first by the
degree of the beginning object and later by the degree of the
ending object, we obtain:

\begin{Prop}\label{spectral_sequence_2}
For a (decreasing) graded poset and a functor $F:\C\rightarrow
\Ab$:
\begin{itemize}
\item There exists a homological type spectral sequence
$E_{*,*}^*$ with target $\limdir_*(F)$. \item There exists a
cohomological type spectral sequence $(E_p)^{*,*}_*$ with target
the column $E_{p,*}^1$ for each $p$.
\end{itemize}
\end{Prop}

If the degree function we take is increasing then the appropriate
conditions for the filtrations are $deg(\sigma_n)\leq p$ and
$deg(\sigma_0)\geq p'$, and the spectral sequences obtained in
Propositions \ref{spectral_sequence_1} and
\ref{spectral_sequence_2} are of homological (cohomological) type
instead of cohomological (homological) type.

For the case of the cochain complex $(C^*(\C,F),d)$ defined in
Section \ref{derivedfunctors} the choices for the filtrations are
$deg(\sigma_n)\leq p$ and $deg(\sigma_0)\geq p'$ for a decreasing
degree function and $deg(\sigma_n)\geq p$ and $deg(\sigma_0)\leq
p'$ for an increasing one. Analogously we obtain spectral
sequences with target $\liminv^i F$ which columns in the first
page are computed by another spectral sequence. In this case we
can use the Normalization Theorem \ref{cosimp_normalization_thm}
to compute cohomology of the cosimplicial abelian groups appearing
in the page $1$ of these spectral sequences.

Table \ref{tabla_ss} shows a summary of the types of the spectral
sequences for all the cases. The statements on weak convergence
and boundedness apply to any of the spectral sequences of the
table.

\begin{center}
\begin{table}
\begin{tabular}{|c|c|c|c|c|c|}
\hline
Complex & Degree    & First       & Second     & First ss. & Second ss. \\
        & function  & filtration  & filtration &           &            \\
\hline
$C_*(\C,F)$ & decreasing    & $deg(\sigma_n)\geq$ &  $deg(\sigma_0)\leq$ & cohomol. type &  homol. type \\
$C_*(\C,F)$ & decreasing    & $deg(\sigma_0)\leq$ &  $deg(\sigma_n)\geq$ & homol. type   &  cohomol. type \\
$C_*(\C,F)$ & increasing    & $deg(\sigma_n)\leq$ &  $deg(\sigma_0)\geq$ & homol. type   &  cohomol. type \\
$C_*(\C,F)$ & increasing    & $deg(\sigma_0)\geq$ &  $deg(\sigma_n)\leq$ & cohomol. type &  homol. type \\
$C^*(\C,F)$ & decreasing    & $deg(\sigma_n)\leq$ &  $deg(\sigma_0)\geq$ & homol. type   &  cohomol. type \\
$C^*(\C,F)$ & decreasing    & $deg(\sigma_0)\geq$ &  $deg(\sigma_n)\leq$ & cohomol. type &  homol. type \\
$C^*(\C,F)$ & increasing    & $deg(\sigma_n)\geq$ &  $deg(\sigma_0)\leq$ & cohomol. type &  homol. type \\
$C^*(\C,F)$ & increasing    & $deg(\sigma_0)\leq$ &  $deg(\sigma_n)\geq$ & homol. type   &  cohomol. type \\
\hline
\end{tabular}
\caption{Filtrations and spectral sequences
obtained\label{tabla_ss}}
\end{table}
\end{center}

\begin{Rmk}\label{normalized_ss}
Recall from Remark \ref{simp_cosimp_normalization} the chain
(cochain) complex chosen to normalize a simplicial (cosimplicial)
abelian group. It is straightforward that normalizing the
simplicial (cosimplicial) abelian groups that computes the page
$1$ of the spectral sequences above has the same effect as
considering the spectral sequences of the normalizations of
$C_*(\C,F)$ ($C^*(\C,F)$).
\end{Rmk}

\section{Examples}\label{examples_ss}
Next there are some examples that show how the spectral sequences
just built work. In this section $F$ denotes a covariant functor
$F:\P\rightarrow \Ab$ where $\P$ is a graded poset.

These examples also serve as preamble to Chapters \ref{section_on
direct limit} and \ref{section_on inverse limit}, where conditions
are found such that $\limdir_i F=0$ and $\liminv^i F=0$ for $i>0$
respectively. The behaviour of the spectral sequences in these
examples resembles the general results in this thesis (Theorems
\ref{pro_acyclic_graded} and \ref{pro_inv_acyclic_graded}).

In particular, in the pushout example $F$ monic implies $\limdir_1
F=0$ (compare with Definition \ref{property_pseudo-projective},
Remark \ref{pseudo_injective_is injective} and Theorem
\ref{pro_acyclic_graded}). For the pullback example we have that
$F$ epic implies $\liminv^1 F=0$ (compare with Definition
\ref{property_pseudo-injective}, Remark
\ref{pseudo_surjective_is_surjective} and Theorem
\ref{pro_inv_acyclic_graded}). The telescope example shows the
importance that in a graded poset every morphism factors as
composition of morphisms of degree $1$ (compare with Step $1$ in
the proof of Theorem \ref{pro_acyclic_graded}). Finally, the
``cycle" example shows that the existence of cycles in $\P$ may
prevent $\limdir_1 F=0$ (compare with Corollaries
\ref{rootedtreemonic_acyclic} and \ref{cor_core_determines_limi}).
In this example also appears a heuristic version of
pseudo-projectiveness (Definition
\ref{property_pseudo-projective}). This property is related
(Theorem \ref{pro_acyclic_graded}) to the vanishing of $\limdir_i
F$.

\begin{Ex}
\textbf{Pushout:} Consider $\P$ the ``pushout category":
$$ \xymatrix{a  & b \ar[l]_{f} \ar[r]^{g} &  c.}$$
This category is a graded poset with increasing degree function
indicated by the subscripts:
$$ \xymatrix{a_1  & b_0 \ar[l]_{f} \ar[r]^{g} &  c_1.}$$
Although in this case it is trivial to compute the derived
functors we apply the earlier propositions to have a taste of it.
So we filter by the ending object ($\sigma_n$) to obtain a
homological type spectral sequence $E^*_{*,*}$ converging to
$\limdir_i F$. We do not filter a second time since this case is
too simple. The column $E^1_{p,*}$ is given by the homology of a
simplicial abelian group. In fact, by Theorem
\ref{simp_normalization_thm}, $E^1_{1,*}$ is the homology of the
quotient chain complex of non-degenerated simplices ending in
degree $1$, that is, the homology of
$$ \xymatrix{ ...0 \ar[r]^<<<<{0}& F_f\oplus F_g \ar[rr]^{F(f)\oplus F(g)} && F_a\oplus F_c. }$$
Analogously, $E^1_{0,*}$ is the homology of the quotient chain complex of non-degenerated simplices ending in degree $0$,
 that is, the homology of
$$ \xymatrix{ ...0 \ar[r]^{0}& F_b. }$$
So the page $E^1$ looks like
\begin{footnotesize}
$$ \xymatrix@=2pt{
 &   &\\
 & 0 & 0 & 0 & 0 & \\
 \ar@<-5pt>[rrrrr] & 0  & F_b & \ar[l]_>>{d_1} ker(F(f)\oplus F(g))  & 0 & \\
 & 0 &  0   & coker(F(f)\oplus F(g))            & 0 & \\
 & 0 & 0 & 0 & 0 & \\
 &   & \ar@<8pt>[uuuuu]  \\
}$$
\end{footnotesize}
and the only nontrivial differential $d_1$ is the restriction of
$d_{F(b)}+id_{F(b)}:F(b)\oplus F(b)\rightarrow F(b)$. It is clear
that the spectral sequence collapses at $E^2$. Clearly $\limdir_i
F=0$ for $i\geq 2$. Notice that if $F(f)$ and $F(g)$ are
monomorphisms then $E^1_{1,0}=0$ and $\limdir_1 F=0$. The
extension problem for $\limdir_0 F=\limdir F$ gives the short
exact sequence:
\begin{footnotesize}
$$
0 \rightarrow F(b)/\ker(F(f))+\ker(F(g)) \rightarrow \limdir
F=coker(F(f)\times F(g)) \rightarrow coker(F(f)\oplus
F(g))\rightarrow 0.
$$
\end{footnotesize}
Using any of the spectral sequences for $\liminv^i F$ it is
straightforward that
\begin{numcases}{{\liminv}^i F=}
F(b) & for $i=0$ \nonumber \\
0 & if $i\geq 1$. \nonumber
\end{numcases}
\end{Ex}

\begin{Ex}
\textbf{Pullback:} Consider $\P$ the ``pullback category":
$$ \xymatrix{ & c \ar[d]^{g} \\ a \ar[r]_{f} & b.}$$
This category is a graded poset with decreasing degree function
indicated by the subscripts:
$$ \xymatrix{ & c_1 \ar[d]^{g} \\ a_1 \ar[r]_{f} & b_0.}$$
We filter by the initial object ($\sigma_0$) to obtain a
cohomological type spectral sequence $E_*^{*,*}$ converging to
$\liminv^i F$. The column $E_1^{p,*}$ is given (Theorem
\ref{cosimp_normalization_thm}) by the cohomology of the
normalized cochain complex of simplices beginning in degree $p$.
So $E_1^{1,*}$ is the cohomology of
$$ \xymatrix{ ...0 & \ar[l]_{0} F^f\oplus F^g && F^a\oplus F^c \ar[ll]_{F(f)\oplus F(g)}, }$$
and $E_1^{0,*}$ the cohomology of
$$ \xymatrix{ ...0 & \ar[l]_{0} F^b.}$$
So the page $E_1$ looks like
\begin{footnotesize}
$$ \xymatrix@=4pt{
 &   &\\
 & 0 & 0 & 0 & 0 & \\
 \ar@<-5pt>[rrrrr] & 0  & F^b \ar[r]^<<{d^1}& coker(F(f)\oplus F(g)) & 0 & \\
 & 0 & 0  & ker(F(f)\oplus F(g)) & 0 & \\
 & 0 & 0 & 0 & 0 & \\
 &   & \ar@<8pt>[uuuuu]  \\
}$$
\end{footnotesize}
and the only nontrivial differential $d_1$ is induced by
$id_{F(b)}\times id_{F(b)}:F(b)\rightarrow F(b)\oplus F(b)$. The
spectral sequence collapses at $E_2$. Notice that if $F(f)$ and
$F(g)$ are epimorphisms then $E_1^{1,0}=0$ and $\liminv^1 F$=0.

The extension problem for $\liminv^0 F=\liminv F$ gives the short
exact sequence:
\begin{footnotesize}
$$
0 \rightarrow ker(F(f)\oplus F(g)) \rightarrow \liminv
F=ker(F(f)-F(g)) \rightarrow \im(F(f))\cap\im(F(g))\rightarrow 0.
$$
\end{footnotesize}

For $\limdir_i F$ it holds that
\begin{numcases}{{\limdir}_i F=}
F(b) & for $i=0$ \nonumber \\
0 & if $i\geq 1$. \nonumber
\end{numcases}
\end{Ex}

\begin{Ex}
\textbf{Telescope:} Consider $\P$ the ``telescope category":
$$ \xymatrix{a_0 \ar[r]^{f_1} & a_1 \ar[r]^{f_2} & a_2 \ar[r]^{f_3} &  a_3 \ar[r]^{f_4} & a_4 ...}$$
where the subscript indicates the name of the object and the value
of an increasing degree function which makes $\P$ a graded poset.
The same homological type spectral sequence of the pushout example
for $\limdir_i F$ has as column $E^1_{p,*}$ the homology of the
normalized chain complex of simplices ending in degree $p$. These
chain complexes become more and more complicated as $p$ grows. To
get more insight we filter a second time by the condition
$deg(\sigma_0)\geq p'$ to obtain for each $p$ a cohomological type
spectral sequence $(E_p)_*^{*,*}$ converging to $E^1_{p,*}$. The
column $(E_p)_*^{p',*}$ is given by the homology of the simplicial
abelian group formed by the simplices that end on degree $p$ and
begin on degree $p'$.

We want to have a look at $\limdir_1 F$. Here we write an informal
discussion, for precise statements look at the proof of Theorem
\ref{pro_acyclic_graded}. The contributions to $\limdir_1 F$ from
the second spectral sequence for $p'<p$ come from the homology at
degree $1$ of normalized chain complexes
$$ \xymatrix{ ...\ar[r]^<<{d} & F_\sigma = F(a_{p'}) \ar[r]^>>{d} & 0 }$$
where $\sigma=\xymatrix{ a_{p'} \ar[r]^{f_p ..\circ f_{p'+1}} & a_p }$ with $p'<p$.
If $p'<p-1$, for $x \in F(p')$ take the $2$-chain
$$y=i_{\sigma'}(-x)$$
where $\sigma'=\xymatrix{ a_{p'} \ar[r]^{f_{p'+1}} & a_{p'+1} \ar[r]^{f_p..\circ f_{p'+2}} &a_p}$.

Then $d(y)=d_0(y)-d_1(y)+d_2(y)=0-(-x)+0=x$, and so if $p'<p-1$
there is no contribution to $\limdir_1 F$. Notice that the key
fact in the argument above is that every morphism of degree
greater than $1$ can be written as composition of morphisms of
degree $1$. Notice that at this point the calculus of $\limdir_1
F$ have been simplified a lot.

The computation of $\liminv^1 F$ is simplified too. In this case
the contributions to $\liminv^1 F$ come from the cohomology of
normalized cochain complexes
$$ \xymatrix{ ... &\ar[l]_<<{d}  F^\sigma=F(a_p) & \ar[l]_<<{d}  0}$$
where $\sigma=\xymatrix{ a_{p'} \ar[r]^{f_p..\circ f_{p'+1}} & a_p }$ with $p'<p$.
For $p'<p-1$, if $x \in F(p)$ is a cohomological class then $d(x)=0$, and so
$\pi^{\sigma'}(d(x))=x=0$
where $\sigma'=\xymatrix{ a_{p'} \ar[r]^{f_{p'+1}} & a_{p'+1} \ar[r]^{f_p..\circ f_{p'+2}} &a_p}$.
Thus if $p'<p-1$ there is no contribution to $\liminv^1 F$.
\end{Ex}

\begin{Ex}\label{ss_example_cycle}
\textbf{Cycle:} Consider $\P$ the following category
$$ 
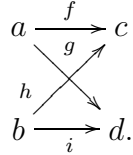
\xymatrix{
a \ar[rd]^<<<<{g} \ar[r]^{f} & c \\
b \ar[ru]^<<{h} \ar[r]_{i} & d. \\
}
$$
It is a graded poset with increasing degree function
$$ 
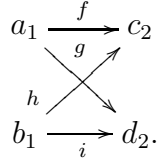
\xymatrix{
a_1 \ar[rd]^<<<<{g} \ar[r]^{f} & c_2 \\
b_1 \ar[ru]^<{h} \ar[r]_{i} & d_2. \\
}$$ The associated undirected graph has a cycle. A direct
computation shows that $\limdir_1 F$ consists of the tuples
$$(x_f,x_g,x_h,x_i)\in F_f\times F_g\times F_h\times F_i=F(a)\times F(a)\times F(b)\times F(b)$$ such that
$x_f+x_g=0$, $F(f)(x_f)+F(h)(x_h)=0$, $x_h+x_i=0$ and
$F(g)(x_g)+F(i)(x_i)=0$. This system of equations may have
non-trivial solutions. For example, if we consider the constant
functor $F=c_\Z$ then $\limdir_1 F=H_1(|N\P|)=H_1(S^1)=\Z$, and so
the solution is cyclic infinite. What happens if we add an initial
object $e$ to $\P$?
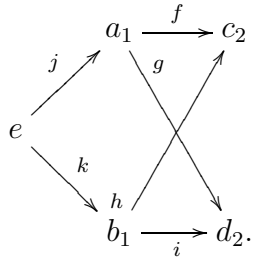
$$ \xymatrix{
 & a_1 \ar[rdd]^<<<<{g} \ar[r]^{f} & c_2 \\
e \ar[ru]^{j}\ar[rd]^{k} &                               & \\
 & b_1 \ar[ruu]^<{h} \ar[r]_{i} & d_2. \\
}$$ In case $F$ is the constant functor of value $\Z$ then
$\limdir_1 F$ is the first group of homology of the cone over
$S^1$, which is a contractible space, and so $\limdir_1 F=0$. What
happens for an arbitrary $F$? The image by the differential from
$C_2(\P,F)$ to $C_1(\P,F)$ of the tuple
$$(x_{jf},x_{jg},x_{kh},x_{ki})\in F_{jf}\times F_{jg}\times F_{kh}\times F_{ki}=F(e)\times F(e)\times F(e)\times F(e)$$
is the tuple
$$(F(j)(x_{jf}),F(j)(x_{jg}),F(k)(x_{kh}),F(k)(x_{ki}),x_{jg}+x_{jg},x_{kh}+x_{ki},-(x_{jf}+x_{kh}),-(x_{jg}+x_{ki}))$$
of $$F_f\times F_g\times F_h\times F_i\times F_j\times F_k\times
F_{f\circ j}\times F_{i\circ k}$$ which equals $$F(a)\times
F(a)\times F(b)\times F(b)\times F(e)\times F(e)\times F(e)\times
F(e).$$ Applying the arguments of the preceding examples we have
that for a class
$$[(x_f,x_g,x_h,x_i,x_j,x_k,x_{f\circ j},x_{i\circ k})]$$ in the kernel of the differential at $C_1(\P,f)$
we can take
\begin{itemize}
\item $x_{f\circ j}=x_{i\circ k}=0$ because $f\circ j$ and $i\circ k$ can be factored by morphisms of lower degree.
\item $x_j=x_k=0$ if $F(j)$ and $F(k)$ are monomorphisms, as $j$ and $k$ are the only arrows arriving to their
ending objects.
\end{itemize}
So supposing that $F(j)$ and $F(k)$ are monomorphisms we can take as representative
of a class in the kernel a tuple
$(x_f,x_g,x_h,x_i,0,0,0,0)$ such that
$x_f+x_g=0$, $F(f)(x_f)+F(h)(x_h)=0$, $x_h+x_i=0$ and $F(g)(x_g)+F(i)(x_i)=0$, as before.

Now, suppose that for every $x_a\in F(a)$ and every $x_b\in F(b)$
such that $F(f)(x_a)=F(h)(x_b)$ there exists $x_e\in F(e)$ such
that $F(j)(x_e)=x_a$ and $F(k)(x_e)=x_b$. This condition is
natural as it is related with the projectiveness of $F$ in
$\Ab^\P$ (cf. Section \ref{section_projective}).

For a tuple $(x_f,x_g,x_h,x_i,0,0,0,0)$ in the kernel, we have
that $F(f)(x_f)=F(h)(-x_h)$ and so, by hypothesis, exists $x_e\in
F(e)$ with $F(j)(x_e)=x_f$ and $F(k)(x_e)=-x_h$. Take now the
$2$-chain of $C_2(\P,F)$ $y=(x_e,-x_e,-x_e,x_e)$. Then the
differential of $y$ is
$$(x_f,-x_f,x_h,-x_h,x_e-x_e,-x_e+x_e,-(x_e-x_e),-(-x_e+x_e))$$
which equals $(x_f,x_g,x_h,x_i,0,0,0,0)$. Thus with these
assumptions $\limdir_1 F=0$.
\end{Ex}

\chapter{Higher direct
limits}\label{section_on direct limit}

\section{Projective objects in
{${\Ab}^{\mathcal{P}}$}.}\label{section_projective} Consider the abelian category
$\Ab^\P$ for some graded poset $\P$. In this section we shall
determine the projective objects in $\Ab^\P$. One of the main
features of projective objects is that derived functors vanish on
them. Recall that in $\Ab$ the projective objects are well known,
and are exactly the free abelian groups. Along the rest of the
section $\P$ denotes a graded poset.

Suppose $F\in \Ab^\P$ is projective. How does $F$ look? Consider
$i_0\in \Ob(\P)$. We show that the quotient of $F(i_0)$ by the
images of the non-identity morphisms arriving to $i_0$ is
projective. To prove it, write

\begin{Defi}\label{defi_im}
$\im(i_0)=\sum_{i\stackrel{\alpha}\rightarrow i_0,\alpha\neq
1_{i_0}} \im F(\alpha)$ (or $\im(i_0)=0$ if the index set of the
sum is empty) and $\coker(i_0)=F(i_0)/\im(i_0)$.
\end{Defi}

Let $F\in \Ab^\P$ be a projective functor. For any diagram in
$\Ab$ as the following

$$ \xymatrix{& \coker(i_0) \ar[d]^{\sigma_0} \ar@{-->}[dl]^{\rho_0} &\\
A_0 \ar[r]^{\pi_0} & B_0 \ar[r] & 0 }
$$

we want to find $\rho_0$ that makes it commutative. Consider the
atomic functors $A,B:\P\rightarrow \Ab$ which take the values on
objects
\begin{numcases}{A(i)=}
A_0 & for $i=i_0$ \nonumber \\
0 & for $i\neq i_0$ \nonumber
\end{numcases}
\begin{numcases}{B(i)=}
B_0 & for $i=i_0$ \nonumber \\
0 & for $i\neq i_0$ \nonumber
\end{numcases}
and on morphisms
\begin{numcases}{A(\alpha)=}
1_{A_0} & for $\alpha=1_{i_0}$ \nonumber \\
0 & for $\alpha\neq 1_{i_0}$ \nonumber
\end{numcases}
\begin{numcases}{B(\alpha)=}
1_{B_0} & for $\alpha=1_{i_0}$ \nonumber \\
0 & for $\alpha\neq 1_{i_0}$ \nonumber
\end{numcases}
and the natural transformations $\sigma:F\Rightarrow B$ and $\pi:A\Rightarrow B$ given by
\begin{numcases}{\sigma(i)=}
\sigma_0\circ p & for $i=i_0$ \nonumber \\
0 & for $i\neq i_0$ \nonumber
\end{numcases}
\begin{numcases}{\pi(i)=}
\pi_0 & for $i=i_0$ \nonumber \\
0 & for $i\neq i_0$ \nonumber
\end{numcases}
where $p$ is the projection $F(i_0)\twoheadrightarrow
\coker(i_0)$. $A\stackrel{\pi}\Rightarrow B\Rightarrow 0$ is exact
as $A_0\stackrel{\pi_0}\rightarrow B_0\rightarrow 0$ is. It is
straightforward that $\pi$ is a natural transformation. The key
point in checking that $\sigma$ is a natural transformation is
that for $\alpha:i_1\rightarrow i_0$, $\alpha\neq 1_{i_0}$ the
diagram
$$
\xymatrix{ F(i_1)\ar[r]^{F(\alpha)}\ar[d]^{\sigma(i_1)} &F(i_0)\ar[d]^{\sigma_0\circ p}\\
           B(i_1)\ar[r]^{0}& B_0 }
$$
must commute. And it does because $p\circ F(\alpha)=0$ for every $\alpha\neq 1_{i_0}$.

So, as $F$ is projective, this data gives a natural transformation
$\rho$ which makes commutative the diagram of natural
transformations
$$
\xymatrix{
             & F \ar@{=>}[d]^\sigma \ar@{==>}[dl]^\rho &\\
A \ar@{=>}[r]^\pi & B \ar@{=>}[r]                           & 0 }
$$
which restricts over $i_0$ to
$$ \xymatrix{
& F(i_0) \ar[d]^{p}\ar[ddl]_{\rho(i_0)}\\
& \coker(i_0) \ar[d]^{\sigma_0} \ar@{-->}[dl]^{\rho_0} &\\
A_0 \ar[r]^{\pi_0} & B_0 \ar[r] & 0. }
$$
Then $\rho_0$ exists if and only if $\ker(p)=\im(i_0)\leq
\ker\rho(i_0)$. To check that this condition holds take
$x=\sum_{j=1,..,k} F(\alpha_j)(x_j)$ in $\im(i_0)$ for
$1_{i_0}\neq\alpha_j:i_j\rightarrow i_0$ $j=1,..,k$. Then
$$
\rho(i_0)(x)=\sum_{j=1,..,k}\rho(i_0)(F(\alpha_j)(x_j))
 = \sum_{j=1,..,k}A(\alpha_j)(\rho(i_j)(x_j))
 = \sum_{j=1,..,k}A(\alpha_j)(0)
 = 0. \nonumber \\
$$
We have just proven
\begin{Lem}\label{lem_proj_pseudo_graded_coker_projective}
Let $F:\P\rightarrow \Ab$ be a projective functor over a graded
poset $\P$. Then $\coker(i_0)$ is projective for every object
$i_0\in \Ob(\P)$.
\end{Lem}
This means that we can write
$$ F(i_0)=\im(i_0)\oplus\coker(i_0)$$
with $\coker(i_0)$ free for every $i_0\in \Ob(\P)$, and also that
\begin{Ex}
For the category $\P$ with shape
$$\cdot \rightarrow \cdot $$
the functor $F:\P\rightarrow \Ab$ with values
$$\Z \stackrel{\times n}\rightarrow \Z $$ is not projective as $\coker$ on the right object equals
the non-free abelian group $\Z/n$.
\end{Ex}

Now that we know a little about the values that a projective
functor $F:\P\rightarrow \Ab$ takes on objects we can wonder about
the values $F(\alpha)$ for $\alpha\in \Hom(\P)$. Do they have any
special property? Recall that a feature of graded posets is that
there is at most one arrow between any two objects, and also that
\begin{Rmk}\label{im_on_graded}
If $\P$ is graded then for any $i_0\in \Ob(\P)$
$$\im(i_0)=\sum_{i\stackrel{\alpha}\rightarrow i_0,deg(\alpha)=1} \im F(\alpha) $$
because every morphism factors as composition of morphisms of
degree $1$.
\end{Rmk}

We prove that the following property holds for $F$:

\begin{Defi}\label{property_pseudo-projective}
Let $F:\P\rightarrow \Ab$ be a functor over a graded poset $\P$
with degree function $deg$. Given $d\geq 0$ we say that $F$ is
\emph{$d$-pseudo-projective} if for any $i_0\in \Ob(\P)$ and $k$
different objects $i_j\in \Ob(\P)$, arrows
$\alpha_j:i_j\rightarrow i_0$ with $deg(\alpha_j)=d$, and $x_j\in
F(i_j)$ $j=1,..,k$ such that
$$\sum_{j=1,..,k}F(\alpha_j)(x_j)=0$$
we have that $x_j\in\im(i_j)$ $j=1,..,k$. If $F$ is
$d$-pseudo-projective for each $d\geq 0$ we call $F$
\emph{pseudo-projective}.
\end{Defi}

\begin{Rmk}\label{pseudo_injective_is injective}
In case $k=1$ and $\im(i_1)=0$ the condition states that
$F(\alpha_1)$ is a monomorphism. Notice that any functor is
$0$-pseudo-projective as the identity is a monomorphism.
\end{Rmk}

Before proving that projective functors $F$ over a graded poset
verify this property we define two functors $\coker$ and $\coker'$
and natural transformations $\sigma$ and $\pi$ that fit in the
diagram
$$
\xymatrix{
             & F \ar@{=>}[d]^{\sigma} &\\
\coker' \ar@{=>}[r]^\pi & \coker \ar@{=>}[r] & 0 }
$$
for any functor $F:\P\rightarrow \Ab$ with $\P$ a graded poset. We
begin defining $\coker$. Because for every $\alpha:i_1\rightarrow
i_0$ holds that $F(\alpha)(\im(i_1))\leq \im(i_0)$ we can factor
$F(\alpha)$ as in the diagram
$$
\xymatrix{ F(i_1)\ar@{->>}[d]\ar[r]^{F(\alpha)} & F(i_0)\ar@{->>}[d]\\
\coker(i_1)\ar[r]^{\overline{F(\alpha)}} & \coker(i_0).
 }
$$
In fact, if $\alpha\neq 1_{i_1}$, then $\overline{F(\alpha)}\equiv
0$ by definition. Because the identity $1_{i_0}$ cannot be
factorized (by non-identity morphisms) in a graded poset then we
have a functor $\coker$ with value $\coker(i)$ on the object $i$
of $\P$ and which maps the non-identity morphisms to zero.
$\coker$ is a kind of ``discrete" functor. Also it is clear that
there exists a natural transformation $\sigma:F\Rightarrow \coker$
with $\sigma(i)$ the projection $F(i)\twoheadrightarrow
\coker(i)$.

Now we define $\coker'$ from $\coker$ in a similar way as free
diagrams are constructed. Let $\coker'$ be defined on objects by
$$\coker'(i_0)=\bigoplus_{\alpha:i\rightarrow i_0} \coker(i).$$
For $\beta\in \Hom(\P)$, $\beta:i_1\rightarrow i_0$,
$\coker'(\beta)$ is the only homomorphism which makes commute the
diagram
$$
\xymatrix{ \coker'(i_1)\ar[r]^{\coker'(\beta)} & \coker'(i_0) \\
\coker(i)\ar@{^{(}->}[u]\ar[r]^{1}& \coker(i)\ar@{^{(}->}[u]
 }
$$
for each $\alpha:i\rightarrow i_1$. In the bottom row of the
diagram, the direct summands $\coker(i)$ of $\coker'(i_1)$ and
$\coker'(i_0)$ correspond to $\alpha:i\rightarrow i_1$ and to the
composition $i\stackrel{\alpha}\rightarrow
i_1\stackrel{\beta}\rightarrow i_0$ respectively.

Then there exists a candidate to natural transformation
$\pi:\coker'\Rightarrow \coker$ which value $\pi(i)$ is the
projection $\pi(i):\coker'(i) \twoheadrightarrow \coker(i)$ onto
the direct summand corresponding to $1_i:i\rightarrow i$. Thus,
$\pi$ is a natural transformation if for every
$\beta:i_1\rightarrow i_0$ with $i_1\neq i_0$ the following
diagram is commutative
$$
\xymatrix{\coker'(i_1)\ar@{->>}[d]^{\pi(i_1)}\ar[r]^{\coker'(\beta)} & \coker'(i_0)\ar@{->>}[d]^{\pi(i_0)}\\
\coker(i_1)\ar[r]^{0} & \coker(i_0).}
$$
It is clear that this square commutes if the identity $1_{i_0}$
cannot be factorized (by non-identity morphisms), and this holds
in a graded poset.

Now we have the commutative triangle
$$
\xymatrix{
             & F \ar@{=>}[d]^{\sigma} \ar@{==>}[dl]^\rho &\\
\coker' \ar@{=>}[r]^\pi & \coker \ar@{=>}[r] & 0 }
$$
where the natural transformation $\rho$ exists because $F$ is
projective. To prove that $F$ is $d$-pseudo-projective for some
$d\geq 0$ take $i_0\in \Ob(\P)$, $k$ objects $i_1,..,i_k$, arrows
$\alpha_j:i_j\rightarrow i_0$ with $deg(\alpha_j)=d$ and elements
$x_j\in F(i_j)$ for $j=1,..,k$ such that
$$\sum_{j=1,..,k}F(\alpha_j)(x_j)=0.$$

To visualize what is going on consider the diagram above near
$i_0$ for $k=2$
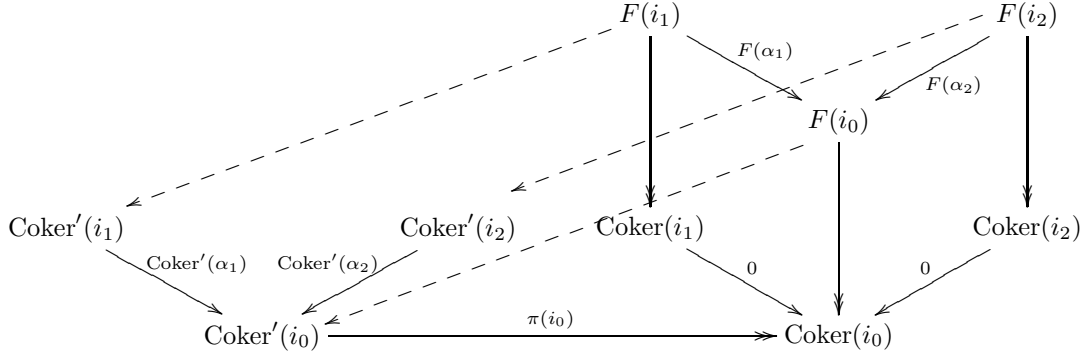
\begin{footnotesize}
$$
\xymatrix{
&&&  F(i_1)\ar@{-->}[ddlll]\ar@{->>}[dd]\ar[rd]^{F(\alpha_1)}& & F(i_2)\ar@<-6pt>@{-->}[ddlll]\ar@{->>}[dd]\ar[dl]^{F(\alpha_2)}\\
&&&  & F(i_0)\ar@<4pt>@{-->}[ddlll]\ar@{->>}[dd] &\\
\coker'(i_1)\ar[rd]      ^{\coker'(\alpha_1)} & &
\coker'(i_2)\ar[ld]_{\coker'(\alpha_2)}&
 \coker(i_1)\ar[rd]^{0}& & \coker(i_2)\ar[dl]_{0}\\
& \coker'(i_0)\ar@{->>}[rrr]^{\pi(i_0)} &&
 & \coker(i_0) &
}
$$
\end{footnotesize}
where $\pi$ is not drawn completely for clarity. Recall that we
are supposing that $\{x_1,..,x_k\}$ is such that
$\sum_{j=1,..,k}F(\alpha_j)(x_j)=0$. Then
$$
0=\rho(i_0)(0)=\sum_{j=1,..,k}\rho(i_0)(F(\alpha_j)(x_j))
 = \sum_{j=1,..,k}\coker'(\alpha_j)(\rho(i_j)(x_j)).
$$

Now consider the projection $p_{j_0}$ for $j_0\in\{1,..,k\}$ from
$\coker'(i_0)$ onto the direct summand
$\coker(i_{j_0})\hookrightarrow \coker'(i_0)$ which corresponds to
$\alpha_{j_0}:i_{j_0}\rightarrow i_0$
$$\coker'(i_0)\stackrel{p_{j_0}}\twoheadrightarrow \coker(i_{j_0}).$$
Then
\begin{equation}\label{projiszero}
 0=p_{j_0}(0)=p_{j_0}(\rho(i_0)(0))=\sum_{j=1,..,k}
p_{j_0}(\coker'(\alpha_j)(\rho(i_j)(x_j))).
\end{equation}

For any $y=\bigoplus_{\alpha:i\rightarrow i_j}
y_\alpha\in\coker'(y_j)$
$$p_{j_0}(\coker'(\alpha_j)(y))=\sum_{\alpha:i\rightarrow i_j,\alpha_j\circ\alpha=\alpha_{j_0}}
y_\alpha.$$

So if $y_j=\rho(i_j)(x_j)=\bigoplus_{\alpha:i\rightarrow i_j}
y_{j,\alpha} \in \coker'(i_j)$ then
$$p_{j_0}(\coker'(\alpha_j)(\rho(i_j)(x_j)))=\sum_{\alpha:i\rightarrow i_j,\alpha_j\circ\alpha=\alpha_{j_0}}
y_{j,\alpha}.$$

This last sum runs over $\alpha:i_{j_0}\rightarrow i_j$ such that
the following triangle commutes
$$
\xymatrix{ i_{j_0}\ar[rd]^{\alpha}\ar[rr]^{\alpha_{j_0}} && i_0.\\
& i_j \ar[ru]^{\alpha_j} &}
$$
Because we are in a graded poset and $deg(i_j)=d$ for each
$j=1,..,k$ then the only chance is $i_j=i_{j_0}$ and
$\alpha=1_{i_{j_0}}$. Because the objects $i_1,..,i_k$ are
different this implies that $j=j_0$ too. Thus
\begin{numcases}{p_{j_0}(\coker'(\alpha_j)(\rho(i_j)(x_j)))=}
y_{j_0,1_{i_{j_0}}} & for $j=j_0$ \nonumber \\
0 & for $j\neq j_0$ \nonumber
\end{numcases}
and Equation (\ref{projiszero}) becomes
$$
0=p_{j_0}(0)=y_{j_0,1_{i_{j_0}}}.
$$
Notice now that $y_{j_0,1_{i_{j_0}}}$ is the evaluation of
$\pi(i_{j_0})$ on $y_{j_0}=\rho(i_{j_0})(x_{j_0})$ and then
$$
0=y_{j_0,1_{i_{j_0}}}=\pi(i_{j_0})(\rho(i_{j_0})(x_{j_0}))=\sigma_{i_{j_0}}(x_{j_0}).
$$
This last equation means that $x_{j_0}$ goes to zero by the
projection $F(i_{j_0})\twoheadrightarrow
\coker(i_{j_0})=F(i_{j_0})/\im(i_{j_0})$, and then
$$x_{i_{j_0}}\in \im(i_{j_0}).$$
As $j_0$ was arbitrary this completes the proof of

\begin{Lem}\label{lem_proj_pseudo_graded}
Let $F:\P\rightarrow \Ab$ be a projective functor over a graded
poset $\P$. Then $F$ is pseudo-projective.
\end{Lem}

\begin{Ex}
For the category $\P$ with shape
$$\cdot \rightarrow \cdot $$
the functor $F:\P\rightarrow \Ab$ with values
$$\Z \stackrel{red_n}\rightarrow \Z/n $$ is not projective as
$red_n$ is not injective, in spite of the $\coker$'s are $\Z$ and
$0$, which are free abelian.
\end{Ex}

Now we define \emph{pre-projective} objects
\begin{Defi}\label{Defi_preprojective}
Let $F:\P\rightarrow \Ab$ be a functor over a graded poset $\P$.
We call $F$ \emph{pre-projective} if
\begin{enumerate}
\item for any $i_0\in \Ob(\P)$ $\coker(i_0)$ is projective.
\label{Defi_preprojective_1}\item $F$ is pseudo-projective.
\label{Defi_preprojective_2}
\end{enumerate}
\end{Defi}

Till now we have obtained that projective functors $\P\rightarrow
\Ab$ over graded posets are pre-projective. In fact, as the next
proposition shows, the restriction we did to graded posets is
worthwhile:

\begin{Prop}\label{pro_projective_bounded_graded}
Let $F:\P\rightarrow \Ab$ be a pre-projective functor over a
graded poset $\P$. If $\P$ is bounded below then $F$ is
projective.
\end{Prop}
\begin{proof}
We can suppose that the degree function $deg$ on $\P$ is
increasing and takes values $\{0,1,2,3,...\}$, and that
$\Ob_0(\P)\neq \emptyset$.

To see that $F$ is proyective in $\Ab^\P$, given a diagram of
functors with exact row as shown, we must find a natural
transformation $\rho:F\Rightarrow A$ making the diagram
commutative:
$$
\xymatrix{
             & F \ar@{=>}[d]^\sigma \ar@{==>}[dl]^\rho &\\
A \ar@{=>}[r]^\pi & B \ar@{=>}[r]                           & 0. }
$$
We define $\rho$ inductively, beginning on objects of degree $0$
and successively on object of degrees $1,2,3,..$.

So take $i_0\in \Ob_0(\P)$ of degree $0$, and restrict to the
diagram in $\Ab$ over $i_0$. By Definition
\ref{Defi_preprojective}, as $\im(i_0)=0$, $F(i_0)=\coker(i_0)$ is
projective. So we can close the following triangle with a
homomorphism $\rho(i_0)$
$$
\xymatrix{
             & F(i_0) \ar[d]^{\sigma(i_0)} \ar@{-->}[dl]_{\rho(i_0)} &\\
A(i_0) \ar[r]^{\pi(i_0)} & B(i_0) \ar[r] & 0. }
$$
As there are no arrows between degree $0$ objects we do not worry
about $\rho$ being a natural transformation. Now suppose that we
have defined $\rho$ on all objects of $\P$ of degree less than $n$
($n\geq 1$), and that the restriction of $\rho$ to the full
subcategory generated by these objects is a natural transformation
and verifies $\pi\circ \rho=\sigma$.

The next step is to define $\rho$ on degree $n$ objects. So take
$i_0\in \Ob_n(\P)$ and consider the splitting
$$F(i_0)=\im(i_0)\oplus\coker(i_0)$$
where
$$\im(i_0)=\sum_{i\stackrel{\alpha}\rightarrow i_0,deg(\alpha)=1} \im F(\alpha).$$
To define $\rho(i_0)$ such that it makes commutative the diagram
$$
\xymatrix{
             & \im(i_0)\oplus\coker(i_0) \ar[d]^{\sigma(i_0)} \ar@{-->}[dl]_{\rho(i_0)} &\\
A(i_0) \ar[r]^{\pi(i_0)} & B(i_0) \ar[r] & 0,}
$$
we define it on $\im(i_0)$ and $\coker(i_0)$ separately. For
$\coker(i_0)$, as it is a projective abelian group, we define it
by any homomorphism that makes commutative the diagram above when
restricted to $\coker(i_0)$. For $\im(i_0)$ take
$x=\sum_{j=1,..,k} F(\alpha_j)(x_j)$ where $\{i_1,..,i_k\}$ are
$k$ different objects, $\alpha_j:i_j\rightarrow i_0$,
$deg(\alpha_j)=1$ and $x_j\in F(i_j)$ for $j=1,..,k$ (see Remark
\ref{im_on_graded}). Then define
$$
\rho(i_0)(x)=\sum_{j=1,..,k} (A(\alpha_j)\circ\rho(i_j))(x_j).
$$
To check that $\rho(i_0)(x)$ does not depend on the choice of the
$i_j$'s, $\alpha_j$'s and $x_j$'s we have to prove that
$$
\sum_{j=1,..,k} F(\alpha_j)(x_j)=0\Rightarrow \sum_{j=1,..,k}
(A(\alpha_j)\circ \rho(i_j))(x_j)=0.
$$

So suppose that
\begin{equation}\label{equ_pre_pro_0}
\sum_{j=1,..,k}F(\alpha_j)(x_j)=0.
\end{equation}
Then using that $F$ is $1$-pseudo-projective and Remark
\ref{im_on_graded} we obtain objects $i_{j,j'}$, arrows
$\alpha_{j,j'}$ of degree $1$, and elements $x_{j,j'}$ for
$j=1,..,k$, $j'=1,..,k_j$ such that
\begin{equation}\label{equ_pre pro_1}
\sum_{j'=1,..,k_j}F(\alpha_{j,j'})(x_{j,j'})=x_j
\end{equation}
for every $j\in\{1,..,k\}$. Notice that possibly not all the
objects $i_{j,j'}$ are different.  Replacing Equation
(\ref{equ_pre pro_1}) in Equation (\ref{equ_pre_pro_0}) we obtain
\begin{equation}\label{equ_pre_pro_2}
\sum_{\textit{$j=1,..,k$, $j'=1,..,k_j$}}
F(\alpha_j\circ\alpha_{j,j'})(x_{j,j'})=0.
\end{equation}
Because in a graded poset there is at most one arrow between two
objects, the condition $i_{j,j'}=i_{j',{j'}'}=i$ implies
$\alpha_j\circ\alpha_{j,j'}=\alpha_{j'}\circ\alpha_{j',{j'}'}:i\rightarrow
i_0$. So, considering objects $i\in \Ob(\P)$, we can rewrite
(\ref{equ_pre_pro_2}) as
\begin{equation}\label{equ pre pro 3}
\sum_{i\in \Ob(\P)} F(\alpha_j\circ\alpha_{j,j'})(\sum_{j,j' |
i_{j,j'}=i} x_{j,j'})=0.
\end{equation}

Call $\{i'_1,..,i'_m\}=\{i_{j,j'}|\textit{$j=1,..,k$,
$j'=1,..,k_j$}\}$ where these sets have $m$ elements. Call
$\beta_l=\alpha_j\circ \alpha_{j,j'}$ if $i'_l=i_{j,j'}$ and
$y_l=\sum_{j,j' | i_{j,j'}=i'_l} x_{j,j'}$ for $l=1,..,m$. Notice
that $deg(\beta_l)=2$ for each $l$. Then Equation (\ref{equ pre
pro 3}) becomes
\begin{equation}\label{equ_pre pro 4}
\sum_{l=1,..,m} F(\beta_l)(y_l)=0.
\end{equation}

Now we repeat the same argument: applying that $F$ is
$2$-pseudo-projective and the Remark \ref{im_on_graded} to
Equation (\ref{equ_pre pro 4}) we obtain objects $i'_{l,l'}$,
arrows $\beta_{l,l'}$ of degree $1$, and elements $y_{l,l'}$ for
$l=1,..,m$, $l'=1,..,k'_l$ such that
\begin{equation}\label{equ_pre pro_5}
\sum_{l'=1,..,k'_l}F(\beta_{l,l'})(y_{l,l'})=y_l
\end{equation}
for every $l\in\{1,..,m\}$. Substituting (\ref{equ_pre pro_5}) in
(\ref{equ_pre pro 4})
$$
\sum_{\textit{$l=1,..,m$,
$l'=1,..,k'_l$}}F(\beta_l\circ\beta_{l,l'})(y_{l,l'})=0.
$$
Now proceed as before regrouping the terms in this last equation.

In a finite number of steps, after a regrouping of terms as above,
we find objects $i''_s$, arrows $\gamma_s$, and elements $z_s$ of
degree $0$ for $s=1,..,r$ which verify an equation
\begin{equation}\label{equ_pre pro_6}
\sum_{s=1,..,r}F(\gamma_s)(z_s)=0.
\end{equation}
Then pseudo-injectivity gives that $z_s\in \im (i''_s)$ for each
$s$. As $deg(i''_s)=0$ then $\im(i''_s)=0$ and so $z_s=0$ (notice
that $z_s=0$ for $s=1,..,r$ does not imply $x_j=0$ for any $j$).

Recall that we want to prove that
\begin{equation}\label{equ pre pro 7}
\sum_{j=1,..,k} (A(\alpha_j)\circ\rho(i_j))(x_j)=0.
\end{equation}
Substituting (\ref{equ_pre pro_1}) in
$\sum_{j=1,..,k}(A(\alpha_j)\circ\rho(i_j))(x_j)$ we obtain
\begin{align}
\sum_{j=1,..,k}(A(\alpha_j)\circ\rho(i_j))(x_j))&=\sum_{j=1,..,k}\sum_{j'=1,..,k_j}(A(\alpha_j)\circ\rho(i_j)\circ F(\alpha_{j,j'}))(x_{j,j'})\notag\\
&=\sum_{j=1,..,k}\sum_{j'=1,..,k_j}(A(\alpha_j)\circ
A(\alpha_{j,j'})\circ\rho(i_{j,j'}) )(x_{j,j'})\notag\\
&=\sum_{j=1,..,k}\sum_{j'=1,..,k_j}(A(\alpha_j\circ
\alpha_{j,j'})\circ\rho(i_{j,j'}) )(x_{j,j'}),\notag
\end{align}
as $\rho$ is natural up to degree less than $n$. Then regrouping
terms
\begin{align}
\sum_{j=1,..,k}\sum_{j'=1,..,k_j}(A(\alpha_j\circ
\alpha_{j,j'})\circ\rho(i_{j,j'}) )(x_{j,j'})&=\sum_{i\in
\Ob(\P)}(A(\alpha_j\circ \alpha_{j,j'})\circ
\rho(i_{j,j'}))(\sum_{j,j' | i_{j,j'}=i} x_{j,j'})\notag\\
&=\sum_{l=1,..,m}(A(\beta_l)\circ \rho(i'_l))(y_l).\notag
\end{align}

Then, after a finite number of steps, we obtain
$$
\sum_{j=1,..,k}(A(\alpha_j)\circ\rho(i_j))(x_j))=\sum_{\textit{$s=1,..,r$}}(A(\gamma_s)\circ\rho(i''_s))(z_s)=0
$$
as $z_s=0$ for each $z=1,..,r$.

So we have checked that $\rho(i_0)(x)$ does not depend on the
choice of $i_j$, $\alpha_j$ and $x_j$. It is straightforward that
$\rho(i_0)$ on $\im(i_0)$ defined in this way is a homomorphism of
abelian groups.

It remains to prove that $\pi(i_0)\rho(i_0)=\sigma(i_0)$ when
restricted to $\im(i_0)$. So take
$x=\sum_{j=1,..,k}F(\alpha_j)(x_j)$ in $\im(i_0)$. Then
\begin{align}
\pi(i_0)(\rho(i_0)(x)) =&\sum_{j=1,..,k}(\pi(i_0)\circ A(\alpha_j)\circ \rho(i_j))(x_j)  \notag \\
=&\sum_{j=1,..,k}(B(\alpha_j)\circ\pi(i_j)\circ\rho(i_j))(x_j) \text{, $\pi$ is a natural transformation}\notag\\
=&\sum_{j=1,..,k}(B(\alpha_j)\circ\sigma(i_j))(x_j) \text{, by the inductive hypothesis} \notag\\
=&\sum_{j=1,..,k}(\sigma(i_0)\circ F(\alpha_j))(x_j) \text{, $\sigma$ is a natural transformation} \notag\\
=&\sigma(i_0)(x) \notag
\end{align}

Defining $\rho(i_0)$ in this way for every $i_0\in \Ob_n(\P)$ we
have now $\rho$ defined on all objects of $\P$ of degree less or
equal than $n$. Finally, to complete the inductive step we have to
prove that $\rho$ restricted to the full subcategory over these
objects is a natural transformation. Take $\alpha:i\rightarrow
i_0$ in this full subcategory. If the degree of $i_0$ is less than
$n$ then the commutativity of
$$
\xymatrix{
F(i)\ar[r]^{F(\alpha)}\ar[d]^{\rho(i)}   &   F(i_0)\ar[d]^{\rho(i_0)} \\
A(i)\ar[r]^{A(\alpha)}                  &   A(i_0) \\
}
$$
is granted by the inductive hypothesis. Suppose that the degree of
$i_0$ is $n$. Take $x'\in F(i)$. Because $\P$ is graded there
exists $\alpha_1:i_1\rightarrow i_0$ of degree $1$ and
$\alpha':i\rightarrow i_1$ such that
$\alpha=\alpha_1\circ\alpha'$:
$$
\xymatrix{ i\ar[rd]^{\alpha'}\ar[rr]^{\alpha} && i_0.\\
& i_1 \ar[ru]^{\alpha_1} &}
$$
Write $x=F(\alpha)(x')=F(\alpha_1)(x_1)$ where
$x_1=F(\alpha')(x')$. Then, by definition of $\rho(i_0)$ on
$\im(i_0)$,
\begin{align}
\rho(i_0)(x)&=(A(\alpha_1)\circ\rho(i_1))(x_1)\notag\\
&=(A(\alpha_1)\circ\rho(i_1))(F(\alpha')(x'))\notag\\
&=(A(\alpha_1)\circ\rho(i_1)\circ F(\alpha'))(x')\notag\\
&=(A(\alpha_1 \circ \alpha')\circ \rho(i))(x')\text{, $\rho$ is natural up to degree less than $n$ }\notag\\
&=(A(\alpha)\circ \rho(i))(x')\notag
\end{align}
and so the diagram commutes.
\end{proof}

\begin{Rmk}
As the following example shows the condition of lower \linebreak
boundedness of $\P$ in Theorem \ref{pro_projective_bounded_graded}
cannot be dropped:

Consider the inverse `telescope category' $\P$ with shape
$$
...\rightarrow \cdot \rightarrow \cdot \rightarrow \cdot
$$
It is a graded poset which is not bounded below. Consider the
functor of constant value $\Z/p$, $c_{\Z/p}:\P\rightarrow \Ab$:
$$
...\rightarrow \Z/p \rightarrow \Z/p \rightarrow \Z/p
$$
It is straightforward that it is a pre-projective functor as all
the cokernels are zero and all the arrows are injective. But it is
not a projective object of $\Ab^\P$ because, in that case, the
adjoint pair $c_\cdot:\Ab\leftrightarrow \Ab^\P:\liminv$ would
give that $\Z/p$ is projective in $\Ab$ (right adjoints preserve
projectives, see \cite[3.2, Ex7]{cohn}).
\end{Rmk}

Thus for the categories that are graded and bounded below we have
the useful

\begin{Cor}\label{cor_projective_bounded_graded}
Let $\P$ be a bounded below graded poset and let $F:\P\rightarrow
\Ab$ be a functor. Then $F$ is projective if and only if it is
pre-projective.
\end{Cor}

This corollary yields the following examples. The degree functions
$deg$ for the bounded below graded posets appearing in the
examples are indicated by subscripts $i_{deg(i)}$ on the objects
$i\in \Ob(\P)$ and take values $\{0,1,2,3,...\}$.

\begin{Ex}\label{examples_projective}
For the `pushout category' $\P$ with shape
$$\xymatrix {a_0 \ar[r]^{f}\ar[d]^{g} & b_1 \\ c_1}$$
a functor $F:\P\rightarrow \Ab$ is projective if and only if
\begin{itemize}
\item $F(a)$, $F(b)/{\im F(f)}$ and $F(c)/{\im F(g)}$ are free
abelian. \item $F(f)$ and $F(g)$ are monomorphisms.
\end{itemize}

\end{Ex}

\begin{Ex}\label{examples_projective_bis}
For the `telescope category' $\P$ with shape
$$ \xymatrix{a_0 \ar[r]^{f_1} & a_1 \ar[r]^{f_2} & a_2 \ar[r]^{f_3} &  a_3 \ar[r]^{f_4} & a_4 ...}$$
a functor $F:\P\rightarrow \Ab$ is projective if and only if
\begin{itemize}
\item $F(a_0)$ is free abelian. \item $F(a_i)/\im F(f_i)$ is free
abelian, $F(f_{i}\circ f_{i-1} \circ ..\circ f_0)$ is a
monomorphism and $\ker F(f_{i}\circ f_{i-1} \circ ..\circ
f_{i-d})\subseteq \im F(f_{i-d-1})$ for $d=0,1,..,i-1$  for each
$i=1,2,3,4,...$.
\end{itemize}
\end{Ex}

\section{Pseudo-projectivity}
Consider a functor $F:\P\rightarrow \Ab$ over a graded poset $\P$.
In this section we look for, and find, conditions on $F$ such that
$\limdir_i F=0$ for $i\geq 1$, i.e., we want conditions such that
the left derived functors of the right exact functor $\limdir$
vanish on $F$. Fix the following notation
\begin{Defi}\label{Defi_direct_acyclic}
Let $\P$ be a graded poset and $F:\P\rightarrow \Ab$. We say $F$
is \emph{$\limdir$-acyclic} if $\limdir_i F=0$ for $i\geq 1$.
\end{Defi}

Recall that for projective objects it holds that any left derived
functor vanishes. So, from Proposition
\ref{pro_projective_bounded_graded}, we obtain firstly that

\begin{Prop}\label{pro_acyclic_projective}
Let $F:\P\rightarrow \Ab$ be a pre-projective functor over a
bounded below graded poset $\P$. Then $F$ is $\limdir$-acyclic.
\end{Prop}

Because being $\limdir$-acyclic is clearly weaker that being
projective we can wonder if is it possible to weaken the
hypothesis on Proposition \ref{pro_acyclic_projective} keeping the
thesis of $\limdir$-acyclicity. The answer is yes and the
following theorem states the appropriate conditions. Notice that
we have removed the condition (\ref{Defi_preprojective_1}) of
Definition \ref{Defi_preprojective}.

\begin{Thm}\label{pro_acyclic_graded}
Let $F:\P\rightarrow \Ab$ be a pseudo-projective functor over a
bounded below graded poset $\P$. Then $F$ is $\limdir$-acyclic.
\end{Thm}
\begin{proof}
We can suppose that the degree function $deg$ on $\P$ is
increasing and takes values $\{0,1,2,3,...\}$, and that
$\Ob_0(\P)\neq \emptyset$. To compute $\limdir_t F$ we use the
(normalized, Remark \ref{normalized_ss}) spectral sequences
corresponding to the third row of Table \ref{tabla_ss} in Chapter
\ref{spectral}. That is, we first filter by the degree of the end
object of each simplex to obtain a homological type spectral
sequence $E^*_{*,*}$. To compute the column $E^1_{p,*}$ we filter
by the degree of the initial object of each object to obtain
cohomological type spectral sequences $(E_p)_*^{*,*}$.

Fix $t\geq 1$. Notice that to prove that $\limdir_t F=0$ is enough
to show that $E^1_{p,t-p}$ is zero for every $p$. The
contributions to $E^1_{p,t-p}$ come from
$(E_p)_{\infty}^{p',p-p'-t}$ for $p'\leq p-t$ (we are using
normalized (Remark \ref{normalized_ss}) spectral sequences). We
prove that
$$
(E_p)_{r}^{p',p-p'-t}=0
$$
if $r$ is big enough for each $p$ and $p'\leq p-t$. This implies
that $\limdir_t F=0$.

Consider the increasing filtration $L^*$ of $C_*(\P,F)$ that gives
rise to the spectral sequence $E^*_{*,*}$. The $n$-simplices are
$$
L^p_n=L^pC_n(\P,F)=\bigoplus_{\sigma \in {N\P}_n,
deg(\sigma_n)\leq p} F_{\sigma}.
$$
For each $p$ we have a decreasing filtration $M_p^*$ of the
quotient $L^p/L^{p-1}$ that gives rise to the spectral sequence
$(E_p)_*^{*,*}$ and which $n$-simplices are
$$
(M_p)^{p'}_n=\bigoplus_{\sigma \in {N\P}_n, deg(\sigma_0)\geq p',
deg(\sigma_n)=p} F_{\sigma}.
$$
For $p'\leq p-t$ the abelian group $(E_p)_r^{p',q'}$ at the
$t=-(p'+q')+p$ simplices is given by
$$
(E_p)_r^{p',q'}= (M_p)^{p'}_t\cap
d^{-1}((M_p)^{p'+r}_{t-1})/(M_p)^{p'+1}_t\cap
d^{-1}((M_p)^{p'+r}_{t-1})+(M_p)^{p'}_t\cap
d((M_p)^{p'-r+1}_{t+1})
$$
where $d$ is the differential of the quotient $L^p/L^{p-1}$
restricted to the subgroups of the filtration $(M_p)^*$. For
$r>p-p'-(t-1)$ there are not $(t-1)$-simplices beginning in degree
at least $p'+r>p-(t-1)$ and ending in degree $p$, i.e.,
$(M_p)^{p'+r}_{t-1}=0$. Because $\P$ is bounded below  for $r$ big
enough $(M_p)^{p'-r+1}_{t+1}=(M_p)^0_{t+1}=(L^p/L^{p-1})_{t+1}$,
i.e., $(M_p)^{p'-r+1}_{t+1}$ equals all the $(t+1)$-simplices that
end on degree $p$. Thus there exists $r$ such that
\begin{equation}
\label{E_p_r_big_dir} (E_p)_r^{p',q'}= (M_p)^{p'}_t\cap
d^{-1}(0)/(M_p)^{p'+1}_t\cap d^{-1}(0)+(M_p)^{p'}_t\cap
d((M_p)^0_{t+1}).
\end{equation}
Fix such an $r$ and take $[x]\in(E_p)_r^{p',q'}$ where
\begin{equation}
\label{equ_exp_x_dir} x=\bigoplus_{\sigma \in {N\P}_t,
deg(\sigma_0)\geq p' , deg(\sigma_t)=p} x_{\sigma}
\end{equation}
and $d(x)=0$. Notice that by definition there is just a finite
number of summands $x_\sigma\neq 0$ in the expression
(\ref{equ_exp_x_dir}) for $x$. We prove that $[x]=0$ in three
steps:

\textbf{Step 1:} In this first step we find a representative $x'$
for $[x]$
$$
x'=\bigoplus_{\sigma \in {N\P}_t, deg(\sigma_0)\geq p' ,
deg(\sigma_t)=p} x'_{\sigma} $$  such that $\deg(\alpha_1)=1$ for
every $\sigma=\xymatrix{ \sigma_0 \ar[r]^{\alpha_1} & \sigma_1
\ar[r]^{\alpha_2}&...\ar[r]^{\alpha_{t-1}}&
\sigma_{t-1}\ar[r]^{\alpha_t} & \sigma_t}$ with $x'_\sigma\neq 0$.

Take $\sigma$ such that $x_\sigma\neq 0$ and suppose that
$deg(\alpha_1)>1$, i.e., $deg(\sigma_0)<deg(\sigma_1)-1$. Then, as
in a graded poset every morphism factors as composition of degree
$1$ morphisms, there exists an object $\sigma_*$ of degree
$deg(\sigma_0)<deg(\sigma_*)<deg(\sigma_1)$ and arrows
$\beta_1:\sigma_0\rightarrow \sigma_*$ and
$\beta_2:\sigma_*\rightarrow \sigma_1$ with $\alpha_1=\beta_2\circ
\beta_1$.
$$\xymatrix{& \sigma_* \ar[rd]^{\beta_2} & \\
\sigma_0 \ar[rr]^{\alpha_1} \ar[ru]^{\beta_1}& & \sigma_1. }$$

Call $\tilde{\sigma}$ to the $(t+1)$-simplex $\sigma=\xymatrix{
\sigma_0 \ar[r]^{\beta_1} & \sigma_* \ar[r]^{\beta_1} & \sigma_1
\ar[r]^{\alpha_2}&...\ar[r]^{\alpha_{t-1}}&
\sigma_{t-1}\ar[r]^{\alpha_t} & \sigma_t}$ and consider the
$(t+1)$-chain of $(M_p)^0_{t+1}$
$y=i_{\tilde{\sigma}}(-x_\sigma)$. Its differential in
$L^p/L^{p-1}$ equals
$$d(y)=d_0(y)-d_1(y)+\sum_{i=2,..,t}(-1)^i d_i(y)=d_0(y)+i_\sigma(x_\sigma)+\sum_{i=2,..,t}(-1)^i d_i(y).$$
Notice that the first morphisms appearing in the simplices
$d_0(\tilde{\sigma})$ and $d_i(\tilde{\sigma})$ for $i=2,..,t$
have degree $deg(\beta_2)$ and $deg(\beta_1)$ respectively, which
are strictly less than $deg(\alpha_1)$. Also notice that $d(y)\in
(M_p)^{p'}_t\cap d((M_p)^0_{t+1})$ (which is  zero in Equation
(\ref{E_p_r_big_dir}) ).

Taking the (finite) sum of the chains $y$ for each term $x_\sigma$
we find that $[x]=[x']$ where
$$
x'=\bigoplus_{\sigma \in {N\P}_t, deg(\sigma_0)\geq p' ,
deg(\sigma_t)=p} x'_{\sigma} $$ and the maximum of the degrees of
the morphisms $\alpha_1$ of the simplices
$$\sigma=\xymatrix{ \sigma_0 \ar[r]^{\alpha_1} & \sigma_1
\ar[r]^{\alpha_2}&...\ar[r]^{\alpha_{t-1}}&
\sigma_{t-1}\ar[r]^{\alpha_t} & \sigma_t}$$ with $x'_\sigma \neq
0$ is smaller than this maximum computed for $x$. So repeating
this process a finite number of times we find a representative as
wished. For simplicity we write also $x$ for this representative.

\textbf{Step 2:} By Step $1$ we can suppose that
$\deg(\alpha_1)=1$ for every $$\sigma=\xymatrix{ \sigma_0
\ar[r]^{\alpha_1} & \sigma_1
\ar[r]^{\alpha_2}&...\ar[r]^{\alpha_{t-1}}&
\sigma_{t-1}\ar[r]^{\alpha_t} & \sigma_t}$$ with $x_\sigma\neq 0$.
Now our objective is to find a representative $x'$ for $[x]$
$$
x'=\bigoplus_{\sigma \in {N\P}_t, deg(\sigma_0)=p' ,
deg(\sigma_t)=p} x'_{\sigma},$$ i.e., such that the expression for
$x'$ runs over simplices $\sigma$ with begin in degree $p'$. Begin
writing $x$ as
$$
x=\bigoplus_{i=p',..,p-t} x_i
$$
where
$$
x_i=\bigoplus_{\sigma \in {N\P}_t, deg(\sigma_0)=i,
deg(\sigma_t)=p} x_\sigma.
$$
Notice that the index $i$ just goes to $p-t$ (and not to $p$)
because we are using normalized (Remark \ref{normalized_ss})
spectral sequences. Now we prove
\begin{Claim}\label{claim_step_2}
For each $i$ from $i=p-t$ to $i=p'$ there exists a representative
$x'_i$ for $[x]$
$$
x'_i=\bigoplus_{\sigma \in {N\P}_t,i\geq deg(\sigma_0)\geq p' ,
deg(\sigma_t)=p} (x'_i)_{\sigma} $$ such that
\begin{equation}
\label{condition_claim_step_2} \textit{$(x'_i)_\sigma\neq 0$ and
$deg(\sigma_0)<i$ imply $deg(\alpha_1)=1$.}
\end{equation}
\end{Claim}
Notice that taking $i=p'$ in the claim, the step $2$ is finished.
The case $i=p-t$ in the claim is fulfilled taking $x'_{p-t}=x$ (by
step $1$). Suppose the statement of the claim holds for $i$. Then
we prove it for $i-1$. We have $x'_i$ such that
$$ x'_i=\bigoplus_{\sigma \in {N\P}_t,i\geq deg(\sigma_0)\geq p' ,
deg(\sigma_t)=p} (x'_i)_{\sigma}, $$ $d(x'_i)=0$ and $[x]=[x'_i]$.
The differential $d$ on $L^p/L^{p-1}$ restricts to
$$d:(M_p)^{p'}_t\rightarrow (M_p)^{p'}_{t-1}$$
and carries $z\in F_\sigma\hookrightarrow \bigoplus_{\sigma \in
{N\P}_t, deg(\sigma_0)\geq p', deg(\sigma_t)=p}
F_{\sigma}=(M_p)^{p'}_t$ to
$$
d(z)=\sum_{j=0,1,..,t-1} (-1)^j d_j(z)
$$
with $d_j(z)\in F_{d_j(\sigma)}\hookrightarrow (M_p)^{p'}_{t-1}$.
Notice that the initial object of $d_j(\sigma)$ is $\sigma_1$ for
$j=0$ and $\sigma_0$ for $j=1,..,t-1$. Also notice that the final
object of $d_j(\sigma)$ is $\sigma_t$ for $j=0,..,t-1$.

By hypothesis $d(x'_i)=0$. So for every $\epsilon\in {N\P}_{t-1}$
with $deg(\epsilon_0)\geq p'$ and $deg(\epsilon_{t-1})=p$ we can
apply the projection
$$\pi_\epsilon: (M_p)^{p'}_{t-1} \twoheadrightarrow F_\epsilon$$ and
obtain $\pi_\epsilon(d(x'_i))=0$. If $deg(\epsilon_0)>i$ then the
remarks on the differential above and condition
(\ref{condition_claim_step_2}) imply that
$$\pi_\epsilon(d(x'_i))=\sum_{\sigma \in{N\P}_t, deg(\sigma_0)=i, d_0(\sigma)=\epsilon}
F(\alpha_1)((x'_i)_{\sigma})
$$
and thus
\begin{equation}\label{equ_epsilon}
0=\sum_{\sigma \in{N\P}_t, deg(\sigma_0)=i, d_0(\sigma)=\epsilon}
F(\alpha_1)((x'_i)_{\sigma})
\end{equation}
for each $\epsilon\in {N\P}_{t-1}$ with $deg(\epsilon_0)>i$ and
$deg(\epsilon_{t-1})=p$. Notice that each summand $(x'_i)_\sigma$
with $\sigma \in {N\P}_t$, $deg(\sigma_0)=i$ and $deg(\sigma)=p$
appears in one and just one equation as (\ref{equ_epsilon}) (take
$\epsilon=d_0(\sigma)$).

Fix an $\epsilon\in {N\P}_{t-1}$ with $deg(\epsilon_0)>i $ and
$deg(\epsilon_{t-1})=p$ and consider the associated Equation
(\ref{equ_epsilon}). Then, as $F$ is
$(i-deg(\epsilon_0))$-pseudo-projective, $(x'_i)_\sigma\in
\im(\sigma_0)$ for every $\sigma \in{N\P}_t$ with
$deg(\sigma_0)=i$ and $d_0(\sigma)=\epsilon$. This means that for
every such a $\sigma$ there exists $k_\sigma$ objects of degree
$(i-1)$, namely $i_\sigma^1,..,i_\sigma^{k_\sigma}$, arrows
$\beta_\sigma^j:i_\sigma^j\rightarrow \sigma_0$ and elements
$x_\sigma^j\in F(i_\sigma^j)$ for $j=1,..,k_\sigma$ such that
\begin{equation}\label{xenim}
(x'_i)_\sigma=\sum_{j=1,..,k_\sigma}F(\beta_\sigma^j)(x_\sigma^j).
\end{equation}
Consider the $(t+1)$-simplices for $j=1,..,k_\sigma$
$$ \sigma^j=\xymatrix{i_\sigma^j\ar[r]^{\beta_\sigma^j} & \sigma_0 \ar[r]^{\alpha_1}
& \sigma_1 \ar[r]^{\alpha_2}&...\ar[r]^{\alpha_{t-1}}&
\sigma_{t-1}\ar[r]^{\alpha_t} &\sigma_t}$$ and the $(t+1)$-chain
of $(M_p)^{i-1}_{t+1}$
$$ y_\sigma=\oplus_{j=1,..,k_\sigma}i_{\sigma^j}(x_\sigma^j).$$

The differential of $y_\sigma$ is
\begin{align}
d(y_\sigma)=&d_0(y_\sigma)+\sum_{j=1,..,t} (-1)^j d_j(y_\sigma)\notag\\
=&d_0(y_\sigma)+R_\sigma  \textit{, where $R_\sigma=\sum_{j=1,..,t} (-1)^j d_j(y_\sigma)$}\notag\\
=&\sum_{j=1,..,k_\sigma}i_{d_0(\sigma^j)}(F(\beta_\sigma^j)(x_\sigma^j))+R_\sigma
\notag\\
=&\sum_{j=1,..,k_\sigma}i_{\sigma}(F(\beta_\sigma^j)(x_\sigma^j))+R_\sigma\notag\\
=&i_{\sigma}(\sum_{j=1,..,k_\sigma}F(\beta_\sigma^j)(x_\sigma^j))+R_\sigma\notag\\
=&i_{\sigma}((x'_i)_\sigma) +R_\sigma\notag
\end{align}

where the last equality is due to (\ref{xenim}). Notice that
$R_\sigma$ lives in the subgroup $\bigoplus_{\sigma \in {N\P}_t,
deg(\sigma_0)=i-1, deg(\sigma_t)=p} F_{\sigma}\subseteq
(M_p)^{p'}_t$ of simplices beginning at degree $(i-1)$. Repeating
the same construction for each $\sigma \in{N\P}_t$ with
$deg(\sigma_0)=i$ and $d_0(\sigma)=\epsilon$ we obtain
$y_\epsilon=\sum_{\sigma} y_\sigma$ such that
$$
d(y_\epsilon)=\oplus_{\sigma \in{N\P}_t, deg(\sigma_0)=i,
d_0(\sigma)=\epsilon} (x'_i)_{\sigma}+R_\epsilon
$$
where $R_\epsilon$ lives in the subgroup $\bigoplus_{\sigma \in
{N\P}_t, deg(\sigma_0)=i-1, deg(\sigma_t)=p} F_{\sigma}\subseteq
(M_p)^{p'}_t$. Repeating the same argument for every $\epsilon\in
{N\P}_{t-1}$ with $deg(\epsilon_0)>i $ and $deg(\epsilon_{t-1})=p$
we obtain $y=\sum_\epsilon y_\epsilon$ such that
$$
d(y)=\oplus_{\sigma \in{N\P}_t, deg(\sigma_0)=i, d_0(\sigma_t)=p}
(x'_i)_{\sigma}+R
$$
where $R$ lives in the subgroup $\bigoplus_{\sigma \in {N\P}_t,
deg(\sigma_0)=i-1, deg(\sigma_t)=p} F_{\sigma}\subseteq
(M_p)^{p'}_t$. By construction $y\in (M_p)^{i-1}_{t+1}\subseteq
(M_p)^0_{t+1}$ and $d(y)\in (M_p)^{i-1}_t\subseteq
(M_p)^{p'}_{t+1}$. Thus $d(y)\in (M_p)^{p'}_t\cap
d((M_p)^0_{t+1})$. Then, by (\ref{E_p_r_big_dir}),
$[x'_i]=[x'_i-d(y)]=[x'_{i-1}]$ where $$x'_{i-1}=\oplus_{\sigma
\in{N\P}_t, i>deg(\sigma_0)\geq p', d_0(\sigma_t)=p}
(x'_i)_{\sigma}+ R$$ is a representative that lives in
$$\bigoplus_{\sigma \in {N\P}_t, i-1\geq deg(\sigma_0)\geq p',
deg(\sigma_t)=p} F_{\sigma}\subseteq (M_p)^{p'}_t$$ as wished.
That condition (\ref{condition_claim_step_2}) holds is clear from
the definition of $x'_{i-1}$.

\textbf{Step 3:} By Step $2$ we can suppose that
$$
x=\bigoplus_{\sigma \in {N\P}_t, deg(\sigma_0)=p' ,
deg(\sigma_t)=p} x_{\sigma}.$$ Our objective now is to see that
there exists $y\in (M_p)^0_{t+1}$ with $d(y)=x$. This implies that
$[x]=0$ and finishes the proof of the theorem. We need the
\begin{Claim}
\label{claim_step_3} There exist chains $x_i\in (M_p)^0_t$ for
$i=p',..,0$ and $y_i \in (M_p)^0_{t+1}$ for $i=p',..,1$ such that
\begin{equation}\label{equ_claim_step_3}
d(y_i)=x_i+x_{i-1}
\end{equation}
for $i=p',..,1$ with $x_{p'}=x$ and $x_0=0$ such that
\begin{enumerate} \item \label{claim_step_3_cond_1} $x_i$ lives on
$\bigoplus_{\sigma \in {N\P}_t, deg(\sigma_0)=i, deg(\sigma_t)=p}
F_{\sigma}\subseteq (M_p)^0_t$ for $i=p',..,0$. \item
\label{claim_step_3_cond_2} $d(x_i)=0$ for $i=p',..,0$.
\end{enumerate}
\end{Claim}
Notice that the claim finishes Step $3$: as $x_0=0$ then
$x_1=d(y_1)$, $x_2=d(y_2)-x_1=d(y_2-y_1)$,
$x_3=d(y_3)-x_2=d(y_3-y_2+y_1)$,..,
$x=x_{p'}=d(y_{p'})-x_{p'-1}=d(y_{p'}-y_{p'-1}+...+(-1)^{p'+1}y_1)$
where $y_{p'}-y_{p'-1}+...+(-1)^{p'+1}y_1\in (M_p)^0_{t+1}$.

Define $x_{p'}\definicio x$. Then condition
(\ref{claim_step_3_cond_1}) and (\ref{claim_step_3_cond_2}) are
satisfied for $i=p'$. We construct $y_i$ and $x_{i-1}$ from $x_i$
recursively beginning on $i=p'$. The arguments are similar to
those used in step $2$.

The differential $d$ on $L^p/L^{p-1}$ restricts to
$$d:(M_p)^0_t\rightarrow (M_p)^0_{t-1}.$$
As $d(x_{p'})=d(x)=0$, for every $\epsilon\in {N\P}_{t-1}$ with
$deg(\epsilon_{t-1})=p$ we can apply the projection
$$\pi_\epsilon: (M_p)^0_{t-1} \twoheadrightarrow F_\epsilon$$ and
obtain $\pi_\epsilon(d(x))=0$. If $deg(\epsilon_0)>p'$ then
$$\pi_\epsilon(d(x))=\sum_{\sigma \in{N\P}_t, d_0(\sigma)=\epsilon}
F(\alpha_1)(x_{\sigma})
$$
and thus
\begin{equation}\label{equ_epsilon_step_3}
0=\sum_{\sigma \in{N\P}_t, d_0(\sigma)=\epsilon}
F(\alpha_1)(x_{\sigma})
\end{equation}
for each $\epsilon\in {N\P}_{t-1}$ with $deg(\epsilon_0)>p'$ and
$deg(\epsilon_{t-1})=p$. Notice that each summand $x_\sigma$ with
$\sigma \in {N\P}_t$, $deg(\sigma_0)=p'$ and $deg(\sigma)=p$
appears in one and just one equation as (\ref{equ_epsilon_step_3})
(take $\epsilon=d_0(\sigma)$). Using now pseudo-injectivity we
build as before $y_\sigma$, $y_\epsilon=\sum_{\sigma} y_\sigma$
and $y=\sum_{\epsilon} y_\epsilon$, where $\epsilon$ runs over
$\epsilon\in {N\P}_{t-1}$ with $deg(\epsilon_0)>p'$ and
$deg(\epsilon_{t-1})=p$, such that
$$
d(y)=x+R
$$
with $R$ living in $\bigoplus_{\sigma \in {N\P}_t,
deg(\sigma_0)=p'-1, deg(\sigma_t)=p} F_{\sigma}\subseteq
(M_p)^0_t$. Call $y_{p'}\definicio y$ and $x_{p'-1}=R$. Then
Equation (\ref{equ_claim_step_3}) is satisfied. Condition
(\ref{claim_step_3_cond_1}) for $i=p'-1$ holds by the construction
of $R$ and condition (\ref{claim_step_3_cond_2}) for $i=p'-1$
holds because $d(x_{p'-1})=d(R)=d(d(y)-x)=d^2(y)-d(x)=0-0=0$ as
$d$ is a differential and $d(x)=0$ by hypothesis. The construction
of $y_i$ and $x_{i-1}$ from $x_i$ is totally analogous to the
construction of $y_{p'}$ and $x_{p'-1}$ from $x_{p'}$ that we have
just made.

After we have built $y_1$ and $x_0$ if we try to build
$y=\sum_\epsilon y_\epsilon$ and $R$ from $x_0$ we find that,
because there are not objects of negative degree (thus if $z\in
Im(i')$ where $deg(i')=0$ then $z=0$), $x_0=0$.
\end{proof}

The following examples come from Example
(\ref{examples_projective}). They show the weaker conditions that
are needed for $\limdir$-acyclicity instead of projectiveness.

\begin{Ex}\label{limacyI_expushout}
For the ``pushout category" $\P$ with shape
$$\xymatrix {a_0 \ar[r]^{f}\ar[d]^{g} & b_1 \\ c_1}$$
a functor $F:\P\rightarrow \Ab$ is $\limdir$-acyclic if $F(f)$ and
$F(g)$ are monomorphisms.

For the ``telescope category" $\P$ with shape
$$ \xymatrix{a_0 \ar[r]^{f_1} & a_1 \ar[r]^{f_2} & a_2 \ar[r]^{f_3} &  a_3 \ar[r]^{f_4} & a_4 ...}$$
a functor $F:\P\rightarrow \Ab$ is $\limdir$-acyclic if
$F(f_{i}\circ f_{i-1} \circ ..\circ f_1)$ is a monomorphism and
$\ker F(f_{i}\circ f_{i-1} \circ ..\circ f_{i-d+1})\subseteq \im
F(f_{i-d})$ for $d=1,2,3,..,i-1$  for each $i=2,3,4,...$

Notice that for this is enough that $F(f_i)$ is a monomorphism for
each $i=1,2,3,..$.
\end{Ex}

\section{Computing higher
limits}\label{limII} Theorem \ref{pro_acyclic_graded} shows that
over a bounded below graded poset pseudo-projectivity is enough
for $\limdir$-acyclicity. But it turns out that
pseudo-projectivity is not necessary for $\limdir$-acyclicity:

\begin{Ex}\label{limacyII_ex1}
For the ``pullback category" $\P$ with shape
$$\xymatrix { & a_0 \ar[d]^{f} \\  b_0\ar[r]^{g} &  c_1}$$
a functor $F:\P\rightarrow \Ab$ is pseudo-projective if
\begin{itemize}
\item $F(f)$ and $F(g)$ are monomorphisms. \item $\im F(f) \cap
\im F(g)=0$.
\end{itemize}
But a straightforward calculus shows $\limdir_i F=0$ for $i\geq 1$
for any $F$.
\end{Ex}

This shows that pseudo-projectivity is not necessary for
$\limdir$-acyclicity. However, we shall see how
pseudo-projectivity allows us to obtain a better knowledge of the
higher limits $\limdir_i F$. We begin with

\begin{Defi}
Let $F:\P\rightarrow \Ab$ be a functor over a graded poset $\P$.
$F':\P\rightarrow \Ab$ is the functor which takes values on
objects
$$F'(i_0)=\bigoplus_{\alpha:i\rightarrow i_0} F(i)$$
for $i_0\in \Ob(\P)$. For $\beta\in \Hom(\P)$,
$\beta:i_1\rightarrow i_0$, $F'(\beta)$ is the only homomorphism
which makes commute the diagram
$$
\xymatrix{ F'(i_1)\ar[r]^{F'(\beta)} & F'(i_0) \\
F(i)_{i\stackrel{\alpha}\rightarrow i_1}\ar@{^{(}->}[u]\ar[r]^{1}&
F(i)_{i\stackrel{\alpha}\rightarrow i_1\stackrel{\beta}\rightarrow
i_0}\ar@{^{(}->}[u]
 }
$$
for each $\alpha:i\rightarrow i_1$.
\end{Defi}

$F'$ is built from $F$ as $Coker'$ was built from $Coker$ in
Section \ref{section_projective}. It mimics the construction of
free objects in $\Ab^\P$. Notice that $Coker_{F'}(i)=F(i)$ for
each $i\in \Ob(\P)$. A nice property of $F'$ is
\begin{Lem}\label{lema F'adjunto}
Let $F:\P\rightarrow \Ab$ be a functor over a graded poset. Then
for each $G\in \Ab^\P$ there is a bijection

$$\xymatrix{\Hom_{\Ab^\P}(F',G)\ar[r]^<<<<{\varphi}_<<<<{\cong} & \prod_{i\in \Ob(\P)} \Hom_{\Ab}(F(i),G(i))}.$$
\end{Lem}
\begin{proof}
$\varphi$ is given by
$$\varphi(\nu:F'\Rightarrow G)_i=(F(i)_{i\stackrel{1_i}\rightarrow i}\hookrightarrow F'(i)\stackrel{\nu_i}\rightarrow G(i)).$$
For a family $\tau=\{\tau_i\}_{i\in \Ob(\P)}\in \prod_{i\in
\Ob(\P)} \Hom_{Set}(F(i),G(i))$ define the natural transformation
$\psi(\tau):F'\Rightarrow G$ on the object $i_0\in \Ob(\P)$ as the
only homomorphism which makes commute the diagram
$$
\xymatrix{ F'(i_0)\ar[r]^{\psi(\tau)_{i_0}} & F(i_0) \\
F(i)_{i\stackrel{\alpha}\rightarrow
i_0}\ar@{^{(}->}[u]\ar[r]^{\tau_i} & G(i)\ar[u]^{G(\alpha)} }
$$
for every $\alpha:i\rightarrow i_0$. Then both compositions
$\varphi\circ \psi$ and  $\psi\circ \varphi$ are the identity.
\end{proof}

Another interesting property of $F'$ is the following

\begin{Lem}\label{F'lim}
Let $F:\P\rightarrow \Ab$ be a functor over a graded poset. Then
$$\limdir F'\cong\bigoplus_{i\in \Ob(\P)} F(i).$$
\end{Lem}
\begin{proof}
It is straightforward using the previous lemma. Notice that the
cone $\eta:F'\Rightarrow \limdir F'$ is given by the homomorphisms
$\eta_{i_0}$ for each $i_0\in \Ob(\P)$ which make commutative the
diagrams
$$
\xymatrix{
 F'(i_0)\ar[r]^>>>{\eta_{i_0}}& \bigoplus_{i\in \Ob(\P)} F(i)\\
 F(i)_{\alpha:i\rightarrow i_0}\ar@{^(->}[u]\ar[r]^{1}&
 F(i)\ar@{^(->}[u]}
$$
for each $\alpha:i\rightarrow i_0$. This description shall be
useful later.
\end{proof}

The main feature of $F'$ we shall use is

\begin{Lem}\label{F'_pseudo}
Let $F:\P\rightarrow \Ab$ be a functor over a graded poset. Then
$F'$ is pseudo-projective.
\end{Lem}
\begin{proof}
Take $k$ objects $i_j\in \Ob(\P)$, arrows $\alpha_j:i_j\rightarrow
i_0$ with $deg(\alpha_j)=d$ and $y_j\in F'(i_j)$ ($j=1,..,k$) such
that
$$\sum_{j=1,..,k}F'(\alpha_j)(y_j)=0 \text{    (in $F'(i_0)$)}.$$
We want that $y_j\in\im_{F'} (i_j)$ for $j=1,..,k$. Write
$y_j=\bigoplus_ {\alpha:i\rightarrow i_j} y_{j,\alpha}$.

Fix $j_0\in\{1,..,k\}$ and consider the projection $p_{j_0}$ from
$F'(i_0)$ onto the direct summand $F(i_{j_0})\hookrightarrow
F'(i_0)$ which corresponds to $\alpha_{j_0}:i_{j_0}\rightarrow
i_0$
$$F'(i_0)\stackrel{p_{j_0}}\twoheadrightarrow
F(i_{j_0}).$$

For any $y=\bigoplus_{\alpha:i\rightarrow i_j} y_\alpha\in
F'(y_j)$,
$$p_{j_0}(F'(\alpha_j)(y))=\sum_{\alpha:i\rightarrow i_j,\alpha_j\circ\alpha=\alpha_{j_0}}
y_\alpha.$$ So, for $y_j=\bigoplus_{\alpha:i\rightarrow i_j}
y_{j,\alpha} \in F'(i_j)$ we have
$$p_{j_0}(F'(\alpha_j)(y_j))=\sum_{\alpha:i\rightarrow i_j,\alpha_j\circ\alpha=\alpha_{j_0}}
y_{j,\alpha}.$$

This last sum runs over $\alpha:i_{j_0}\rightarrow i_j$ such that
the following triangle commutes
$$
\xymatrix{ i_{j_0}\ar[rd]^{\alpha}\ar[rr]^{\alpha_{j_0}} && i_0.\\
& i_j \ar[ru]^{\alpha_j} &}
$$
Because we are in a graded poset and $deg(i_j)=d$ for each
$j=1,..,k$ then the only chance is $i_j=i_{j_0}$ and
$\alpha=1_{i_{j_0}}$. Because the objects $i_1,..,i_k$ are
different this implies that $j=j_0$ too. Thus
\begin{numcases}{p_{j_0}(F'(\alpha_j)(y_j))=}
y_{j_0,1_{i_{j_0}}} & for $j=j_0$ \nonumber \\
0 & for $j\neq j_0$. \nonumber
\end{numcases}
Then
\begin{align}
 0=&p_{j_0}(\sum_{j=1,..,k}F'(\alpha_j)(y_j))\notag\\
  =&p_{j_0}(\sum_{j=1,..,k}F'(\alpha_j)(y_j))\notag\\
  =&\sum_{j=1,..,k}p_{j_0}(F'(\alpha_j)(y_j))=y_{j_0,1_{i_{j_0}}}.\notag
\end{align}

As $j_0$ was arbitrary this means that $y_{j,1_{i_j}}=0$ for each
$j\in\{1,..,k\}$. Now it is clear by the definition of $F'$ that
$y_j\in\im_{F'} (i_j)$ for $j=1,..,k$.

\end{proof}

\begin{Rmk}\label{classical_enough projectives}
The epic natural transformation $G'\Rightarrow F$, where
$G=\Z\circ\U\circ F$ with $\U:\Ab\rightarrow \Set$ the forgetful
functor and $\Z:\Set\rightarrow \Ab$ the free abelian group on a
set, is the usual way to prove that $\Ab^\C$ has enough
projectives for any small category $\C$.
\end{Rmk}

By Lemma \ref{lema F'adjunto} for the family of homomorphisms
$\{F(i)\stackrel{1_{F(i)}}\rightarrow F(i)\}_{i\in \Ob(\P)}$ we
have a natural transformation $\pi:F'\Rightarrow F$. It is clear
that $\pi:F'\Rightarrow F\Rightarrow 0$ is exact in $\Ab^\P$. Thus
we can consider the object-wise kernel of $\pi:F'\Rightarrow F$ to
obtain a short exact sequence of functors
$$0\Rightarrow K_F\Rightarrow F'\stackrel{\pi}\Rightarrow F\Rightarrow
0.$$

If $\P$ is bounded below then the long exact sequence (see Section
\ref{derivedfunctors}) associated to this short exact sequence
gives
\begin{numcases}{\limdir_j F=}
\limdir_{j-1} K_F & $j>1$ \nonumber \\
\ker\{\limdir K_F \rightarrow \limdir F'\}   & $j=1$ \nonumber
\end{numcases}
because $F'$ is $\limdir$-acyclic (it is pseudo-projective by
Lemma \ref{F'_pseudo} and apply Theorem \ref{pro_acyclic_graded}).
So writing $K_0\definicio F$ and $K_1\definicio K_F$, we have
$$\limdir_1 F=\ker\{\limdir K_1\rightarrow \limdir K'_0\}$$
where the map $\limdir K_1\rightarrow \limdir K'_0$ comes from the
long exact sequence of derived functors associated to a short
exact sequence
$$0\Rightarrow K_1\Rightarrow K'_0\Rightarrow K_0\Rightarrow 0.$$
Also we obtained that $\limdir_j F=\limdir_{j-1} K_1$ for $j\geq
2$ and that $\limdir K'_0=\bigoplus_{i\in \Ob(\P)} K_0(i)$. Thus
applying the same machinery to the functor $K_1$ we have a short
exact sequence
$$0\Rightarrow K_2\Rightarrow K'_1\Rightarrow K_1\Rightarrow 0$$
and
$$\limdir_2 F=\limdir_1 K_1=\ker\{\limdir K_2\rightarrow \limdir K'_1\}$$
with $\limdir K'_1 =\bigoplus_{i\in \Ob(\P)} K_1(i)$. Recursively
we obtain short exact sequences
$$0\Rightarrow K_j\Rightarrow K'_{j-1}\Rightarrow K_{j-1}\Rightarrow 0$$
and
$$\limdir_j F=\limdir_{j-1} K_1=\limdir_{j-2} K_2=..=\limdir_1 K_{j-1}=\ker\{\limdir K_j\rightarrow \limdir K'_{j-1}\}$$
for every $j\geq 1$, where $\limdir K'_{j-1}=\bigoplus_{i\in
\Ob(\P)} K_{j-1}(i)$.
\begin{Lem}
Let $\P$ be a bounded below graded poset and $F:\P\rightarrow \Ab$
a functor. Then there are functors $K_j:\P\rightarrow \Ab$ for
$j=0,1,2,...$ with $K_0=F$ and $K_1=\ker(F'\Rightarrow F)$ such
that
$$\limdir_j F=\limdir_{j-1} K_1=\limdir_{j-2} K_2=..=\limdir_1 K_{j-1}=\ker\{\limdir K_j\rightarrow \limdir K'_{j-1}\}$$
for each $j=0,1,2,...$.
\end{Lem}

The values $F'(i_0)$ and $K_F(i_0)$ can be very big because they
contain a copy of $F(i)$ for each $i\rightarrow i_0$. This can be
improved considering the functor $\coker:\P\rightarrow \Ab$ in
Section \ref{section_projective}. Suppose that for every $i_0$
there is a section $s_{i_0}:\coker(i_0)\rightarrow F(i_0)$ to the
projection $F(i_0)\twoheadrightarrow \coker(i_0)$ (for example if
$\coker(i_0)$ is free for each $i_0$ or if $F$ is an epic
functor). Then by Lemma \ref{lema F'adjunto} there is a natural
transformation $\coker'\Rightarrow F$. If $\P$ is bounded below
then it is easy to see by induction on the degree of objects that
this natural transformation is object-wise surjective. Notice that
$\coker'(i_0)$ is, in general, smaller than $F'(i_0)$.

\begin{Lem}\label{compute_limi_with_coker}
Let $\P$ be a bounded below graded poset and let $F:\P\rightarrow
\Ab$ be an epic functor. Then there is a short exact sequence of
functors
$$0\Rightarrow K\Rightarrow \coker'\Rightarrow F\Rightarrow 0$$
where $\coker'$ is $\limdir$-acyclic.
\end{Lem}

We finish this section with some examples of $\limdir$-acyclic
functors:

\begin{Ex}
Let $\P$ be a graded poset with initial object $i_0$. Then $\P$ is
contractible and thus $H_i(\P;M)=0$ for $i\geq 1$ and any trivial
coefficients $M\in \Ab$. We can prove this by taking the functor
$F_M:\P\rightarrow \Ab$ which takes the value $M$ on $i_0$ and $0$
otherwise. Then $(F_M)':\P\rightarrow \Ab$ is the functor of
constant value $M$ and thus $H_i(\P;M)=\limdir_i (F_M)'$ for
$i\geq 0$. These higher limits vanish for $i\geq 1$ because
$(F_M)'$ is pseudo-projective (Lemma \ref{F'_pseudo}). Finally,
$H_0(\P;M)=\limdir_0 (F_M)'=M$ by Lemma \ref{F'lim}.
\end{Ex}

\begin{Defi}
Let $\P$ be graded poset. We say $\P$ is a \emph{rooted tree} if
$\P$ is a tree and $\P$ has a initial object.
\end{Defi}

\begin{Ex}
The category with shape
$$
\xymatrix{a\ar[r]&c\ar[r]&e\\
b\ar[r]&d\ar[ru]\ar[r]&f\ar[r]&g }
$$
is a tree but not a rooted tree. The category with shape
$$
\xymatrix{&&d\\
a\ar[r]\ar[rd]&b\ar[ru]\ar[r]&e\\
&c&&
 }
$$
is a rooted tree.
\end{Ex}

\begin{Cor}\label{rootedtreemonic_acyclic}
Let $\P$ be a rooted tree. Then any monic functor $F:\P\rightarrow
\Ab$ is $\limdir$-acyclic.
\end{Cor}
\begin{proof}
Just check that a monic functor $F:\P\rightarrow \Ab$ over a
rooted tree $\P$ is pseudo-projective.
\end{proof}

\begin{Ex}\label{ex_pushout_telescope limdir_monic acyclic}
Any monic functor over the push-out category or the ``telescope
category" is $\limdir$-acyclic.
\end{Ex}






\section{$Lim_1$ as a flow problem.}\label{section_lim1_flow} In
this section we give an interpretation of the first derived limit
$\limdir_1 F$ in terms of flow problems. While in classical flows
on directed graphs (see \cite[III.1]{bela}) a non-negative integer
is associated to each edge $i_0\rightarrow i_1$, we associate a
value from the abelian group $F(i_0)$. Here, $F$ is a functor
$F:\P\rightarrow \Ab$ and $\P$ is a graded poset. Notice that, by
Section \ref{limII}, all the higher limits $\limdir_i F$ can be
reduced to first derived functors $\limdir_1 K_{i-1}$ for suitable
$K_{i-1}$.

Recall Definition \ref{Defi:limi}. It is straightforward that
$\limdir_1 F$ equals a quotient $M/N$ where $N\subseteq M$ are
abelian subgroups of
$$\bigoplus_{i_0\stackrel{\alpha}\rightarrow i_1 \in
{N\P}_1} F(i_0).$$

Notice that ${N\P}_1=\Hom(\P)$ and that an element
$x=\{x_\alpha\}_{\alpha\in \Hom(\P)}$ of
$\bigoplus_{i_0\stackrel{\alpha}\rightarrow i_1 \in {N\P}_1}
F(i_0)$ belongs to $M$ just in case
$$\sum_{i\stackrel{\alpha}\rightarrow i_0}
F(\alpha)(x_\alpha)=\sum_{i_0\stackrel{\alpha}\rightarrow i}
x_\alpha$$ for every $i_0\in \Ob(\P)$. The subgroup $N$ is the
image of the differential
$$C_2(\P,F)\stackrel{d}\rightarrow C_1(\P,F)$$ and it is generated
by
$$x_\alpha \oplus (-x)_{\beta\circ\alpha} \oplus (F(\alpha)(x))_\beta$$
where $i_0\stackrel{\alpha}\rightarrow i_1$ and
$i_1\stackrel{\beta}\rightarrow i_2$ are any two composable
morphisms and $x\in F(i_0)$.

\begin{Defi}
Let $\P$ be a graded poset and let $F:\P\rightarrow \Ab$ be a
functor. A \emph{generalized flow} is an element
$x=\{x_\alpha\}_{\alpha\in \Hom(\P)}$ with finitely many terms
different from zero such that
$$\sum_{i\stackrel{\alpha}\rightarrow i_0}
F(\alpha)(x_\alpha)=\sum_{i_0\stackrel{\alpha}\rightarrow i}
x_\alpha$$ for every $i_0\in \Ob(\P)$. A minimal trivial flow is a
generalized flow
$$x_\alpha \oplus (-x)_{\beta\circ\alpha} \oplus (F(\alpha)(x))_\beta$$
where $i_0\stackrel{\alpha}\rightarrow i_1$ and
$i_1\stackrel{\beta}\rightarrow i_2$ are any two composable
morphisms and $x\in F(i_0)$.
\end{Defi}

With these definitions it is clear that
\begin{Lem}
$\limdir_1 F=0$ if and only if any generalized flow can be written as sum of minimal
trivial flows.
\end{Lem}

There is also an easy geometrical interpretation of $M$ on the
directed graph associated to $\P$ (see \ref{graded_poset_graphs}):
because in the graded poset $\P$ any morphism factors as
composition of morphisms of degree $1$  , then in the class
$\overline{x}\in \limdir_1 F=M/N$ of a generalized flow
$x=\{x_\alpha\}_{\alpha\in \Hom(\P)}$ there is a representantive
$x'=\{x'_\alpha\}_{\alpha\in \Hom_1(\P)}$ which takes nonzero
values just in degree $1$ morphisms. To obtain $x'$ is enough to
use the cycles
$$(-x_\gamma)_{\gamma_1} \oplus (x_\gamma)_{\gamma} \oplus (F(\alpha)(-x_\gamma))_{\gamma_2}$$
where $\gamma=\gamma_2\circ\gamma_1$ is a morphism of degree
greater than $1$, arguing by induction on
$\max\{deg(\alpha)|x_\alpha\neq 0 \}$ (morphisms of degree $0$ do
not appear as we assume that we are using normalized chain
complexes \ref{normalization}).

\begin{Defi}
Let $\P$ be a graded poset and let $F:\P\rightarrow \Ab$ be a
functor. A \emph{flow} is an element $x=\{x_\alpha\}_{\alpha\in
\Hom_1(\P)}$ with finitely many terms different from zero such
that
$$\sum_{i\stackrel{\alpha}\rightarrow i_0,deg(\alpha)=1}
F(\alpha)(x_\alpha)=\sum_{i_0\stackrel{\alpha}\rightarrow
i,deg(\alpha)=1} x_\alpha$$ for every $i_0\in \Ob(\P)$.
\end{Defi}

Again we have
\begin{Lem}
$\limdir_1 F=0$ if and only if any flow can be written as sum of minimal trivial flows.
\end{Lem}

Consider the directed graph associated to $\P$ (see
\ref{graded_poset_graphs}). Then a flow corresponds to a choice of
a value $x_\alpha\in F(i_0)$ for each edge
$i_0\stackrel{\alpha}\rightarrow i_1$ of this graph, i.e., for
each morphism of degree $1$ of $\P$:
$$
\xymatrix{i_0 \ar[r]^{x_\alpha} & i_1 }.
$$

In order for $x$ to be a flow, for each $i_0\in \Ob(\P)$ we must
have the equality
$$\sum_{i\stackrel{\alpha}\rightarrow i_0, deg(\alpha)=1}
F(\alpha)(x_\alpha)=\sum_{i_0\stackrel{\alpha}\rightarrow i,
deg(\alpha)=1} x_\alpha.$$ For example, for a vertex $i_0$ as
$$
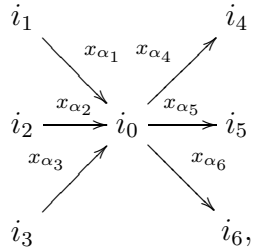
\xymatrix{i_1 \ar[rd]^{x_{\alpha_1}} &  & i_4\\
 i_2 \ar[r]^{x_{\alpha_2}} & i_0 \ar[ru]^{x_{\alpha_4}}\ar[r]^{x_{\alpha_5}}\ar[rd]^{x_{\alpha_6}} & i_5\\
  i_3 \ar[ru]^{x_{\alpha_3}} &  & i_6, }
$$
where $\Hom_1(\P,i_0)=\{\alpha_1,\alpha_2,\alpha_3\}$ and
$\Hom_1(i_0,\P)=\{\alpha_4,\alpha_5,\alpha_6\}$, this condition
reads
$$
F(\alpha_1)(x_{\alpha_1})+F(\alpha_1)(x_{\alpha_2})+F(\alpha_3)(x_{\alpha_3})=
x_{\alpha_4}+x_{\alpha_5}+x_{\alpha_6}.
$$

A minimal trivial flow can be represented as the generalized flow
$$
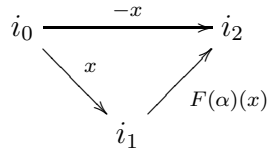
\xymatrix{i_0\ar[rd]^{x}\ar[rr]^{-x}& &i_2\\
&i_1\ar[ru]_{F(\alpha)(x)} &
}
$$
(but notice that the arrows are not edges of the directed graph
associated to $\P$ in general).

\begin{Ex}
Consider the graded poset $\P$ with shape (it was considered
before on Example \ref{ss_example_cycle})
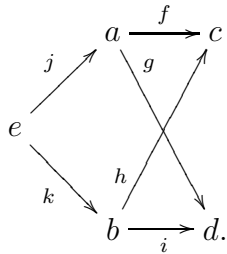
$$ \xymatrix{
 & a \ar[rdd]^<<<<{g} \ar[r]^{f} & c \\
e \ar[ru]^{j}\ar[rd]_{k} &         & \\
 & b \ar[ruu]^<<<<{h} \ar[r]_{i} & d. \\
}$$ This representation corresponds exactly with the directed
graph associated to $\P$, because only the degree $1$ morphisms
are displayed. Consider the functor \linebreak $F:\P\rightarrow
\Ab$ with constant value $\Z$. Then the following is a flow
$$ 
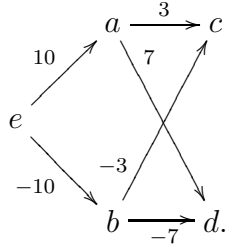
\xymatrix{
 & a \ar[rdd]^<<<<{7} \ar[r]^{3} & c \\
e \ar[ru]^{10}\ar[rd]_{-10} &         & \\
 & b \ar[ruu]^<<<<{-3} \ar[r]_{-7} & d. \\
}$$ This flow can be written as sum of the following minimal trivial flows:
$$ 
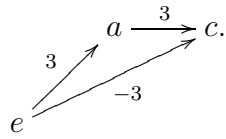
\xymatrix{
 & a  \ar[r]^{3} & c. \\
e \ar[ru]^{3}\ar[rru]_{-3} &         &}
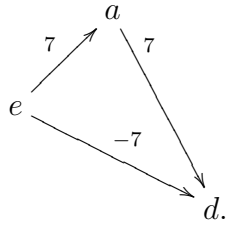
$$
$$ \xymatrix{
 & a  \ar[rdd]^<<<<{7}\\
e \ar[ru]^{7}\ar[rrd]^{-7} &         & \\
 & & d. \\
}
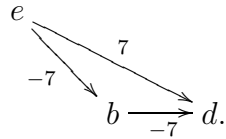
$$
$$ \xymatrix{
e \ar[rd]_{-7}\ar[rrd]^7 &         & \\
 & b \ar[r]_{-7} & d. \\
}
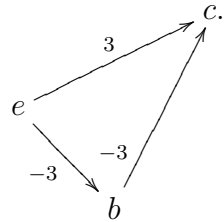
$$
$$ \xymatrix{
 &  & c. \\
e \ar[rru]^{3}\ar[rd]_{-3} &         & \\
 & b \ar[ruu]^<<<<{-3} \\
}$$
\end{Ex}

This naive interpretation of $\limdir_1 F$ has some consequences,
which will be useful later. Recall (\ref{section graded_posets})
that $(p_0\downarrow \P)_*$ is the full subcategories of $\P$ with
objects $\{p| \exists p_0\rightarrow p,p\neq p_0\}$.

\begin{Defi}\label{defi_core}
Let $\P$ be a graded poset. Write $\P'$ for the full subcategory
of $\P$ with objects all but those objects $i_0$ such that there
are no arriving arrows to $i_0$ and such that $(i_0\downarrow
\P)_*$ is empty or connected. Write $\P^n=(\P^{n-1})'$ for $n\geq
1$ and $\P^0=\P$. We have inclusions
$$\P^n\subseteq \P^{n-1}\subseteq..\subseteq \P^1=\P'\subseteq \P^0=\P$$

Then define the \emph{core} of $\P$, denoted $core(\P)$, as the
inverse limit $\liminv_n \P^n$ in the complete category of small
categories.
\end{Defi}

Notice that if $\P$ is finite then we obtain $core(\P)$ from $\P$
after a finite number of times of applications of the operator
$(\cdot)'$.

\begin{Ex}\label{ex_core_some_category}
Consider the graded poset $\P$ with shape
$$ 
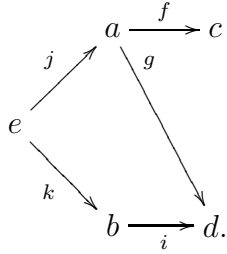
\xymatrix{
 & a \ar[rdd]^<<<<{g} \ar[r]^{f} & c \\
e \ar[ru]^{j}\ar[rd]_{k} &         & \\
 & b  \ar[r]_{i} & d. \\
}$$ Then $\P^1=\P'$ is the category
$$ 
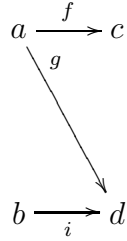
\xymatrix{
 & a \ar[rdd]^<<<<{g} \ar[r]^{f} & c \\
 &         & \\
 & b  \ar[r]_{i} & d \\
}$$ and $\P^2=\P''=core(\P)$ is
$$ \xymatrix{
 & a \ar[rdd]^<<<<{g} \ar[r]^{f} & c \\
 &         & \\
 &  & d. \\
}$$
\end{Ex}

\begin{Lem}\label{lem_lim1_core}
Let $\P$ be a finite graded poset and let $F:\P\rightarrow \Ab$ be
a functor. Then any class $\overline{x}\in \limdir_1 F$ has a
representative $x'=\{x'_\alpha\}_{\alpha\in \Hom(core(\P))}$ which
takes nonzero values just in morphisms of the core of $\P$.
\end{Lem}
\begin{proof}
By the definition of the core of $\P$ it is enough to prove that
any generalized flow $x$ is in the same class that a generalized
flow $x'=\{x'_\alpha\}_{\alpha\in \P'}$ which takes nonzero values
just in morphisms of $\P'$.

Take any generalized flow $x$ over $\P$ and its class
$\overline{x}\in \limdir_1 F$. By the above remarks we can suppose
that $x$ is in fact a flow, i.e., takes nonzero values just in
degree $1$ morphisms.

Consider any object $i_0$ which is in $\P$ but not in $\P'$. Then
there are no arriving arrows to $i_0$ and $(i_0\downarrow \P)_*$
is empty or connected. If it is empty then there is nothing to do.
Thus, assume  that $(i_0\downarrow \P)_*$ is non-empty and
connected. The equation for $x$ on $i_0$ becomes
\begin{equation}\label{lema_core_equ_x}
0=\sum_{i_0\stackrel{\alpha}\rightarrow i, deg(\alpha)=1}
x_\alpha.
\end{equation}
Fix $i_0\stackrel{\alpha}\rightarrow i_*$ with $deg(\alpha)=1$.
Now, take any other arrow ($\alpha\neq\beta$)
$i_0\stackrel{\beta}\rightarrow j_*$ of degree $1$ (if it does not
exist such a $\beta$ then we are done). Recall that
$(i_0\downarrow \P)_*$ is connected. Then there exists a zigzag of
morphisms connecting $i_*$ and $j_*$. Because $\P$ is a graded
poset we can assume that the zig-zag has the following shape:
$$
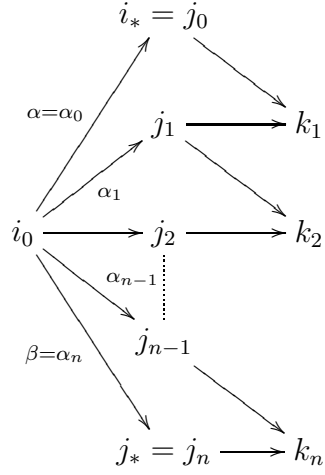
\xymatrix{  & i_*=j_0\ar[rd]& \\
& j_1\ar[r]\ar[rd]&k_1\\
i_0\ar[ruu]^{\alpha=\alpha_0}\ar[ru]_{\alpha_1}\ar[rd]^>>>>>{\alpha_{n-1}}\ar[rdd]_{\beta=\alpha_n}\ar[r]&j_2\ar@{..}[d]\ar[r]&k_2\\
&j_{n-1}\ar[rd]\\
&j_*=j_n\ar[r]&k_n}
$$
where $deg(\alpha_i)=1$ for $i=0,...,n$. For each $i=0,...,n-1$
consider the following diamond with commutative triangles
$$
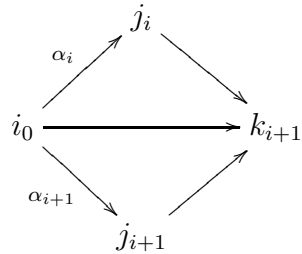
\xymatrix{  & j_i\ar[rd]& \\
i_0\ar[ru]^{\alpha_i}\ar[rd]_{\alpha_{i+1}}\ar[rr]&&k_{i+1}\\
&j_{i+1}\ar[ru]}
$$
and the two minimal trivial flows
$$
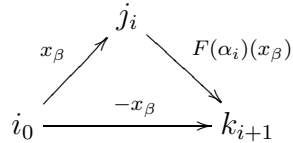
\xymatrix{ &j_i\ar[rd]^{F(\alpha_i)(x_\beta)}\\
i_0\ar[ru]^{x_\beta}\ar[rr]^{-x_\beta}&&k_{i+1} }
$$
and
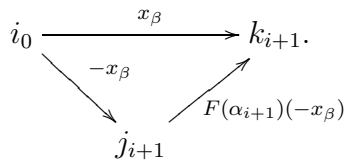
$$
\xymatrix{ i_0\ar[rd]^{-x_\beta}\ar[rr]^{x_\beta}&&k_{i+1}.\\
&j_{i+1}\ar[ru]_{F(\alpha_{i+1})(-x_\beta)}&}
$$
If we add up these two minimal trivial flows for all $i=0,...,n-1$
we obtain a generalized flow $y_{j_*}$ with $\overline{y_{j_*}}=0$
and such that, for $\gamma:i_0\rightarrow i$,
\begin{numcases}{{y_{j_*}}_\gamma=}
 x_\beta & if $\gamma=\alpha$ \nonumber \\
-x_\beta & if $\gamma=\beta$ \nonumber \\
 0 & otherwise. \nonumber
\end{numcases}

Thus, $x'=x+y_{j_*}$ is a representative for $x$ which verifies,
for $\gamma:i_0\rightarrow i$,
\begin{numcases}{x'_\gamma=}
x_\alpha+x_\beta & if $\gamma=\alpha$ \nonumber \\
 0 & if $\gamma=\beta$ \nonumber\\
 x_\gamma & otherwise. \nonumber
\end{numcases}

Notice that, by construction, if $x'_\gamma\neq 0$ then either
$x_\gamma\neq 0$ or $\gamma$ has origin in $i$ and $i\in
\Ob(\P')$. Doing the same construction for each arrow
$i_0\rightarrow j$ different from $\alpha$ (and of degree $1$)
with $x_{i_0\rightarrow j}\neq 0$ (recall that there is a finite
number of these arrows) we obtain a representative $x'$ for $x$
such that, for $\gamma:i_0\rightarrow i$,
\begin{numcases}{x'_\gamma=}
x_\alpha+\sum_{\beta:i_0\rightarrow j,deg(\beta)=1,\beta\neq\alpha} x_\beta & if $\gamma=\alpha$ \nonumber \\
 0 & otherwise \nonumber
\end{numcases}

Then, by Equation (\ref{lema_core_equ_x}), $x'_\alpha=0$ and
$x'_\gamma=0$ for each arrow $\gamma$ with origin in $i_0$. Again
by construction if $x'_\gamma\neq 0$ then either $x_\gamma\neq 0$
either $\gamma$ has origin in $i$ and $i\in \Ob(\P')$. As it was
described before the lemma, there is a representative $x''$ for
$\overline{x'}=\overline{x}$ which takes non-zero values just in
morphisms of degree $1$. Moreover, if $x''_\gamma\neq 0$ then
either $x_\gamma\neq 0$ either $\gamma$ has origin in $i$ and
$i\in \Ob(\P')$.

Then, we repeat the process (but with $x''$ instead of $x$) for
any object $i_0$ which is in $\P$ but no in $\P'$ and for which
there is a morphism $\alpha$ with origin in $i_0$ and with
$x_\alpha\neq 0$. There is a finite number of such objects.

Thus, finally, we obtain a representative for $\overline{x}$ which
takes values different from zero on morphism which do not begin on
objects from $\Ob(\P)\setminus\Ob(\P')$, i.e., on morphisms from
$\Hom(\P')$.

Because $\P$ is finite by hypothesis, then $core(\P)$ is reached
in a finite number of computations $\P,\P',(\P')',...,core(\P)$
and the lemma is proven.
\end{proof}

Next, we present some conditions on vanishing of higher limits
derived from the previous lemma:

\begin{Cor}\label{cor_core_empty_limdir_acyclic}
Let $\P$ be a finite graded poset with $core(\P)=\emptyset$. Then
every functor $F:\P\rightarrow \Ab$ is $\limdir$-acyclic.
\end{Cor}
\begin{proof}
That $\limdir_1 F=0$ results directly from the lemma above. The
limit $\limdir_i F$ with $i\geq 2$ is given by $\limdir_1 K_{i-1}$
where $K_{i-1}:\P\rightarrow \Ab$ is a functor (see section
\ref{limII}).
\end{proof}

\begin{Cor}\label{cor_core_determines_limi}
Let $\P$ be a finite graded poset such that $core(\P)$ is a tree.
If $F:\P\rightarrow \Ab$ is a functor with $F|_{core(\P)}$ monic
then $F$ is $\limdir$-acyclic.
\end{Cor}
\begin{proof}
For $\limdir_1 F$ notice that the class $\overline{x}\in \limdir_1
F$ of any generalized flow $x$ has as representative a flow $x'$
living in $core(\P)$. Then the equations for $x'$ to be a flow
imply that $x'\equiv 0$ as $core(\P)$ is a tree and $F$ restricted
to $core(\P)$ is monic.

The limit $\limdir_i F$ with $i\geq 2$ is given by $\limdir_1
K_{i-1}$ where $K_{i-1}:\P\rightarrow \Ab$ is a monic functor (see
Section \ref{limII}). Then use the same argument.
\end{proof}

This implies, in particular, that (cf. Corollary
\ref{rootedtreemonic_acyclic}):

\begin{Cor}\label{Ex_monic on tree acyclic direct}
Let $\P$ be a finite tree. If $F:\P\rightarrow \Ab$ is a monic
functor then $F$ is $\limdir$-acyclic.
\end{Cor}
\begin{proof}
If $\P$ is a tree then $core(\P)$ is a subcategory of $\P$ and
thus it is a tree too. Then use the previous corollary.
\end{proof}





\chapter{Higher inverse
limits}\label{section_on inverse limit}

\section{Injective objects in $\Ab^{\mathcal{I}}$.}\label{section_injective}
Consider the abelian category $\Ab^\P$ for some graded poset $\P$.
In this section we shall determine the injective objects in
$\Ab^\P$ following dual arguments to those of Chapter
\ref{section_projective}. Recall that in $\Ab$ the injective
objects are well known, and are exactly the direct sums of $\Q$
and $\Z[p^{\infty}]$ for various primes $p$. Along the rest of the
section, $\P$ denotes a graded poset.

Suppose $F\in \Ab^\P$ is injective. How does $F$ look? We show
that the intersection of the kernels of the non-identity morphisms
with source $i_0$ is an injective abelian group. To prove it,
write

\begin{Defi}\label{defi_ker}
$\ker(i_0)=\bigcap_{i_0\stackrel{\alpha}\rightarrow i,\alpha\neq
1_{i_0}} \ker F(\alpha)$ (or $\ker(i_0)=F(i_0)$ if the index set
of the intersection is empty) and $\coim(i_0)=F(i_0)/\ker(i_0)$.
\end{Defi}

For any diagram in $\Ab$ as the following

$$ \xymatrix{& \ker(i_0)  &\\
A_0  \ar@{-->}[ur]^{\rho_0} & B_0
\ar[l]^{\lambda_0}\ar[u]^{\sigma_0} & 0\ar[l] }
$$
we want to find $\rho_0$ that makes it commutative. Consider the
atomic functors $A,B:\P\rightarrow \Ab$ which take the values on
objects
\begin{numcases}{A(i)=}
A_0 & for $i=i_0$ \nonumber \\
0 & for $i\neq i_0$ \nonumber
\end{numcases}
\begin{numcases}{B(i)=}
B_0 & for $i=i_0$ \nonumber \\
0 & for $i\neq i_0$ \nonumber
\end{numcases}
and on morphisms
\begin{numcases}{A(\alpha)=}
1_{A_0} & for $\alpha=1_{i_0}$ \nonumber \\
0 & for $\alpha\neq 1_{i_0}$ \nonumber
\end{numcases}
\begin{numcases}{B(\alpha)=}
1_{B_0} & for $\alpha=1_{i_0}$ \nonumber \\
0 & for $\alpha\neq 1_{i_0}$ \nonumber
\end{numcases}
and the natural transformations $\sigma:B\Rightarrow F$ and
$\lambda:A\Rightarrow B$ given by
\begin{numcases}{\sigma(i)=}
j\circ \sigma_0 & for $i=i_0$ \nonumber \\
0 & for $i\neq i_0$ \nonumber
\end{numcases}
\begin{numcases}{\lambda(i)=}
\lambda_0 & for $i=i_0$ \nonumber \\
0 & for $i\neq i_0$ \nonumber
\end{numcases}
where $j$ is the inclusion $\ker(i_0)\hookrightarrow F(i_0)$. Then
$0\Rightarrow B\stackrel{\lambda}\Rightarrow A$ is exact as
$0\rightarrow B_0\stackrel{\lambda_0}\rightarrow A_0$ is so. It is
straightforward that $\lambda$ is a natural transformation. The
key point in checking that $\sigma$ is a natural transformation is
that for $\alpha:i_0\rightarrow i_1$, $\alpha\neq 1_{i_0}$ the
diagram
$$
\xymatrix{ F(i_0)\ar[r]^{F(\alpha)} &F(i_1)\\
           B_0\ar[u]^{j\circ \sigma_0}\ar[r]^{0}& B(i_1)\ar[u]^{\sigma(i_1)} }
$$
must commute. And it does because $F(\alpha)\circ j=0$ for every
$\alpha\neq 1_{i_0}$.

So, as $F$ is injective, this data gives a natural transformation
$\rho$ which makes commutative the diagram of natural
transformations
$$
\xymatrix{
             & F  &\\
A  \ar@{==>}[ru]^\rho& B \ar@{=>}[l]^\lambda \ar@{=>}[u]^\sigma  &
0\ar@{=>}[l] }
$$
which restricts over $i_0$ to
$$ \xymatrix{
& F(i_0)\\
& \ker(i_0)  \ar[u]^{j}   &\\
A_0  \ar[uur]^{\rho(i_0)}\ar@{-->}[ur]^{\rho_0}& B_0
\ar[l]^{\lambda_0}\ar[u]^{\sigma_0} & \ar[l]0. }
$$
Then $\rho_0$ exists if and only if $\im(\rho(i_0))\subseteq
\ker(i_0)$. To check that this condition holds take
$\alpha:i_0\rightarrow i_1$ with $\alpha\neq 1_{i_0}$ and $a\in
A_0$. Then
$$
F(\alpha)(\rho(i_0)(a))=\rho(i_1)(A(\alpha)(a))=\rho(i_1)(0)=0.
$$
We have just proven
\begin{Lem}\label{lem_proj_pseudo_graded_ker_injective}
Let $F:\P\rightarrow \Ab$ be an injective functor over a graded
poset $\P$. Then $\ker(i_0)$ is injective for every object $i_0\in
\Ob(\P)$.
\end{Lem}
This means that we can write
$$ F(i_0)=\ker(i_0)\oplus\coim(i_0)$$
with $\ker(i_0)$ injective for every $i_0\in \Ob(\P)$, and also
that
\begin{Ex}
For the category $\P$ with shape
$$\cdot \rightarrow \cdot $$
the functor $F:\P\rightarrow \Ab$ with values
$$\Q \twoheadrightarrow \Q/\Z=\bigoplus_{\textit{$p$ prime}}\Z[p^\infty],$$ is not injective as $\ker$ on the left object equals
the non-injective abelian group $\Z$.
\end{Ex}

Now that we know a little about the values that an injective
functor $F:\P\rightarrow \Ab$ takes on objects we can wonder about
the values $F(\alpha)$ for $\alpha\in \Hom(\P)$. Do they have any
special property? Recall that a feature of graded posets is that
there is at most one arrow between any two objects, and also that
\begin{Rmk}\label{ker_on_graded}
If $\P$ is graded then for any $i_0\in \Ob(\P)$
$$\ker(i_0)=\bigcap_{i_0\stackrel{\alpha}\rightarrow i,deg(\alpha)=1} \ker F(\alpha) $$
because every morphism factors as composition of morphisms of
degree $1$.
\end{Rmk}

We prove that the following property holds for $F$:

\begin{Defi}\label{property_pseudo-injective}
Let $F:\P\rightarrow \Ab$ be a functor over a graded poset $\P$
with degree function $deg$. Fix an integer $d\geq 0$. If for any
$i_0\in \Ob(\P)$, different objects $\{i_j\}_{j\in J}\subseteq
\Ob(\P)$, arrows $\alpha_j:i_0\rightarrow i_j$ with
$deg(\alpha_j)=d$ and elements $x_j\in \ker(i_j)$ for each $j\in
J$, there is $y\in F(i_0)$ with
$$F(\alpha_j)(y)=x_j$$
for each $j\in J$, we call $F$ \emph{$d$-pseudo-injective}. If $F$
is $d$-pseudo-injective for each $d\geq 0$ we call $F$
pseudo-injective.
\end{Defi}

\begin{Rmk}\label{pseudo_surjective_is_surjective}
In case $J=\{1\}$ and $\ker(i_{1})=F(i_{1})$ we are claiming that
the homomorphism $F(\alpha_{1}):F(i_0)\rightarrow F(i_{1})$ is
surjective. Notice that any functor is $0$-pseudo-injective as the
identity is an epimorphism.
\end{Rmk}

Before proving that injective functors $F$ over a graded poset
verify this property we define two functors $\ker$ and $\ker'$ and
natural transformations $\sigma$ and $\lambda$ that fit in the
diagram
$$
\xymatrix{
             & F &\\
\ker'  & \ker \ar@{=>}[l]^\lambda \ar@{=>}[u]^{\sigma}  &
\ar@{=>}[l]0 }
$$
for any functor $F:\P\rightarrow \Ab$ with $\P$ a graded poset. We
begin defining $\ker$. Because for every $\alpha:i_1\rightarrow
i_0$ holds that $F(\alpha)(\ker(i_1))\leq \ker(i_0)$ we can factor
$F(\alpha)$ as in the diagram
$$
\xymatrix{ F(i_1)\ar[rr]^{F(\alpha)} && F(i_0)\\
\ker(i_1)\ar@{^(->}[u]\ar[rr]^{F(\alpha)|_{\ker(i_1)}} &&
\ker(i_0). \ar@{^(->}[u]}
$$
In fact, if $\alpha\neq 1_{i_1}$, then
$F(\alpha)|_{\ker(i_1)}\equiv 0$ by definition. Because the
identity $1_{i_0}$ cannot be factorized (by non-identity
morphisms) in a graded poset then we have a functor $\ker$ with
value $\ker(i)$ on the object $i$ of $\P$ and which maps the
non-identity morphisms to zero. $\ker$ is a kind of ``discrete"
functor. Also it is clear that exists a natural transformation
$\sigma:\ker\Rightarrow F$ with $\sigma(i)$ the inclusion
$\ker(i)\hookrightarrow F(i)$.

Now we define $\ker'$ from $\ker$ in a dual way as free diagrams
are constructed. Let $\ker'$ be defined on objects by
$$\ker'(i_0)=\bigoplus_{\alpha:i_0\rightarrow i} \ker(i).$$
For $\beta\in \Hom(\P)$, $\beta:i_1\rightarrow i_0$,
$\ker'(\beta)$ is the only homomorphism which makes commute the
diagram
$$
\xymatrix{ \ker'(i_1)\ar[r]^{\ker'(\beta)} & \ker'(i_0) \\
\ker(i)_{i_1\rightarrow i}\ar@{^{(}->}[u]\ar[r]^{1}&
\ker(i)_{i_0\rightarrow i}\ar@{^{(}->}[u]
 }
$$
for each $\alpha:i_1\rightarrow i$ that factors through $\beta$
$$
\xymatrix{ &i_0\ar[rd]& \\
i_1\ar[ru]^\beta \ar[rr]&&i }
$$
and the diagram
$$
\xymatrix{ \ker'(i_1)\ar[r]^{\ker'(\beta)} & \ker'(i_0) \\
\ker(i)_{i_1\rightarrow i}\ar@{^{(}->}[u]\ar[ru]^{0} }
$$
for each $\alpha:i_1\rightarrow i$ that does not factor through $\beta$.

Then there exists a candidate to natural transformation
$\lambda:\ker\Rightarrow \ker'$ whose value $\lambda(i)$ is the
inclusion of $\ker(i)$ into the direct summand
$\ker(i)_{1_i:i\rightarrow i}$ corresponding to the identity on
$i$. $\lambda$ is a natural transformation if for every
$\beta:i_1\rightarrow i_0$ with $i_1\neq i_0$ the following
diagram is commutative
$$
\xymatrix{\ker'(i_1)\ar[r]^{\ker'(\beta)} & \ker'(i_0)\\
\ker(i_1)\ar@{^(->}[u]^{\lambda(i_1)}\ar[r]^{0} &
\ker(i_0)\ar@{^(->}[u]^{\lambda(i_0)}.}
$$
It is straightforward that this square commutes because
$1_{i_1}:i_1\rightarrow i_1$ cannot be factored through $\beta$.

Now we have the commutative triangle
$$
\xymatrix{
             & F &\\
\ker' \ar@{==>}[ru]^\rho & \ker \ar@{=>}[l]^\lambda
\ar@{=>}[u]^{\sigma} & \ar@{=>}[l]0 }
$$
where the natural transformation $\rho$ exists because $F$ is
injective. To prove that $F$ is pseudo-injective take $i_0\in
\Ob(\P)$, different objects $\{i_j\}_{j\in J}\subseteq \Ob(\P)$,
arrows $\alpha_j:i_0\rightarrow i_j$ with $deg(\alpha_j)=d$ and
elements $x_j\in \ker(i_j)$ for each $j\in J$. To visualize what
is going on consider the diagram above near $i_0$ in case
$J=\{1,2\}$
\begin{footnotesize}
$$
\xymatrix{
&&& & F(i_0)\ar[ld]_{F(\alpha_1)}\ar[rd]^{F(\alpha_2)}& \\
&&&F(i_1)&  & F(i_2) &\\
&\ker'(i_0)\ar@{-->}[uurrr]\ar[rd] ^{\ker'(\alpha_2)}\ar[ld]
_{\ker'(\alpha_1)} & &&
 \ker(i_0)\ar@{_(->}[uu]\ar@{_(->}[lll]^{\lambda(i_0)}\ar[rd]^{0}\ar[ld]^{0}& \\
\ker'(i_1)\ar@{-->}[uurrr]&& \ker'(i_2)\ar@<4pt>@{-->}[rrruu]
&\ker(i_1)\ar@{_(->}[uu]& &\ker(i_2).\ar@{_(->}[uu] & }
$$
\end{footnotesize}

Consider the element $x\in
\ker'(i_0)=\bigoplus_{\alpha:i_0\rightarrow i} \ker(i)$ given in
components by
\begin{numcases}{x_\alpha=}
x_j & if $\alpha=\alpha_j:i_0\rightarrow i_j$ for some $j\in J$ \nonumber \\
0 &   otherwise. \nonumber
\end{numcases}
Fix $j_0\in J$ and consider the map
$\ker'(\alpha_{j_0}):\ker'(i_0)\rightarrow \ker'(i_{j_0})$.
Because all the maps $\{\alpha_j\}_{j\in J}$ have degree $d$ then
we have a commutative triangle
$$
\xymatrix{ &i_{j_0}\ar[rd]^{\beta}& \\
i_0\ar[ru]^{\alpha_{j_0}} \ar[rr]^{\alpha_j}&&i_j }
$$
only in case $j=j_0$ and $\beta=1_{j_0}$. This means that
$\ker'(\alpha_{j_0})$ maps $x\in \ker'(i_0)$ to $x_{j_0}\in
\ker(i_{j_0})_{1_{j_0}:i_{j_0}\rightarrow i_{j_0}}\hookrightarrow
\ker'(i_{j_0})$.

Thus for each $j\in J$,
\begin{equation}\label{inject_eqn1}
\ker'(\alpha_j)(x)=\lambda(i_j)(x_j).
\end{equation}
Set $y\definicio \rho(i_0)(x)$ and take $j\in J$. Then
$$
F(\alpha_j)(y)=F(\alpha_j)(\rho(i_0)(x))=\rho(i_0)(\ker'(\alpha_j)(x))=\rho(i_0)(\lambda(i_j)(x_j))=\sigma(i_0)(x_j)=x_j
$$
and we are done.

This completes the proof of

\begin{Lem}\label{lem_inj_pseudo_graded}
Let $F:\P\rightarrow \Ab$ be an injective functor over a graded
poset $\P$. Then $F$ is pseudo-injective.
\end{Lem}

\begin{Ex}
For the category $\P$ with shape
$$\cdot \rightarrow \cdot $$
the functor $F:\P\rightarrow \Ab$ with values
$$\Z[p^\infty] \stackrel{\iota}\hookrightarrow \Q/\Z,$$ where $p$ is prime, is not
injective as the inclusion $\iota$ is not surjective, in spite of
the $\ker$'s are $\Q/\Z$ and $0$, which are injective objects of
$\Ab$.
\end{Ex}

Now we define \emph{pre-injective} objects
\begin{Defi}\label{Defi_preinjective}
Let $F:\P\rightarrow \Ab$ be a functor over a graded poset $\P$.
We call $F$ \emph{pre-injective} if
\begin{enumerate}
\item for any $i_0\in \Ob(\P)$ $\ker(i_0)$ is injective.
\label{Defi_preinjective_1}\item $F$ is pseudo-injective.
\label{Defi_preinjective_2}
\end{enumerate}
\end{Defi}

Till now we have obtained that injective functors $\P\rightarrow
\Ab$ over graded posets are pre-injective. In fact, as the next
proposition shows, the restriction we did to graded posets is
worthwhile:

\begin{Prop}\label{pro_injective_bounded_graded}
Let $\P$ be a bounded above graded poset and $F:\P\rightarrow \Ab$
be a functor. Then $F$ is injective if and only if it is
pre-injective.
\end{Prop}
\begin{proof}
If $F$ is injective then Lemmas
\ref{lem_proj_pseudo_graded_ker_injective} and
\ref{lem_inj_pseudo_graded} prove that $F$ is pre-injective. So
assume that $F$ is an pre-injective functor.

We can suppose that the degree function $deg$ on $\P$ is
decreasing and takes values $\{...,3,2,1,0\}$, and that
$\Ob_0(\P)\neq \emptyset$.

To see that $F$ is injective in $\Ab^\P$, given a diagram of
functors with exact row as shown, we must find a natural
transformation $\rho:A\Rightarrow F$ making the triangle
commutative
$$
\xymatrix{
             & F  &\\
A  \ar@{==>}[ur]^\rho & B \ar@{=>}[l]^\lambda \ar@{=>}[u]^\sigma &
\ar@{=>}[l] 0. }
$$
We define $\rho$ inductively, beginning on objects of degree $0$
and successively on objects of degrees $1,2,3,..$.

So take $i_0\in \Ob_0(\P)$ of degree $0$, and restrict to the
diagram in $\Ab$ over $i_0$. By Definition
\ref{Defi_preinjective}, as $\ker(i_0)=F(i_0)$, $F(i_0)$ is an
injective abelian group. So we can close the following triangle
with a homomorphism $\rho(i_0)$
$$
\xymatrix{
             & F(i_0)  &\\
A(i_0) \ar@{-->}[ru]^{\rho(i_0)} & B(i_0)\ar[l]^{\lambda(i_0)}
\ar[u]_{\sigma(i_0)} & \ar[l]0.}
$$
As there are no arrows between degree $0$ objects then $\rho$
restricted to the full subcategory with objects of degree $0$ is
trivially a natural transformation. Now suppose that we have
defined $\rho$ on all objects of $\P$ of degree less than $n$
($n\geq 1$), and that the restriction of $\rho$ to the full
subcategory generated by these objects is a natural transformation
and verifies $\rho \circ \lambda=\sigma$.

The next step is to define $\rho$ on degree $n$ objects. So take
$i_0\in \Ob_n(\P)$ and consider the splitting
$$F(i_0)=\ker(i_0)\oplus\coim(i_0)$$
where
$$\ker(i_0)=\bigcap_{i_0\stackrel{\alpha}\rightarrow i,\alpha\neq
1_{i_0}} \ker F(\alpha).$$ To define $\rho(i_0)$ such that it
makes commutative the diagram
$$
\xymatrix{
             & \ker(i_0)\oplus\coim(i_0)    &\\
A(i_0) \ar@{-->}[ru]^{\rho(i_0)} &
\ar[l]^{\lambda(i_0)}B(i_0)\ar[u]^{\sigma(i_0)}  & \ar[l]0,}
$$
we define by components $\rho(i_0)=\rho(i_0)_K\oplus\rho(i_0)_C$
with $\rho(i_0)_K:A(i_0)\rightarrow \ker(i_0)$ and
$\rho(i_0)_C:A(i_0)\rightarrow \coim(i_0)$. Define $\rho(i_0)_K$
as a homomorphism which closes the diagram
$$
\xymatrix{   & \ker(i_0)\\
             & F(i_0)\ar@{->>}[u]_{p_{i_0}}    \\
A(i_0) \ar@{-->}[ruu]^{\rho(i_0)_K} &
\ar[l]^{\lambda(i_0)}B(i_0)\ar[u]^{\sigma(i_0)}  & \ar[l]0 }
$$
where $p_{i_0}$ is the projection. This homomorphism $\rho(i_0)_C$
does exist because $\ker(i_0)$ is an injective abelian group.

Defining $\rho(i_0)_C$ needs a little bit more of work:

\begin{Claim}\label{claim_injective_defining_rhoK}
Fix $a\in A(i_0)$. Then for each $l=0,..,n-1$ there is $y_l\in
F(i_0)$ such that
\begin{equation}\label{injective_claim_eqn1}
F(\alpha)(y_l)=\rho(i)(A(\alpha)(a))
\end{equation}
for each $\alpha:i_0\rightarrow i$ with $deg(i)\leq l$.
\end{Claim}

Notice that $\rho(i)$ in Equation (\ref{injective_claim_eqn1}) is
defined by the induction hypothesis and because $deg(i)\leq l<n$.

We prove Claim \ref{claim_injective_defining_rhoK} by induction.
The base case reduces to find $y_0\in F(i_0)$ such that
$$
F(\alpha)(y_0)=\rho(i)(A(\alpha)(a))
$$
for each $\alpha:i_0\rightarrow i$ with $deg(i)=0$. Notice that,
as there are no objects of negative degree, then $\ker(i)=F(i)$
for each $i$ with $deg(i)=0$. Then $\rho(i)(A(\alpha)(a))\in
F(i)=\ker(i)$ for each $i$ with $deg(i)=0$ and thus $y_0$ exists
because $F$ is $n$-pseudo-injective (taking $J\definicio\{i|
\exists \alpha:i_0\rightarrow i \textrm{ and $deg(i)=0$}\}$ ).

Now we prove the inductive step of Claim
\ref{claim_injective_defining_rhoK}. Suppose we have $y_{l-1}$
verifying the hypothesis of the claim for $l<n$. Then we construct
$y_l$. Consider any object $i$ of degree $l$($<n$) such that there
is an arrow $\alpha:i_0\rightarrow i$. We have the elements
$\rho(i)(A(\alpha)(a))$ and $F(\alpha)(y_{l-1})$ in the abelian
group $F(i)$. Moreover, for any $\alpha':i\rightarrow i'$ with
$\alpha'\neq 1_i$ it holds that
\begin{equation}\label{injective_claim_eqn2}
F(\alpha')(F(\alpha)(y_{l-1}))=F(\alpha'\circ
\alpha)(y_{l-1})=\rho(i')(A(\alpha'\circ \alpha)(a))
\end{equation}
by the induction hypothesis on $y_{l-1}$ and because
$deg(i')<deg(i)=l$. Also we have that
\begin{equation}\label{injective_claim_eqn3}
F(\alpha')(\rho(i)(A(\alpha)(a)))=\rho(i')A(\alpha'\circ
\alpha)(a))
\end{equation}
because $\rho$ is a natural transformation on objects of degree
less than $n$ and because $deg(i)=l<n$.

Equations (\ref{injective_claim_eqn2}) and
(\ref{injective_claim_eqn3}) give that
$\rho(i)(A(\alpha)(a))-F(\alpha)(y_{l-1})\in \ker(i)$. Then
considering $J\definicio\{i| \exists \alpha:i_0\rightarrow i
\textrm{ and $deg(i)=l$}\}$ and applying that $F$ is
$(n-l)$-pseudo-injective we obtain $y'\in F(i_0)$ such that
$$
F(\alpha)(y')=\rho(i)(A(\alpha)(a))-F(\alpha)(y_{l-1})
$$
for each $\alpha:i_0\rightarrow i$ with $deg(i)=l$. Define
$$y_l\definicio y_{l-1}+y'.$$
To check that $y_l$ satisfies the statement of the claim take
$\alpha:i_0\rightarrow i$ with $deg(i)\leq l$. If $deg(i)=l$ then
\begin{align}
F(\alpha)(y_l)&=F(\alpha)(y_{l-1}+y')\notag\\
&=F(\alpha)(y_{l-1})+(\rho(i)(A(\alpha)(a))-F(\alpha)(y_{l-1}))\notag\\
&=\rho(i)(A(\alpha)(a)).\notag
\end{align}
If $deg(i)<l$ then
$$F(\alpha)(y_l)=F(\alpha)(y_{l-1}+y')=\rho(i)(A(\alpha)(a))+F(\alpha)(y').$$
In fact, $F(\alpha)(y')=0$: factor $\alpha$ as
$$
\xymatrix{ &i'\ar[rd]^{\beta}& \\
i_0\ar[ru]^{\alpha'} \ar[rr]^{\alpha}&&i }
$$
where $i'$ has degree $l$. Then
\begin{align}
&F(\alpha)(y')=F(\beta)(F(\alpha')(y'))=F(\beta)(\rho(i')(A(\alpha')(a))-F(\alpha')(y_{l-1}))=\notag
\\
=&\rho(i)(A(\alpha)(a))-F(\alpha)(y_{l-1})=0 \notag
\end{align}
by using the induction hypothesis on $y_{l-1}$ and on the
naturality of $\rho$. This finishes the proof of the claim.

Now we define $\rho(i_0)_C$. Take $a\in A(i_0)$ and define
$\rho(i_0)_C(a)=\overline{y}\in \coim(i_0)$, where $y\in F(i_0)$
is such that
$$
F(\alpha)(y)=\rho(i)(A(\alpha)(a))
$$
for each $\alpha:i_0\rightarrow i$ with $\alpha\neq 1_{i_0}$. The
element $y$ is obtained by taking $y=y_l$ with $l=n-1$ in Claim
\ref{claim_injective_defining_rhoK}.

It is straightforward that $\rho(i_0)_C$ is well defined because
if $y'\in F(i_0)$ also verifies that
$F(\alpha)(y)=\rho(i)(A(\alpha)(a))$  for each
$\alpha:i_0\rightarrow i$ with $\alpha\neq 1_{i_0}$ then $y-y'\in
\ker(i_0)$ and thus $\overline{y}=\overline{y'}$ in $\coim(i_0)$.
Also it is clear that $\rho(i_0)_C$ is a homomorphism of abelian
groups.

It remains to prove that $\rho(i_0)\circ\lambda(i_0)=\sigma(i_0)$.
Take $b\in B(i_0)$ and write
$\sigma(i_0)(b)=p_{i_0}(\sigma(i_0)(b))\oplus
\overline{\sigma(i_0)(b)}\in F(i_0)$. We want to see that
$\rho(i_0)_K(\lambda(i_0)(b))=p_{i_0}(\sigma(i_0)(b))$ and that
$\rho(i_0)_C(\lambda(i_0)(b))=\overline{\sigma(i_0)(b)}$. The
equation
$$
\rho(i_0)_K(\lambda(i_0)(b))=p_{i_0}(\sigma(i_0)(b))
$$
holds by definition of $\rho(i_0)_K$. Moreover, for each
$\alpha:i_0\rightarrow i$ with $\alpha\neq 1_{i_0}$, we have that
$$
F(\alpha)(\sigma(i_0)(b))=\sigma(i)(B(\alpha)(b))
$$
as $\sigma$ is a natural transformation, and that
$$
\rho(i)(A(\alpha)(\lambda(i_0)(b)))=\rho(i)(\lambda(i)(B(\alpha)(b)))=\sigma(i)(B(\alpha)(b))
$$
as $\lambda$ is a natural transformation and as $\rho\circ
\lambda=\sigma$ holds on objects of degree less than $n$ by the
induction hypothesis. Then, by the definition of $\rho(i_0)_C$,
$$\rho(i_0)_C(\lambda(i_0)(b))=\overline{\sigma(i_0)(b)}.$$

Defining $\rho(i_0)$ in this way for every $i_0\in \Ob_n(\P)$ we
have now $\rho$ defined on all objects of $\P$ of degree less or
equal than $n$. The last thing to do in order to complete the
inductive step is to prove that $\rho$ restricted to the full
subcategory over these objects is a natural transformation. Take
$\alpha:i_0\rightarrow i$ of degree different from zero in this
full subcategory. If the degree of $i_0$ is less than $n$ then the
commutativity of
$$
\xymatrix{
F(i_0)\ar[r]^{F(\alpha)}   &   F(i) \\
A(i_0)\ar[r]^{A(\alpha)}\ar[u]^{\rho(i_0)}                  &   A(i)\ar[u]^{\rho(i)}\\
}
$$
is granted by the inductive hypothesis. Suppose that the degree of
$i_0$ is $n$. Take $a\in A(i_0)$. Then
$\rho(i_0)(a)=\rho(i_0)_K(a)\oplus
\rho(i_0)_C(a)=\rho(i_0)_K(a)\oplus \overline{y}$ where $y\in
F(i_0)$ is such that
$$F(\alpha')(y)=\rho(i')(A(\alpha')(a))$$  for
each $\alpha':i_0\rightarrow i'$ with $\alpha'\neq 1_{i_0}$. Then
$$F(\alpha)(\rho(i_0)_K(a))=0$$
because $\rho(i_0)_K(a)\in \ker(i_0)$. Thus
$$F(\alpha)(\rho(i_0)(a))=F(\alpha)(y)$$
and $$F(\alpha)(y)=\rho(i)(A(\alpha)(a))$$ by the construction of
$y$.
\end{proof}

This proposition yields the following examples. The degree
functions $deg$ for the bounded above graded posets appearing in
the examples are indicated by subscripts $i_{deg(i)}$ on the
objects $i\in \Ob(\P)$ and take values $\{...,3,2,1,0\}$.

\begin{Ex}\label{examples_injective}

For the ``pullback category'' $\P$ with shape
$$\xymatrix { & a_0 \ar[d]^{f} \\  b_0\ar[r]^{g} &  c_1}$$
a functor $F:\P\rightarrow \Ab$ is injective if and only if
\begin{itemize}
\item $F(c)$, $\ker(F(f))$ and $\ker(F(g))$ are injective abelian
groups. \item $F(f)$ and $F(g)$ are epimorphisms.
\end{itemize}

For the inverse ``telescope category" $\P$ with shape
$$ \xymatrix{... \ar[r]& a_4 \ar[r]^{f_4} & a_3 \ar[r]^{f_3} & a_2 \ar[r]^{f_2} &  a_1 \ar[r]^{f_1} & a_0}$$
a functor $F:\P\rightarrow \Ab$ is injective if and only if
\begin{itemize}
\item $F(a_0)$ is an injective abelian group. \item $\ker(F(f_i))$
is an injective abelian group and $F(f_{i-d}\circ f_{i-d+1} \circ
..\circ f_{i})$ is an epimorphism for $d=0,1,..,i-1$ for each
$i=1,2,3,4,...$
\end{itemize}
\end{Ex}

\section{Pseudo-injectivity}
Consider a functor $F:\P\rightarrow \Ab$ over a small category
$\P$. In this section we look for, and find, conditions on $F$
such that $\liminv^i F=0$ for $i\geq 1$, i.e., we want conditions
such that the right derived functors of the left exact functor
$\liminv$ vanishes on $F$. We restrict along this section to
graded posets $\P$. Fix the following notation
\begin{Defi}\label{Defi_inverse_acyclic}
Let $\P$ be a graded poset and $F:\P\rightarrow \Ab$. We say $F$
is \emph{$\liminv$-acyclic} if $\liminv^i F=0$ for $i\geq 1$.
\end{Defi}

Recall that for injective objects it holds that any right derived
functor vanishes. So, from Proposition
\ref{pro_injective_bounded_graded}, we obtain firstly that

\begin{Prop}\label{pro_acyclic_injective}
Let $F:\P\rightarrow \Ab$ be a pre-injective functor over a
bounded above graded poset $\P$. Then $F$ is $\liminv$-acyclic.
\end{Prop}

Because being $\liminv$-acyclic is clearly weaker than being
injective we can wonder if is it possible to weaken the hypothesis
on Proposition \ref{pro_acyclic_injective} keeping the thesis of
$\liminv$-acyclicity. The answer is yes and the following theorem
states the appropriate conditions. Notice that we have removed the
condition (\ref{Defi_preinjective_1}) of Definition
\ref{Defi_preinjective}.

\begin{Thm}\label{pro_inv_acyclic_graded}
Let $F:\P\rightarrow \Ab$ be a pseudo-injective functor over a
bounded above graded poset $\P$. Then $F$ is $\liminv$-acyclic.
\end{Thm}
\begin{proof}
We can suppose that the degree function $deg$ on $\P$ is
decreasing and takes values $\{...,3,2,1,0\}$, and that
$\Ob_0(\P)\neq \emptyset$. To compute $\liminv^t F$ we use the
normalized (see \ref{normalized_ss}) spectral sequences
corresponding to the sixth row of Table \ref{tabla_ss} of Chapter
\ref{spectral}.

Fix $t\geq 1$. To prove $\liminv^t F=0$ it is enough to show that
$E_1^{p,t-p}$ is zero for every $p$. The contributions to
$E_1^{p,t-p}$  comes from $(E^p)^{\infty}_{p',p-p'-t}$ for $p'\leq
p-t$ (we are using normalized spectral sequences, Remark
\ref{normalized_ss}). We prove that
$$
(E^p)^{r}_{p',p-p'-t}=0
$$
if $r$ is big enough for each $p$ and $p'\leq p-t$. This implies
that $\liminv^t F=0$.

Consider the decreasing filtration $L^*$ of $C^*(\P,F)$ that gives
rise to the spectral sequence $E_*^{*,*}$. The $n$-simplices are
$$
L^p_n=L^pC^n(\P,F)=\prod_{\sigma \in {N\P}_n, deg(\sigma_0)\geq p}
F^{\sigma}.
$$
For each $p$ we have an increasing filtration $M_p^*$ of the
quotient $L^p/L^{p+1}$ that gives rise to the spectral sequence
$(E^p)^*_{*,*}$ and which $n$-simplices are
$$
(M_p)^{p'}_n=\prod_{\sigma \in {N\P}_n, deg(\sigma_0)=p,
deg(\sigma_n)\leq p'} F^{\sigma}.
$$
For $p'\leq p-t$ the abelian group $(E^p)^r_{p',q'}$ at the
$t=-(p'+q')+p$ simplices is given by
$$
(E^p)^r_{p',q'}= (M_p)^{p'}_t\cap
d^{-1}((M_p)^{p'-r}_{t+1})/(M_p)^{p'-1}_t\cap
d^{-1}((M_p)^{p'-r}_{t+1})+(M_p)^{p'}_t\cap
d((M_p)^{p'+r-1}_{t-1})
$$
where $d$ is the differential of the quotient $L^p/L^{p+1}$
restricted to the subgroups of the filtration $(M_p)^*$. For
$r>p'-p+t+1$ there are not $(t+1)$-simplices beginning in degree
$p$ and ending in degree less or equal to $p'-r<p-(t+1)$, i.e.,
$(M_p)^{p'-r}_{t+1}=0$. Also $r$ big enough implies
$(M_p)^{p'+r-1}_{t-1}=(M_p)^p_{t-1}=(L^p/L^{p+1})_{t+1}$, i.e.,
$(M_p)^{p'+r-1}_{t-1}$ equals all the $(t-1)$-simplices that begin
on degree $p$. Thus there exists $r$ such that
\begin{equation}
\label{E_p_r_big} (E^p)^r_{p',q'}= (M_p)^{p'}_t\cap d^{-1}(0)/(M_p)^{p'-1}_t\cap
d^{-1}(0)+(M_p)^{p'}_t\cap d((M_p)^p_{t-1})
\end{equation}
Fix such an $r$ and take $[x]\in(E^p)^r_{p',q'}$ where
\begin{equation}
\label{equ_exp_x} x\in (M_p)^{p'}_t = \prod_{\sigma \in {N\P}_t,
deg(\sigma_0)=p, deg(\sigma_t)\leq p'} F^{\sigma}
\end{equation}
and $d(x)=0$. We prove that $[x]=0$ in two steps:

\textbf{Step 1:} In this first step we find a representative
$x'\in (M_p)^{p'}_t\cap d^{-1}(0)$ for $[x]$ such that
$x'_\sigma=0$ for every $\sigma=\xymatrix{ \sigma_0
\ar[r]^{\alpha_1} & \sigma_1
\ar[r]^{\alpha_2}&...\ar[r]^{\alpha_{t-1}}&
\sigma_{t-1}\ar[r]^{\alpha_t} & \sigma_t}$ with
$deg(\sigma_{t-1})\leq p'$.

Notice that because $\P$ is bounded above then
$$(M_p)^{p'}_t = \prod_{\sigma \in {N\P}_t, deg(\sigma_0)=p, deg(\sigma_t)\leq p'}
F^{\sigma}=\prod_{\sigma \in {N\P}_t, deg(\sigma_0)=p, 0\leq
deg(\sigma_t)\leq p'} F^{\sigma}.$$

We need the following:

\begin{Claim}\label{claim_inj_step_1}
For each $l=0,...,p'$ there exists a representative $x'_l\in
(M_p)^{p'}_t\cap d^{-1}(0)$ for $[x]$ such that $(x'_l)_\sigma=0$
for every $$\sigma=\xymatrix{ \sigma_0 \ar[r]^{\alpha_1} &
\sigma_1 \ar[r]^{\alpha_2}&...\ar[r]&
\sigma_{t-1}\ar[r]^{\alpha_t} & \sigma_t}$$ with
$deg(\sigma_{t-1})\leq l$.
\end{Claim}

Notice that taking $l=p'$ in the claim the step $1$ is finished.
The case $l=0$ in the claim is fulfilled by taking $x'_0=x$
because there are no objects of negative degree.

Now we build $x'_l$ from $x'_{l-1}$ for $1\leq l\leq p'$. For that we need the

\begin{Claim}\label{claim_inj_step_2}
For each $m=0,...,l$ there exists a representative $x'_{l,m}\in
(M_p)^{p'}_t\cap d^{-1}(0)$ for $[x]$ such that
$(x'_{l,m})_\sigma=0$ for every $$\sigma=\xymatrix{ \sigma_0
\ar[r]^{\alpha_1} & \sigma_1 \ar[r]^{\alpha_2}&...\ar[r]&
\sigma_{t-1}\ar[r]^{\alpha_t} & \sigma_t}$$ with
$deg(\sigma_{t-1})<l$ or $deg(\sigma_{t-1})\leq l$ and
$deg(\sigma_t)<m$.
\end{Claim}

Notice that taking $m=l$ in Claim \ref{claim_inj_step_2} we obtain
$x'_l\definicio x'_{l,l}$ such that $(x'_l)_\sigma=0$ if
$deg(\sigma_{t-1})\leq l$, i.e., such that $x'_l$ satisfies Claim
\ref{claim_inj_step_1}. Define $x'_{l,0}\definicio x'_{l-1}$. This
satisfies the claim for $m=0$ as there is no object of negative
degree. Now we build $x'_{l,m}$ from $x'_{l,m-1}$ for $1\leq m\leq
l$.

By hypothesis $d(x'_{l,m-1})=0$. The differential $d$ is the restriction of the induced
differential on $L^p/L^{p+1}$ to
$$d:(M_p)^{p'}_t\rightarrow (M_p)^{p'}_{t+1}.$$
Recall that
$$d=\sum_{j=1,..,t+1} (-1)^j d^j$$
and that, for $\sigma\in N\P_{t+1}$, the final object of
$d_j(\sigma)$ is $\sigma_{t+1}$ for each $j=1,2,..,t$ and
$\sigma_{t}$ for $j=t+1$.

Thus for every $\epsilon\in N\P_{t+1}$ with $deg(\epsilon_0)=p$
and $deg(\epsilon_{t+1})\leq p'$, we can apply the projection
$$\pi^\epsilon:(M_p)^{p'}_{t+1}\twoheadrightarrow F^\epsilon$$
and obtain that $\pi^\epsilon(d(x'_{l-1}))=0$.

Now fix $\sigma\in N\P_t$ with $deg(\sigma_0)=p$,
$deg(\sigma_{t-1})=l$ and $deg(\sigma_t)=m-1$. Consider the object
$\sigma_t\in \Ob_{m-1}(\P)$ and any morphism
$\alpha:\sigma_t\rightarrow i$ with $\alpha\neq 1_{\sigma_t}$.
Write
$$\epsilon\definicio\xymatrix{
\sigma_0 \ar[r]^{\alpha_1} & \sigma_1 ...\ar[r]&
\sigma_{t-1}\ar[r]^{\alpha_t} & \sigma_t \ar[r]^\alpha & i.}$$

Notice that $deg(\epsilon)=p$ and
$deg(\epsilon_{t+1})=deg(i)<deg(\sigma_t)=m-1\leq l-1<p'$. By
hypothesis $\pi^\epsilon(d(x'_{l-1}))=0$. In fact, as $x'_{l,m-1}$
fulfills Claim \ref{claim_inj_step_2}, then
$\pi^\epsilon(d(x'_{l,m-1}))=F(\alpha)((x'_{l,m-1})_\sigma)$ and
so
$$F(\alpha)((x'_{l,m-1})_\sigma)=0.$$

Because $\alpha:\sigma_t\rightarrow i$ was arbitrary we have proven that
$$(x'_{l,m-1})_\sigma\in \ker(\sigma_t)$$
for each $\sigma\in N\P_t$ with $deg(\sigma_0)=p$,
$deg(\sigma_{t-1})=l$ and $deg(\sigma_t)=m-1$.

Consider now any $\tau\in N\P_{t-1}$ with $deg(\tau_0)=p$ and
$deg(\tau_{t-1})=l$. For the object $\tau_{t-1}$ of degree $l$
consider all the non-trivial morphisms from $\tau_{t-1}$ to
objects of degree $m-1$. Define $\Delta_\tau=\{\alpha\in
\Hom_\P(\tau_{t-1},i)|deg(i)=m-1\}$.

Notice that for each $\alpha\in \Delta_\tau$ and the $t$-simplex
$$\tau_\alpha\definicio\xymatrix{
\tau_0 \ar[r]^{\alpha_1} & \tau_1 ...\ar[r]&
\tau_{t-2}\ar[r]^{\alpha_{t-1}} & \tau_{t-1} \ar[r]^\alpha & i}
$$
we have proven before that $(x'_{l,m-1})_{\tau_\alpha}\in
\ker(i)$. Then, as $F$ is $(m-1-l)$-pseudo-injective, there is an
element $x_\tau\in F(\tau_{t-1})$ such that
$$F(\alpha)(x_\tau)=(x'_{l,m-1})_{\tau_\alpha}$$
for every $\alpha\in\Delta_\tau$.

Consider now the $(t-1)$-chain $y\in (M_p)^p_{t-1}$ given by
\begin{numcases}{y_\tau=}
x_\tau & if $deg(\tau_{t-1})=l$ \nonumber \\
0      & otherwise. \nonumber
\end{numcases}

If fact, as $y$ takes non zero values just on $\tau$ with
$deg(\tau_{t-1})=l\leq p'$, then $y\in (M_p)^{p'}_{t-1}$ and thus
$d(y)\in (M_p)^{p'}_t$ ($(M_p)^*$ is a filtration of a
differential complex). This means that $d(y)\in (M_p)^{p'}_t\cap
d((M_p)^p_{t-1})$.

Define $x'_{l,m}=x'_{l,m-1}-d(y)$. Then $x'_{l,m}\in (M_p)^{p'}$
and $d(x'_{l,m})=0$. Thus $[x'_{l,m}]=[x'_{l,m-1}]=[x]$. That
$x'_{l,m}$ fulfills Claim \ref{claim_inj_step_2} is clear by
construction. This finishes the proofs of Claims
\ref{claim_inj_step_2} and \ref{claim_inj_step_1}.

\textbf{Step 2:} By the step $1$ we can suppose that
$$
x_\sigma=0 $$ for $\sigma \in {N\P}_t$ with $deg(\sigma_0)=p$ and
$deg(\sigma_{t-1})\leq p$. Our objective now is to see that there
exists $y\in (M_p)^p_{t-1}$ with $d(y)=x$. This implies that
$[x]=0$ and finishes the proof of the theorem. We need the
\begin{Claim}\label{claim_inj_step_3} There are chains $x_i\in (M_p)^p_t$ for
$i=p',p'+1,..,p-t-1,p-t,p-t+1$ and $y_i \in (M_p)^p_{t-1}$ for
$i=p',p'+1,..,p-t-1,p-t$ such that
\begin{equation}
d(y_i)=x_i+x_{i+1}
\end{equation}
for $i=p',p'+1,..,p-t-1,p-t$ with $x_{p'}=x$ and $x_{p-t+1}=0$.
Moreover
\begin{enumerate}[(I)] \item \label{claim_inj_step_2_cond_1} $(x_i)\sigma=0$ if $deg(\sigma_{t-1})\leq i$
for $i=p',..,p-t+1$ \item \label{claim_inj_step_2_cond_2} $d(x_i)=0$ for $i=p',..,p-t+1$.
\end{enumerate}
\end{Claim}
Notice that the claim finishes the step $2$: as $x_{p-t+1}=0$ then
$x_{p-t}=d(y_{p-t})$,
$x_{p-t+1}=d(y_{p-t+1})-x_{p-t}=d(y_{p-t+1}-y_{p-t})$,
$x_{p-t+2}=d(y_{p-t+2})-x_{p-t+1}=d(y_{p-t+2}-y_{p-t+1}+y_{p-t})$,..,
$x=x_{p'}=d(y_{p'})-x_{p'+1}=d(y_{p'}-y_{p'+1}+...+(-1)^{p'-(p+t)}y_{p-t})$
where $y_{p'}-y_{p'+1}+...+(-1)^{p'-(p+t)}y_{p-t}\in
(M_p)^p_{t+1}$.

Define $x_{p'}\definicio x$. Then conditions
(\ref{claim_inj_step_2_cond_1}) and
(\ref{claim_inj_step_2_cond_2}) are satisfied for $i=p'$. We
construct $y_i$ and $x_{i+1}$ from $x_i$ recursively beginning
with $i=p'$. The element $y_i$ is built as in Claim
\ref{claim_inj_step_1}, the only difference being that now $y_i$
lies in $(M_p)^i_{t-1}\subseteq (M_p)^p_{t-1}$, but not in
$(M_p)^{p'}_{t-1}$. Then write $x_{i+1}\definicio d(y_i)-x_i$.
Notice that $d(x_{i+1})=d^2(y)-d(x_i)=0-0=0$ by the induction
hypothesis on $i$. Also notice that $x_{i+1}$ satisfies
(\ref{claim_inj_step_2_cond_1}) of the Claim
\ref{claim_inj_step_3} by the construction of $y_i$.

When constructing $y_{p-t}$ and $x_{p-t+1}$ from $x_{p-t}$ notice
that, because $\tau \in {N\P}_{t-1}$ with $deg(\tau_0)=p$ and
$deg(\tau_{t-1})=p-(t-1)$ imply that every morphism in $\tau$ is
of degree $1$, then the condition $d_j(\sigma)=\tau$ for
$\sigma\in{N\P}_{t}$ only holds for $j=t$. This implies that
$x_{p-t+1}=0$ by the construction of $y_{p-t}$.
\end{proof}

\begin{Ex}\label{examples_liminv acyclic}
For the '`pullback category'' $\P$ with shape
$$\xymatrix { & a_0 \ar[d]^{f} \\  b_0\ar[r]^{g} &  c_1}$$
a functor $F:\P\rightarrow \Ab$ is $\liminv$-acyclic if $F(f)$ and
$F(g)$ are epimorphisms.

For the inverse `telescope category' $\P$ with shape
$$ \xymatrix{... & a_4 \ar[r]^{f_4} & a_3 \ar[r]^{f_3} & a_2 \ar[r]^{f_2} &  a_1 \ar[r]^{f_1} & a_0}$$
a functor $F:\P\rightarrow \Ab$ is $\liminv$-acyclic if
$F(f_{i-d}\circ f_{i-d+1} \circ ..\circ f_{i})$ is an epimorphism
for $d=0,1,..,i-1$ for each $i=1,2,3,4,...$
\end{Ex}

\section{Computing higher limits}\label{section_liminvII}
Theorem \ref{pro_inv_acyclic_graded} shows that over a bounded
above graded poset pseudo-injectivity is enough for
$\liminv$-acyclicity. But it turns out that pseudo-injectivity is
not necessary for $\limdir$-acyclicity:

\begin{Ex}\label{liminvacyII_ex1}
For the `pushout category' $\P$ with shape
$$\xymatrix {a_0 \ar[r]^{f}\ar[d]^{g} & b_1 \\ c_1}$$
a functor $F:\P\rightarrow \Ab$ is pseudo-injective if and only if
$F(a)\stackrel{F(f)\oplus F(g)}\longrightarrow F(b)\oplus F(c)$ is
an epimorphism. But a straightforward calculus shows that
$\liminv^i F=0$ for $i\geq 1$ for any $F$.
\end{Ex}

Anyway, we shall see how pseudo-injectivity allows us to obtain a
better knowledge of the higher limits $\liminv^i F$. We begin with

\begin{Defi}
Let $F:\P\rightarrow \Ab$ be a functor over a graded poset $\P$.
Then $F':\P\rightarrow \Ab$ is the functor which takes values on
objects
$$F'(i_0)=\bigoplus_{\alpha:i_0\rightarrow i} F(i)$$
for $i_0\in \Ob(\P)$. For $\beta\in \Hom(\P)$,
$\beta:i_1\rightarrow i_0$, $F'(\beta)$ is the only homomorphism
which makes the diagram
$$
\xymatrix{ F'(i_1)\ar[r]^{F'(\beta)} & F'(i_0) \\
F(i)_{i_1\rightarrow i}\ar@{^{(}->}[u]\ar[r]^{1}&
F(i)_{i_0\rightarrow i}\ar@{^{(}->}[u]
 }
$$
to commute for each $\alpha:i_1\rightarrow i$ that factors through
$\beta$
$$
\xymatrix{ &i_0\ar[rd]& \\
i_1\ar[ru]^\beta \ar[rr]&&i, }
$$
and the diagram
$$
\xymatrix{ F'(i_1)\ar[r]^{F'(\beta)} & F'(i_0) \\
F(i)_{i_1\rightarrow i}\ar@{^{(}->}[u]\ar[ru]^{0} }
$$
for each $\alpha:i_1\rightarrow i$ that does not factor through
$\beta$.
\end{Defi}

Notice that $F'$ is built from $F$ as $\ker'$ was built from
$\ker$ in Section \ref{section_injective}, and that
$\ker_{F'}(i)=F(i)$ for each $i\in \Ob(\P)$. A nice property of
$F'$ is

\begin{Lem}\label{lema F'invadjunto}
Let $F:\P\rightarrow \Ab$ be a functor over a graded poset. Then
for each $G\in \Ab^\P$ there is a bijection

$$\xymatrix{\Hom_{\Ab^\P}(G,F')\ar[r]^<<<<{\varphi}_<<<<{\cong} & \prod_{i\in \Ob(\P)} \Hom_{\Ab}(G(i),F(i)).}$$
\end{Lem}
The proof is analogous to that of Lemma \ref{lema F'adjunto}.
Another interesting property of $F'$ is the following

\begin{Lem}\label{F'invlim}
Let $F:\P\rightarrow \Ab$ be a functor over a graded poset. Then
$$\liminv F'\cong\prod_{i\in \Ob(\P)} F(i).$$
\end{Lem}
\begin{proof}
It is straightforward using the previous lemma.
\end{proof}

The main feature of $F'$ we shall use is

\begin{Lem}\label{F'inv_pseudo}
Let $F:\P\rightarrow \Ab$ be a functor over a graded poset. Then
$F'$ is pseudo-injective.
\end{Lem}
\begin{proof}
It is straightforward.
\end{proof}

\begin{Rmk}\label{classical_enough injectives}
The monic natural transformation $F\Rightarrow G'$, where
$G=\Q/\Z^{\U\circ F}$ with $\U:\Ab\rightarrow \Set$ the forgetful
functor and $\Q/\Z^{-}:\Set\rightarrow \Ab$ the power on a set, is
a way to prove that $\Ab^\C$ has enough injectives for any small
category $\C$ (see \cite[243ff.]{cohn} on how to construct the
maps $F(i)\rightarrow G(i)$ for $i\in \Ob(\C)$).
\end{Rmk}
By Lemma \ref{lema F'invadjunto} there is a unique natural
transformation $\lambda:F\Rightarrow F'$ corresponding to the
family of homomorphisms $\{F(i)\stackrel{1_{F(i)}}\rightarrow
F(i)\}_{i\in \Ob(\P)}$. It is clear that $\lambda$ is a monic
natural transformation. Thus we can consider the object-wise
co-kernel of $\lambda$ to obtain a short exact sequence of
functors
$$
0\Rightarrow F\stackrel{\lambda}\Rightarrow F'\Rightarrow
C_F\Rightarrow 0.
$$

If $\P$ is bounded above then the long exact sequence (see Section
\ref{derivedfunctors}) associated to this short exact sequence
gives
\begin{numcases}{\liminv_j F=}
\liminv_{j-1} C_F & $j>1$ \nonumber \\
\coim\{\liminv F' \rightarrow \liminv C_F\}   & $j=1$ \nonumber
\end{numcases}
because $F'$ is $\liminv$-acyclic (it is pseudo-injective by Lemma
\ref{F'inv_pseudo} and apply Theorem
\ref{pro_inv_acyclic_graded}). Arguing as in the direct limit case
we have
\begin{Lem}
Let $\P$ be a bounded above graded poset and let $F:\P\rightarrow
\Ab$ be a functor. Then there are functors $C_j:\P\rightarrow \Ab$
for $j=0,1,2,...$ with $C_0=F$ and $C_1=\coim(F\Rightarrow F')$
such that
$$\liminv_j F=\liminv_{j-1} C_1=\liminv_{j-2} C_2=..=\liminv_1 C_{j-1}=\coim\{\liminv C'_{j-1}\rightarrow \liminv C_j\}$$
for each $j=0,1,2,...$.
\end{Lem}

The values $F'(i_0)$ and $C_F(i_0)$ can be very big. This can be
improved considering the functor $\ker:\P\rightarrow \Ab$ of
Section (\ref{section_injective}). Suppose that $\P$ is bounded
above and choose a family of homomorphisms $F(i)\rightarrow
\ker(i)$ such that the choice for the objects $i_0$ which do not
have any departing arrow is the identity $F(i_0)\rightarrow
F(i_0)=\ker(i_0)$. Then the natural transformation
$\lambda:F\Rightarrow \ker'$ given by Lemma \ref{lema
F'invadjunto} is monic in case $F$ is monic. We have

\begin{Lem}\label{compute_limiinv}
Let $\P$ be a bounded above graded poset and let $F:\P\rightarrow
\Ab$ be a monic functor. Then there is a short exact sequence of
functors
$$0\Rightarrow F\Rightarrow \ker'\Rightarrow C\Rightarrow 0$$
where $\ker'$ is $\liminv$-acyclic.
\end{Lem}

\chapter{Cohomology}\label{section_cohomology}

In this chapter we develop further the tools of Chapter
\ref{section_on inverse limit} in order to compute the integer
cohomology of graded posets $\P$ under some structural
assumptions, i.e., the existence of a ``covering family".
Throughout this chapter $\P$ shall be a bounded above graded poset
with decreasing degree function $deg$ which takes the values
$\{...,3,2,1,0\}$, and with $\Ob_0(\P)\neq \emptyset$. Conditions
\ref{graded_posets_finiteness assumption on starting arrows} are
assumed for $\P$ in the whole chapter.

The first two sections are devoted to the computation of higher
limits $\liminv^i F$, where $F:\P\rightarrow \Ab$ takes free
groups as values. After this, we specialize in compute
$H^i(\P;\Z)=\liminv^i c_\Z$, whereas $c_\Z:\P\rightarrow \Ab$ is
the constant functor of value $\Z$. Finally, we apply it to
simplex-like posets.

The key concepts defined in Sections \ref{section_liminvIII},
and \ref{section_liminvV}, i.e., $p$-condensed functors, free
functors, covering families and adequate covering families, are
all of local nature in the sense that they depend on
$(i_0\downarrow \P)$ and on $F|_{(i_0\downarrow \P)}$ for all
$i_0\in \Ob(\P)$.

\section{$p$-condensed
functors}\label{section_liminvIII} Recall that any family of
homomorphisms $\{F(i)\rightarrow \ker(i)\}_{\{i\in \Ob(\P)\}}$
gives a functor $\lambda:F\Rightarrow \ker'_F$ by Lemma \ref{lema
F'invadjunto}. If $\lambda$ is monic then we have a short exact
sequence of natural transformations
$$
\xymatrix{ 0\ar@{=>}[r]&F\ar@{=>}[r]^{\lambda}&\ker_F'\ar@{=>}[r]&
G\ar@{=>}[r]&0 }
$$
where $\ker_F'$ is $\liminv$-acyclic. This implies that $\liminv^1
F=\coim(\liminv \ker_F'\rightarrow \liminv G)$ and $\liminv^i
F=\liminv^{i-1} G$ for $i\geq 2$. Notice that $\liminv \ker_F'$ is
known by Lemma \ref{F'invlim}.

%
Next we find conditions on $F$ such that we can build a monic
natural transformation $\lambda$.

\begin{Defi}\label{Defi_p-condensed}
Let $F:\P\rightarrow \Ab$ be a functor. Suppose $\P$ has
decreasing degree function $deg:\Ob(\P)\rightarrow \{..,3,2,1,0\}$
and let $0\leq p\in \Z$. We say that $F$ is \emph{$p$-condensed}
if
\begin{enumerate}[(a)]
\item \label{condensed_b} $F(i)=0$ if $deg(i)<p$, and\item
\label{condensed_c} $\ker_F(i)=0$ if $deg(i)>p$.
\end{enumerate}
\end{Defi}

Notice that constant functors are $0$-condensed. If the functor
$F$ is $p$-condensed then we can consider the natural
transformation $\lambda:F\Rightarrow \ker'_F$ given by Lemma
\ref{lema F'invadjunto} for the maps $\tau_i:F(i)\rightarrow
\ker_F(i)$
\begin{numcases}{\tau_i=}
1_F(i)& if $deg(i)=p$ \nonumber \\
0     & otherwise. \nonumber
\end{numcases}
Notice that we have
\begin{numcases}{\ker_F(i)=}
F(i)& if $deg(i)=p$ \nonumber \\
0     & otherwise \nonumber
\end{numcases}
by hypothesis (\ref{condensed_b}) and (\ref{condensed_c}) in
Definition \ref{Defi_p-condensed}. In fact, the functor $\ker_F'$
takes values on objects
\begin{equation}\label{equ_expresionkerF'}
\ker_F'(i_0)=\prod_{i\in (i_0\downarrow \P)_{p}} F(i)
\end{equation}
 where
$(i_0\downarrow \P)_{p}=\{i\in \Ob(\P)|\exists i_0\rightarrow i,
deg(i)=p\}$ and $\ker_F(\beta)$ for $\beta:i_0\rightarrow i_1$ is
induced by the projections determined by $(i_1\downarrow
\P)_{p}\subset (i_0\downarrow \P)_{p}$. The homomorphism
$\lambda_i:F(i)\rightarrow \ker_F'(i)$ is given by
$$
\lambda_i=\prod_{i\in (i_0\downarrow \P)_{p}} F(\alpha_i)
$$
where $\alpha_i:i_0\rightarrow i$ is the unique arrow from $i_0$
to $i\in (i_0\downarrow \P)_{p}$. So $\lambda_i$ is a kind of
``diagonal". An easy induction argument on $deg(i)\in
\{p,p+1,..\}$ shows that $\lambda$ is a monic natural
transformation. We have obtained

\begin{Lem}\label{lem_pcondensed_shortexactsequence}
Let $F:\P\rightarrow \Ab$ be a $p$-condensed functor. Then there
is a short exact sequence
$$ \xymatrix{
0\ar@{=>}[r]&F\ar@{=>}[r]^{\lambda}&\ker_F'\ar@{=>}[r]&
G\ar@{=>}[r]&0. }
$$
\end{Lem}

$G$ is obtained by taking the object-wise co-image of $\lambda$.
On the object $i_0$, $G$ takes the value
\begin{equation}\label{equ_expresionG}
G(i_0)=\prod_{i\in (i_0\downarrow \P)_{p}} F(i)/
\lambda_{i_0}(F(i_0)).
\end{equation}

It is clear that $G$ verifies condition (\ref{condensed_b}) of
Definition \ref{Defi_p-condensed} for $p+1$, but in general
condition (\ref{condensed_c}) does not hold for $G$ and $p+1$.
More precisely, if $deg(i_0)>p+1$ then $\ker_G(i_0)=0$ is
equivalent to the natural map
$$
F(i_0)\rightarrow \liminv_{(i_0\downarrow \P)_*} F
$$
being an isomorphism, where $(i_0\downarrow \P)_*$ is the full
subcategory of $\P$ with objects $\{i|\exists i_0\rightarrow i,
i\neq i_0\}$. This natural map is a monomorphism by condition
(\ref{condensed_c}) of $F$ being $p$-condensed. So,
$\ker_G(i_0)=0$ if and only if $F(i_0)\rightarrow
\liminv_{(i_0\downarrow \P)_*} F $ is surjective. This is a local
property. We summarize these results in the following:

\begin{Lem}\label{lem_pcondensed_p+1condensed}
Let $F:\P\rightarrow \Ab$ be a $p$-condensed functor. Then there
is a short exact sequence
$$ \xymatrix{
0\ar@{=>}[r]&F\ar@{=>}[r]^{\lambda}&\ker_F'\ar@{=>}[r]&
G\ar@{=>}[r]&0. }
$$
Moreover, $G$ is $(p+1)$-condensed if and only if for each object
$i_0$ of degree greater than $p+1$, we have
$F(i_0)\stackrel{\cong}\rightarrow \liminv_{(i_0\downarrow \P)_*}
F$.
\end{Lem}


\section{Covering
families}\label{section_liminvV} In this section we study a bit
further the condition given in Lemma
\ref{lem_pcondensed_p+1condensed} which implies that the $G$ is
$(p+1)$-condensed, where the functor $G$ is defined by the short
exact sequence
$$ \xymatrix{
0\ar@{=>}[r]&F\ar@{=>}[r]^{\lambda}&\ker_F'\ar@{=>}[r]&
G\ar@{=>}[r]&0, }
$$
and $F$ is a $p$-condensed functor. This condition states that $G$
is $(p+1)$-condensed if and only if for each object $i_0$ of
degree greater than $p+1$ it holds that the map
\begin{equation}\label{equ_map_Fi0->liminvIi0*}
F(i_0)\rightarrow \liminv_{(i_0\downarrow \P)_*} F
\end{equation}
is an isomorphism. Recall that this map is a monomorphism if $F$
is $p$-condensed.

Fix $i_0$ of degree greater than $p+1$ and consider the map given
by restriction
$$
\liminv_{(i_0\downarrow \P)_*} F=\Hom_{(i_0\downarrow
\P)_*}(c_\Z,F)\rightarrow \prod_{i\in J} F(i)
$$
over the subset $J\subseteq(i_0\downarrow \P)_p$.  If this
restriction map turns out to be injective (notice that it is
injective for $J=(i_0\downarrow \P)_p$ because $F$ is
$p$-condensed) then the composition
$$F(i_0)\rightarrow \liminv_{(i_0\downarrow \P)_*} F\rightarrow
\prod_{i\in J} F(i)$$ is also injective. If $F$ is a free functor
(Definition \ref{defi_free_functor}) then both groups $F(i_0)$ and
$\prod_{i\in J} F(i)$ are free abelian groups (because we are
assuming Remark \ref{graded_posets_finiteness assumption on
starting arrows}). If the map
$$F(i_0)\rightarrow \prod_{i\in
J} F(i)$$ is pure then, by Lemma \ref{lema_fpuremono_imply_iso},
the condition $\dim F(i_0)=\sum_{i\in J} \dim F(i)$ implies that
this composition is an isomorphism and so
$F(i_0)\stackrel{\cong}\rightarrow \liminv_{(i_0\downarrow \P)_*}
F$. Thus we study the subsets $J\subseteq \Ob(\P)$ that make this
restriction map a pure monomorphism:

\begin{Defi}\label{defi_covering_family}
Let $\P$ be a bounded above graded poset with decreasing degree
function $deg$ which takes values $\{..,3,2,1,0\}$. A family of
subsets \linebreak $\J=\{J^{i_0}_p\}_{i_0\in \Ob(\P)\textit{,
$0\leq p\leq deg(i_0)$}}$ with $J^{i_0}_p\subseteq (i_0\downarrow
\P)_p$ is a \emph{covering family} if
\begin{enumerate}[a)]
\item \label{defi_covering_a}For each $i_0$ and $0\leq p<deg(i_0)$
it holds that $\bigcup_{i\in J^{i_0}_{p+1}} (i\downarrow \P)_p =
(i_0\downarrow \P)_p$ \item \label{defi_covering_b}For each $i_0$,
$0\leq p<deg(i_0)$ and $i\in J^{i_0}_{p+1}$ it holds that
$J^i_p\subseteq J^{i_0}_p$
\end{enumerate}
\end{Defi}

Notice that the definition above does not depend on a functor
defined over the category $\P$. Also, we have
$J^{i_0}_{deg(i_0)}=\{i_0\}$ by \ref{defi_covering_a}). The next
definition states the relation we expect between a covering family
and a $p$-condensed free functor
\begin{Defi}\label{defi_F_is_J_determined}
Let $\P$ be a bounded above graded poset, $\J$ be a covering
family and $F:\P\rightarrow \Ab$ be a $p$-condensed free functor.
We say that $F$ is \emph{$\J$-determined} if for any object $i_0$
of degree greater than $p+1$ the restriction map
$$
\liminv_{(i_0\downarrow \P)_*} F\rightarrow \prod_{i\in J^{i_0}_p}
F(i)
$$
is a monomorphism and the map
$$
F(i_0)\rightarrow \prod_{i\in J^{i_0}_p} F(i)
$$
is pure. If $deg(i_0)=p+1$ then we require that the last map above
is a pure monomorphism.
\end{Defi}

The main feature of covering families is that allow freeness plus
$\J$-determinacy to pass from $F$ to $G$. For an object $i_0$ with
$deg(i_0)\geq p+1$ notice that the map
$$
F(i_0)\rightarrow \prod_{i\in J^{i_0}_p} F(i)
$$
is a pure monomorphism as consequence of Definition
\ref{defi_F_is_J_determined}. The condition in Definition
\ref{defi_F_is_J_determined} for $deg(i_0)=p+1$ is added in order
to obtain that $G$ is a free functor. Notice that the following
proposition restricts to functors which take free abelian groups
as values.

\begin{Prop}\label{lem_covering_familiy_determination_hereditary}
Let $\P$ be a bounded above graded poset and $\J$ a covering
family. Assume that $F:\P\rightarrow \Ab$ is $p$-condensed, free
and $\J$-determined  and consider the functor $G$ defined by
$$ \xymatrix{
0\ar@{=>}[r]&F\ar@{=>}[r]^{\lambda}&\ker_F'\ar@{=>}[r]&
G\ar@{=>}[r]&0. }
$$
Then, if for each object $i_0$ with $deg(i_0)\geq p+1$ it holds
that $\dim F(i_0)=\sum_{i\in J^{i_0}_p} \dim F(i)$ then $G$ is
$(p+1)$-condensed, free and $\J$-determined.
\end{Prop}
\begin{proof}
Notice that the hypothesis implies that for any object $i_0$ of
degree $deg(i_0)>p+1$ the two maps
$$
F(i_0)\rightarrow \liminv_{(i_0\downarrow \P)_*} F\rightarrow
\prod_{i\in J^{i_0}_p} F(i)
$$
are isomorphisms. In particular,
$F(i_0)\stackrel{\cong}\rightarrow \liminv_{(i_0\downarrow \P)_*}
F$ and so $G$ is $(p+1)$-condensed.

If $deg(i_0)=p+1$ then the map
$$
F(i_0)\rightarrow \prod_{i\in J^{i_0}_p} F(i)
$$
is an isomorphism by hypothesis. Next we prove that $G$ is a free
functor. Consider any $i\in \Ob(\P)$ with $deg(i)\geq p+1$ (if
$deg(i)<p+1$ then $G(i)=0$) and the short exact sequence of
abelian groups
$$
0\rightarrow F(i)\stackrel{\lambda_i}\rightarrow
\ker_F'(i)\stackrel{\pi_i}\rightarrow G(i)\rightarrow 0.
$$
Then it is straightforward that the map
$$
s_i:\ker_F'(i)=\prod_{j\in (i\downarrow \P)_p}
F(j)\twoheadrightarrow \prod_{j\in J^i_p}
F(j)\stackrel{\cong}\rightarrow F(i)
$$
is a section of $\lambda_i$, i.e. $s_i\circ \lambda_i=1_{F(i)}$
(use that the restriction map $F(i)\rightarrow \prod_{j\in
J^{i}_p} F(j)$ is injective). This implies that the short exact
sequence above splits and so $G(i)$ is a subgroup of the free
abelian group $\ker_F'(i)$, and thus it is free as well. Next we
prove that $G$ is $\J$-determined.

Take $i_0$ of degree $n=deg(i_0)$ greater than $p+2$. We first
check that the restriction map
$$
\liminv_{(i_0\downarrow \P)_*} G\rightarrow \prod_{i\in
J^{i_0}_{p+1}} G(i)
$$
is injective. Consider any element $\psi\in
\liminv_{(i_0\downarrow \P)_*} G$ which is in the kernel of the
restriction map above. Notice that, as $deg(i_0)>p+2$, we can
consider the subset $J^{i_0}_{p+2}\subseteq (i_0\downarrow \P)_*$.
If for any $j\in J^{i_0}_{p+2}$ it holds that $\psi_j(1)=0$ then
$\psi=0$ by Definition
\ref{defi_covering_family}\ref{defi_covering_a}) and because $G$
is $(p+1)$-condensed.

Thus take $j\in J^{i_0}_{p+2}$. We want to see that $x\definicio
\psi_j(1)=0$. Recall the short exact sequence of abelian groups
$$
0\rightarrow F(j)\stackrel{\lambda_j}\rightarrow
\ker_F'(j)\stackrel{\pi_j}\rightarrow G(j)\rightarrow 0
$$
and take $y\in \ker_F'(j)$ such that $\pi_j(y)=x$. Recall that
$\ker_F'(j)=\prod_{i\in (j\downarrow \P)_{p}} F(i)$ and denote by
$\alpha_i:j\rightarrow i$ the unique arrow from $j$ to $i$ for
$i\in (j\downarrow \P)_p$.

Now consider the restriction $y|\in \prod_{i\in J^j_p} F(i)$.
Because $deg(j)=p+2>p+1$ the map $F(j)\rightarrow \prod_{i\in
J^j_p} F(i)$ is an isomorphism by hypothesis. Then there exists a
unique $z\in F(j)$ with $F(\alpha_i)(z)=y_i$ for each $i\in
J^j_p\subseteq (j\downarrow \P)_{p}$. If we prove that
$F(\alpha_i)(z)=y_i$ for each $i\in (j\downarrow \P)_{p}$ then
$\lambda_j(z)=y$. This implies that
$x=\pi_j(y)=\pi_j(\lambda_j(z))=0$ and finishes the proof.

Thus take $i\in (j\downarrow \P)_p$. By Definition
\ref{defi_covering_family}\ref{defi_covering_a}) there is
$\beta_i:i'\rightarrow i$ with $i'\in J^j_{p+1}$. Write
$\beta_{i'}:j\rightarrow i'$ for the unique arrow from $j$ to
$i'$. It holds that $\alpha_i=\beta_i\circ\beta_{i'}$. By
Definition \ref{defi_covering_family}\ref{defi_covering_b}) we
have that $J^j_{p+1}\subseteq J^{i_0}_{p+1}$. Thus
$G(\beta_{i'})(x)=G(\beta_{i'})(\psi_j(1))=\psi_{i'}(1)=0$ as
$\psi$ is in the kernel of the restriction map. The short exact
sequence
$$ 0\rightarrow F(i')\stackrel{\lambda_{i'}}\rightarrow
\ker_F'(i')\stackrel{\pi_{i'}}\rightarrow G(i')\rightarrow 0
$$
implies that there exists $t_{i'}\in F(i')$ such that
$\lambda_{i'}(t_{i'})=\ker_F'(\beta_{i'})(y)$. Consider
$z_{i'}=F(\beta_{i'})(z)$. We have that $z_{i'}$ and $t_{i'}$ have
the same image by the restriction map
$$
\liminv_{\P^{i'}_*} F\rightarrow \prod_{i\in J^{i'}_{p}} F(i)
$$
because $J^{i'}_p\subseteq J^j_p$. Because $F$ is $\J$-determined
then this restriction map is a monomorphism and so
$z_{i'}=t_{i'}$. This implies that
$$F(\alpha_i)(z)=F(\beta_i\circ\beta_{i'})(z)=F(\beta_i)(z_{i'})=F(\beta_i)(t_{i'})=y_i$$
and the proof of the restriction map being injective is finished.

Now we check that the map
$$
\omega:G(i_0)\rightarrow \prod_{i\in J^{i_0}_{p+1}} G(i)
$$
is pure. Take $z\in \prod_{i\in J^{i_0}_{p+1}} G(i)$ such that
there exists $x\in G(i_0)$ with $n\cdot z=\omega(x)$ for some
$n\neq 0$. We have to check that there exists $x'\in G(i_0)$ with
$z=\omega(x')$, or equivalently, that $x=n\cdot x'$ for some
$x'\in G(i_0)$. Recall once more the short exact sequence of
abelian groups
$$
0\rightarrow F(i_0)\stackrel{\lambda_{i_0}}\rightarrow
\ker_F'(i_0)\stackrel{\pi_{i_0}}\rightarrow G(i_0)\rightarrow 0
$$
and take $y\in \ker_F'(i_0)$ with $\pi_{i_0}(y)=x$. We are going
to build $h\in F(i_0)$ such that $y-\lambda_{i_0}(h)=n\cdot y'$,
i.e., such that for any $i\in (i_0\downarrow \P)_p$ the element
$(y-\lambda_{i_0}(h))_i=y_i-F(i_0\rightarrow i)(h)\in F(i)$ is
divisible by $n$. This implies that $x=n\cdot x'$ with
$x'=\pi_{i_0}(y')$.

Notice that by hypothesis for each $j\in J^{i_0}_{p+1}$,
$G(i_0\rightarrow j)(x)=n\cdot z_j \in G(j)$. This implies that
there exist $h_j\in F(j)$ and $y_j\in \ker_F'(j)$ such that
$\ker_F'(i_0\rightarrow j)(y)-\lambda_j(h_j)=n\cdot y_j$, i.e.,
such that for each $i\in (j\downarrow \P)_p\subseteq
(i_0\downarrow \P)_p$  we have that $y_i-F(j\rightarrow
i)(h_j)=n\cdot (y_j)_i\in F(i)$ (it is enough to take $y_j$ with
$\pi_j(y_j)=z_j$).

To build $h$ we use the map
$$
\tau:\prod_{i\in J^{i_0}_p} F(i)\stackrel{\cong}\rightarrow F(i_0)
$$
given by hypothesis, which is the inverse of the map
$$
F(i_0)\rightarrow \prod_{i\in J^{i_0}_p} F(i).
$$
For each $i\in J^{i_0}_p\subseteq (i_0\downarrow \P)_p$ choose, by
Definition \ref{defi_covering_family}\ref{defi_covering_a}),
$j(i)\in J^{i_0}_{p+1}$ such that there is an arrow
$j(i)\rightarrow i$. Then set $\eta_i=F(j(i)\rightarrow
i)(h_{j(i)})\in F(i)$, where $h_{j(i)}$ is built as before. Define
$h\definicio \tau(\eta)$. By construction $F(i_0\rightarrow
i)(h)=F(j(i)\rightarrow i)(h_{j(i)})$ for each $i\in J^{i_0}_p$
(but not for an arbitrary $i\in (i_0\downarrow \P)_p$).

With this definition for $h$ we check now that
$y_i-F(i_0\rightarrow i)(h)$ is divisible by $n$ for each $i\in
(i_0\downarrow \P)_p$. This finishes the proof. Fix $i\in
(i_0\downarrow \P)_p$ and $j_i\in J^{i_0}_{p+1}$ such that there
is an arrow $j_i\rightarrow i$ (we are not assuming that
$j_i=j(i)$ if $i\in J^{i_0}_p$). On the one hand we have by
hypothesis that
$$
y_k-F(j_i\rightarrow k)(h_{j_i})=n\cdot (y_{j_i})_k
$$
for each $k\in \P^{j_i}_p$. In particular,
\begin{equation}\label{equ_covefami_1}
y_k-F(j_i\rightarrow k)(h_{j_i})=n\cdot (y_{j_i})_k
\end{equation}
for each $k\in J^{j_i}_p$. Set $h'\definicio F(i_0\rightarrow
j_i)(h)$. Because $j_i\in J^{i_0}_{p+1}$ then, by Definition
\ref{defi_covering_family}\ref{defi_covering_b}),
$J^{j_i}_p\subseteq J^{i_0}_p$ and thus by construction for any
$k\in J^{j_i}_p$
$$
y_k-F(i_0\rightarrow k)(h)=y_k-F(j(k)\rightarrow
k)(h_{j(k)})=n\cdot (y_{j(k)})_k.
$$
Notice that $F(i_0\rightarrow k)(h)=F(j_i\rightarrow
k)(F(i_0\rightarrow j_i)(h))=F(j_i\rightarrow k)(h')$. On the
other hand, we have obtained
\begin{equation}\label{equ_covefami_2}
y_k-F(j_i\rightarrow k)(h')=n\cdot (y_{j(k)})_k
\end{equation}
for each $k\in \J^{j_i}_p$.

Now write $\eta_k=(y_{j(k)})_k-(y_{j_i})_k$ for each $k\in
J^{j_i}_p$ and write $h''=\tau(\eta)\in F(j_i)$ where $\tau$ is
the inverse of the map
$$
F(j_i)\rightarrow \prod_{i\in J^{j_i}_p} F(i).
$$
By Equations (\ref{equ_covefami_1}) and (\ref{equ_covefami_2}) it
is straightforward that the elements $h_{j_i}-n\cdot h''$ and $h'$
have the same image by this map. Then, as this map is injective by
hypothesis, $h'=h_{j_i}-n\cdot h''$. Then, recall that $i\in
J^{j_i}_p$,
$$
y_i-F(i_0\rightarrow i)(h)=y_i-F(j_i\rightarrow
i)(h')=y_i-F(j_i\rightarrow i)(h_{j_i}-n\cdot h''),
$$
and this equals
$$
n\cdot (y_{j_i})_i+n\cdot F(j_i\rightarrow i)(h'').
$$
Thus $y_i-F(i_0\rightarrow i)(h)$ is divisible by $n$.

If $deg(i_0)=p+2$ we have to see that the map
$$
\omega:G(i_0)\rightarrow \prod_{i\in J^{i_0}_{p+1}} G(i)
$$
is a pure monomorphism. To prove that $\omega$ is a monomorphism
use the proof above starting where $\psi_j$ is considered for an
arbitrary object $j$ of degree $p+2$. The proof of $\omega$ being
pure is exactly the same as above.
\end{proof}

\begin{Rmk}\label{rmk_dimensions}
Notice that, in the conditions of the proposition, and assuming
that $F$ takes finitely generated free abelian groups as values,
we have the following formula for the rank of the free abelian
group $G(i_0)$ for $deg(i_0)\geq p+1$
$$
\dim(G(i_0))=\sum_{i\in (i_0\downarrow \P)_{p}} \dim F(i) - \dim
F(i_0)
$$
(recall that we are assuming \ref{graded_posets_finiteness
assumption on starting arrows}). This is so because of the short
exact sequence of free abelian groups
$$
0\rightarrow F(i_0)\stackrel{\lambda_{i_0}}\rightarrow
\ker_F'(i_0)\stackrel{\pi_{i_0}}\rightarrow G(i_0)\rightarrow 0.
$$
\end{Rmk}
\begin{Rmk}\label{rmk_section}
Consider again the map $s_{i_0}:\ker'_F(i_0)\rightarrow F(i_0)$
with $s_{i_0}\circ \lambda_{i_0}=1_{F(i_0)}$ built in the proof of
the previous proposition for $deg(i_0)\geq p+1$. To $s_{i_0}$
corresponds the monomorphism
$$
\xymatrix{ \ar[r]^{\delta_{i_0}}G(i_0) & \ker'_F(i_0)}
$$
given by
$$
\xymatrix{ \pi_{i_0}(x)\ar@{|->}[r] & x-(\lambda_{i_0}\circ
s_{i_0})(x)},
$$
which satisfies $\pi_{i_0}\circ\delta_{i_0}=1_{G(i_0)}$. It is
straightforward that, by construction, $\im
\delta_{i_0}=\prod_{i\in (i_0\downarrow \P)_p\setminus J^{i_0}_p}
F(i)$, and thus
$$
G(i_0)\stackrel{\delta_{i_0}}\cong \prod_{i\in (i_0\downarrow
\P)_p\setminus J^{i_0}_p} F(i).
$$
Moreover, $x=\delta_{i_0}(y)$ is the only preimage of $y$ by
$\pi_{i_0}$ which verifies $x_i=0$ for $i\in J^{i_0}_p$.
\end{Rmk}

The main consequence of the previous proposition is that it
reduces the problem of whether $G$ is $(p+1)$-condensed to some
integral equations. Moreover, this procedure can be applied
recursively because the resulting functor $G$ turns out to be
$(p+1)$-condensed, free and $\J$-determined, and so the
proposition applies to $G$ too. Notice that the ranks of $G$ are
given by Remark \ref{rmk_dimensions}.

\subsection{Adequate covering families}\label{section_coho_local1}
In this section we apply the work developed in Sections
\ref{section_liminvIII} and \ref{section_liminvV} to compute the
cohomology with integer coefficients of the realization of a
bounded above graded poset $\P$ equipped with a covering family
$\J$.

To compute the abelian group $H^p(\P;\Z)$ for $p\geq 1$ we compute
the higher limit ${\liminv}^p c_\Z$ where $c_\Z:\P\rightarrow \Ab$
is the functor of constant value $\Z$ which sends every morphism
to the identity $1_\Z$. We begin studying the extra conditions
that the covering family $\J$ must satisfy to apply iteratively
the Proposition
\ref{lem_covering_familiy_determination_hereditary} beginning on
$c_\Z$.

First, notice that $c_\Z$ is $0$-condensed (we are assuming
\ref{graded_posets_finiteness assumption on starting arrows}) and
free (Definition \ref{defi_free_functor}). By Definition
\ref{defi_F_is_J_determined}, $c_\Z$ is $\J$-determined as
$0$-condensed functor if and only if for each $i_0\in \Ob(\P)$
with $deg(i_0)\geq 2$ the set $J^{i_0}_0$ intersects each
connected component of $(i_0\downarrow \P)_*$. The dimensional
equation in Proposition
\ref{lem_covering_familiy_determination_hereditary} for $i_0\in
\Ob(\P)$ with $deg(i_0)\geq 1$ becomes $\dim
c_\Z(i_0)=1=|J^{i_0}_0|=\sum_{i\in J^{i_0}_0} \dim c_\Z(i)$. Thus,
$c_\Z$ is $\J$-determined as $0$-condensed functor if and only if
$(i_0\downarrow \P)_*$ is connected for $deg(i_0)\geq 2$ and
$|J^{i_0}_0|=1$ for $deg(i_0)\geq 1$.

The successive applications of Proposition
\ref{lem_covering_familiy_determination_hereditary} give, by the
dimensional equation in the statement of the Proposition
\ref{lem_covering_familiy_determination_hereditary}, the following
structural conditions linking $\P$ and $\J$:

\begin{Defi}\label{defi_structural numbers of graded poset}
Let $\P$ be a bounded above graded poset. Define, inductively on
$p$, the number $R^{i_0}_p$ for each object $i_0$ with
$deg(i_0)\geq p$ by $R^{i_0}_0=1$ and by
$$R^{i_0}_p=\sum_{i\in (i_0\downarrow \P)_{p-1}} R^i_{p-1}-R^{i_0}_{p-1}$$
for $p\geq 1$.
\end{Defi}
\begin{Defi}\label{defi_adequate_covering_family}
Let $\P$ be a bounded above graded poset and $\J$ be a covering
family for $\P$. We say that $\J$ is \emph{adequate} if
$(i_0\downarrow \P)_*$ is connected for $deg(i_0)\geq 2$, and if
we have the equality
$$
R^{i_0}_p=\sum_{i\in J^{i_0}_p} R^i_p
$$
for $p\geq 0$ and $deg(i_0)\geq p+1$.
\end{Defi}

\begin{Prop}\label{cohomology with adequate covering family}
Let $\P$ be a bounded above graded poset and let $\J$ be an
adequate covering family. Then there is a sequence of functors
$F_0,F_1,F_2,..$ defined by $F_0\definicio c_\Z:\P\rightarrow \Ab$
and by the short exact sequence
$$ \xymatrix{
0\ar@{=>}[r]&F_{p-1}\ar@{=>}[r]^{\lambda_{p-1}}&\ker_{F_{p-1}}'\ar@{=>}[r]^{\pi_p}&
F_p\ar@{=>}[r]&0 }
$$
for $p=1,2,3,..$. Moreover, $F_p$ is $p$-condensed, free and
$\J$-determined for $p\geq 0$. For $deg(i_0)\geq p$ we have $\dim
F_p(i_0)=R^{i_0}_p$.
\end{Prop}
\begin{proof}
We prove by induction the following
\begin{Claim}
For $p=0,1,2,..,N$ there exist functors $F_0,F_1,F_2,..,F_N$ given
by $F_0\definicio c_\Z:\P\rightarrow \Ab$ and by a short exact
sequence
$$
\xymatrix{
0\ar@{=>}[r]&F_{p-1}\ar@{=>}[r]^{\lambda_{p-1}}&\ker_{F_{p-1}}'\ar@{=>}[r]^{\pi_p}&
F_p\ar@{=>}[r]&0 }
$$
for $p=1,2,3,..,N$. Moreover, for any $p=0,1,2,..,N$, $F_p$ is
$p$-condensed, free and $\J$-determined and for $deg(i_0)\geq p$
we have
$$\dim F_p(i_0)=R^{i_0}_p.$$
\end{Claim}

The base case of the claim, $N=0$, holds by the arguments before
the proposition and because $\dim F_0(i_0)=\dim
c_\Z(i_0)=1=R^{i_0}_0$ for any $i_0\in\Ob(\P)$.

Suppose that the claim holds for $N\geq 0$. Then we prove it for
$N+1$. We apply Proposition
\ref{lem_covering_familiy_determination_hereditary} to the
$N$-condensed, free and $\J$-determined functor $F_N$ given by the
induction hypothesis. Thus we have to check that for every object
$i_0$ of degree greater or equal to $N+1$ the following equality
holds
$$
\dim F_N(i_0)=\sum_{i\in {J}^{i_0}_N} \dim F_N(i).
$$

By the induction hypothesis this equations is exactly
$$R^{i_0}_N=\sum_{j\in J^{i_0}_N} R^i_N,$$
and this holds by the definition of adequate covering family. By
Remark \ref{rmk_dimensions} the rank of $F_{N+1}(i_0)$ for
$deg(i_0)\geq N+1$ is given by
$$
\dim F_{N+1}(i_0)=\sum_{i\in (i_0\downarrow \P)_N} \dim F_N(i) -
\dim F_N(i_0)
$$
which, by the induction hypothesis again, equals
$$
\dim F_{N+1}(i_0)=\sum_{i\in (i_0\downarrow \P)_N}
R^i_N-R^{i_0}_N.
$$
Thus, $\dim F_{N+1}(i_0)=R^{i_0}_{N+1}$ by Definition
\ref{defi_structural numbers of graded poset}.
\end{proof}

\subsection{Local basis and morphisms.}\label{section_coho_local2} In the previous
section we found a sequence of functors $F_0=c_\Z, F_1, F_2,...$
to compute the cohomology of a graded poset $\P$ equipped with an
adequate covering family $\J$. Moreover, we found the dimension of
$F_p(i)$ for any $i\in \P$. In this section we shall build
inductively an explicit basis for $F_p(i)$.

For $p=0$ and $i_0\in \Ob(\P)$ we choose the basis
$B^{i_0}_0=\{1\}$ of $F_0(i_0)=\Z$. For $p+1\geq 1$ and
$deg(i_0)\geq p+1$ notice that we have an isomorphism (cf.
{\ref{rmk_section})
$$
F_{p+1}(i_0)\stackrel{\delta_{i_0}}\cong \prod_{i\in
(i_0\downarrow \P)_p\setminus J^{i_0}_p} F_p(i)
$$
which sends $y\in F_{p+1}(i_0)$ to the only element $x\in
\ker'_{F_p}(i_0)$ which projects on $y$ by $\pi_{i_0}$ and with
$x_i=0$ for each $i\in J^{i_0}_p$. This proves the following

\begin{Lem}\label{lemma_coho_local_explicit basis}
Let $\P$ be a bounded above graded poset and let $\J$ be an
adequate covering family. Let $c_\Z,F_1,F_2,..$ be the sequence of
functors given by Proposition $\ref{cohomology with adequate
covering family}$. Then there are basis $B^i_p$ of $F_p(i)$ for
$p\geq 0$ and $deg(i)\geq p$ such that, for $p\geq0$ and
$deg(i_0)\geq p+1$, the map $\pi_{i_0}$ of the short exact
sequence
$$
0\rightarrow F_p(i_0)\stackrel{\lambda_{i_0}}\rightarrow \prod_{i\in
(i_0\downarrow \P)_p} F_p(i) \stackrel{\pi_{i_0}}\rightarrow
F_{p+1}(i_0)\rightarrow 0
$$
maps $\bigcup_{i\in (i_0\downarrow \P)_p\setminus J^{i_0}_p} B^i_p
$ bijectively onto $B^{i_0}_{p+1}$.
\end{Lem}

\begin{Rmk}\label{basis_coeff_equal_delta_i}
Notice that for $y=\pi_{i_0}(x)\in F_{p+1}(i_0)$, the expression
of $y$ in terms of the basis of the lemma above,
$$
y=\sum_{i\in (i_0\downarrow \P)_p\setminus J^{i_0}_p}
\sum_{l=1,..,R^i_p} y^l_i \cdot \pi_{i_0}(e^i_{p,l}),
$$
where $B^i_p=\{e^i_{p,1},..,e^i_{p,R^i_p}\}$, corresponds to the
expression of $\delta_{i_0}(y)=x-(\lambda_{i_0}\circ
s_{i_0})(x)\in \prod_{i\in (i_0\downarrow \P)_p\setminus
J^{i_0}_p} F_p(i)$,
$$
\delta_{i_0}(y)=\sum_{i\in (i_0\downarrow \P)_p\setminus
J^{i_0}_p} \sum_{l=1,..,R^i_p} y^l_i \cdot e^i_{p,l}
$$
in terms of the basis $\bigcup_{i\in (i_0\downarrow \P)_p\setminus
J^{i_0}_p} B^i_p$.
\end{Rmk}

Next we are interested in the expression of $F_p(\alpha)$ with
respect to the basis of Lemma \ref{lemma_coho_local_explicit
basis}. Fix the map $\alpha:i_1\rightarrow i_2$ with $deg(i_1)\geq
deg(i_2)\geq p\geq 1$. The basis $B^{i_1}_p$ is in 1-1
correspondence with the set
$$\bigcup_{i\in (i_1\downarrow \P)_p\setminus J^{i_1}_p}
B^i_{p-1},$$ and the basis $B^{i_2}_p$ is in 1-1 correspondence
with the set
$$\bigcup_{i\in (i_2\downarrow \P)_p\setminus J^{i_2}_p}
B^i_{p-1}.$$

Take an element $\pi_{i_1}(e^i_{{p-1},l})$ of $B^{i_1}_p$, where
$i\in (i_1\downarrow \P)_p\setminus J^{i_1}_p$ and $l\in
{1,..,R^i_{p-1}}$. Then
$$F_p(\alpha)(\pi_{i_1}(e^i_{{p-1},l}))=\pi_{i_2}(\ker'_{F_{p-1}}(\alpha)(e^i_{{p-1},l})).$$

We have
\begin{numcases}{\ker'_{F_{p-1}}(\alpha)(e^i_{{p-1},l})=}
e^i_{{p-1},l} & if $i\in (i_2\downarrow \P)_{p-1}$ \nonumber \\
0 & otherwise, \nonumber
\end{numcases}
and thus
\begin{numcases}{F_p(\alpha)(e^i_{{p-1},l})=}
\pi_{i_2}(e^i_{{p-1},l}) & if $i\in (i_2\downarrow \P)_{p-1}$ \nonumber \\
0 & otherwise. \nonumber
\end{numcases}

Suppose that $i\in (i_2\downarrow \P)_{p-1}$. If $i\notin
J^{i_2}_{p-1}$ then $\pi_{i_2}(e^i_{{p-1},l})$ is an element of
$B^{i_2}_p$ and so $F_p(\alpha)$ sends the element
$\pi_{i_1}(e^i_{{p-1},l})$ to the element
$\pi_{i_2}(e^i_{{p-1},l})$. If $i\in J^{i_2}_{p-1}$ then we have
to find the coefficients of $\pi_{i_2}(e^i_{{p-1},l})$ with
respect to the basis $B^{i_2}_p$. By the Remark
\ref{basis_coeff_equal_delta_i} these coefficients equal the
coefficients of $\delta_{i_2}(\pi_{i_2}(e^i_{{p-1},l}))$ with
respect to the basis $\bigcup_{i\in (i_2\downarrow
\P)_{p-1}\setminus J^{i_2}_{p-1}} B^i_{p-1}$ of $\prod_{i\in
(i_2\downarrow \P)_{p-1}\setminus J^{i_2}_{p-1}} F_{p-1}(i)$. This
proves the following

\begin{Lem}\label{lemma_coho_local_explicit morphism}
Let $\P$ be a bounded above graded poset and $\J$ an adequate
covering family. Let $\alpha:i_1\rightarrow i_2$ with
$deg(i_1)\geq deg(i_2)\geq p\geq 1$ be a map in $\P$. Then
$F_p(\alpha)$ is the identity when restricted to
$$
\langle \bigcup_{i\in (i_1\downarrow \P)_{p-1}\setminus
J^{i_1}_{p-1} \cap (i_2\downarrow \P)_{p-1}\setminus
J^{i_2}_{p-1}} B^i_{p-1} \rangle
$$
\end{Lem}

\begin{Rmk}\label{nice_morphism on J^i_p}
Notice that if $i_2\in J^{i_1}_p$ then Definition
\ref{defi_covering_family} implies that $J^{i_2}_{p-1}\subseteq
J^{i_1}_{p-1}$. Then, for $i\in (i_1\downarrow \P)_{p-1}\setminus
J^{i_1}_{p-1}$, if $i\in(i_2\downarrow \P)_{p-1}$ then $i\notin
J^{i_2}_{p-1}$. Thus, in this case, $(i_1\downarrow
\P)_{p-1}\setminus J^{i_1}_{p-1} \cap (i_2\downarrow
\P)_{p-1}\setminus J^{i_2}_{p-1}=(i_2\downarrow \P)_{p-1}\setminus
J^{i_2}_{p-1}$ and the lemma applies to the whole
$F_p(i_2)\subseteq F_p(i_1)$.
\end{Rmk}


\section{Global behaviour.}\label{section_coho_global}
In the previous section we saw (Proposition \ref{cohomology with
adequate covering family}) how the local properties of a graded
poset $\P$ equipped with an adequate covering family $\J$ give
rise to a sequence $F_0=c_\Z$, $F_1$, $F_2$,... of functors. In
this section we study some global properties of $\J$, as well as
we define global covering families.

The first global point to notice is that the sequence of functors
from Proposition \ref{cohomology with adequate covering family}
does not depend on the adequate covering family $\J$. Thus, two o
more adequate covering families can be considered for the the same
graded poset and they still give rise to the same sequence of
functors.

Next we focus on the short exact sequence
$$ \xymatrix{
0\ar@{=>}[r]&F_{p-1}\ar@{=>}[r]^{\lambda_{p-1}}&\ker_{F_{p-1}}'\ar@{=>}[r]^{\pi_p}&
F_p\ar@{=>}[r]&0 }
$$
for $p\geq 1$ of Proposition \ref{cohomology with adequate
covering family}.  The beginning of the long exact sequence of
this short exact sequence is
\begin{equation}\label{eqn_long_exact_sequence_H^p}
\xymatrix{0\ar[r]&\liminv F_{p-1}\ar[r]^{\iota_{p-1}}&\liminv
\ker_{F_{p-1}}'\ar[r]^{\omega_p}& \liminv F_p\ar[r]&H^p(\P;\Z)
\ar[r]&0 }
\end{equation}
by Lemma \ref{F'inv_pseudo}, where
$\iota_{p-1}=\widetilde{\lambda_{p-1}}$ and
$\omega_p=\widetilde{\pi_p}$.

Notice that the three inverse limits appearing above are free
abelian groups as the corresponding functors take free abelian
groups as values. In fact, for the middle term we have the exact
description
\begin{equation}\label{eqn_kerF' directproduct}
\liminv \ker_{F_{p-1}}' \cong\prod_{i\in \Ob_{p-1}(\P)} F_{p-1}(i)
\end{equation}
given by Lemma \ref{F'invlim}. It turns out that also there is a
simpler description for $\liminv F_p$, which can be interpreted as
the analogue to the fact that the cohomology on degree $n$ depends
upon the $n+1$ skeleton (recall that $F_p(i)=0$ if $deg(i)<p$):

\begin{Lem}\label{lemma_liminvFp p+2_esqueleto}
Let $\P$ be a bounded above graded poset and let $\J$ be an
adequate covering family. Let $c_\Z,F_1,F_2,..$ be the sequence of
functors given by Proposition $\ref{cohomology with adequate
covering family}$. Then
$$
\liminv F_p\cong\liminv F_p|_{\P_{\{p+1,p\}}}
$$
for each $p\geq 0$.
\end{Lem}
\begin{proof}
Consider the restriction map
$$
\liminv F_p\rightarrow \liminv F_p|_{\P_{\{p+1,p\}}}.
$$
This map is clearly a monomorphism because $F_p$ is a
$p$-condensed functor. To see that it is surjective take $\psi\in
\liminv F_p|_{\P_{\{p+1,p\}}}$. We want to extend $\psi$ to each
$i\in \Ob(\P)$ with $deg(i)>p+1$. We do it inductively on
$deg(i)$.

Notice that (see the beginning of the proof of Proposition
\ref{lem_covering_familiy_determination_hereditary})
$$
F_p(i)\rightarrow \liminv_{(i\downarrow \P)_*} F_p
$$
is an isomorphism for $deg(i)>p+1$. For $deg(i)=p+2$ we have that
$j\in {(i\downarrow \P)_*}$ implies that $deg(j)\leq p+1$. Then
there is a unique way of extending $\psi$ to $\psi(i)$. Once we
have extended $\psi$ to $\P_{\{p+2,p+1,p\}}$ we proceed with an
induction argument. That the extension that we are building is
actually a functor is again due to that $F_p$ is $p$-condensed.
\end{proof}

Also, from Equations (\ref{eqn_long_exact_sequence_H^p}) and
(\ref{eqn_kerF' directproduct}), we have the following formula,
analogue to that of the euler characteristic for $CW$-complexes:
\begin{Lem}\label{lemma_analogue eulerchar}
Let $\P$ be a bounded graded poset for which exists an adequate
covering family. Then
$$
\sum_i (-1)^i\dim H^i(\P;\Z)=\sum_i (-1)^i\sum_{p\in
\Ob_i(\P)}R^i_p.
$$
\end{Lem}

\subsection{Global covering families.}\label{section_global_coverings}
Recall that covering families were defined as subsets of the local
categories $(i\downarrow \P)$ for $i\in \Ob(\P)$, where $\P$ is a
bounded above graded poset. In this section we define global
covering families by subsets of the whole category $\P$, imitating
some of the local features of the (local) covering families.

\begin{Defi}\label{defi_global covering family}
Let $\P$ be a bounded above graded poset for which there exists an
adequate covering family, and consider the sequence of functors
$F_0=c_\Z$, $F_1$, $F_2$,... given by Proposition \ref{cohomology
with adequate covering family}. A \emph{global covering family} is
a family of subsets $\K=\{K_p\}_{p\geq 0}$ with $K_p\subseteq
\Ob_p(\P)$ such that for $p\geq 1$ the restriction map
$$
\liminv F_p\rightarrow \prod_{i\in K_p} F_p(i)
$$
is a monomorphism and the map
$$
\prod_{i\in \Ob_{p-1}(\P)\setminus K_{p-1}} F_{p-1}(i)\rightarrow
\prod_{i\in K_p} F_p(i)
$$
is pure. For $p=0$ we require that the first map above is a pure
monomorphism.
\end{Defi}

Notice that the definition does not depend on the particular
adequate covering family used to obtain the sequence of functors
(as the sequence of functors does not depend on it). The map
$$
\prod_{i\in \Ob_{p-1}(\P)\setminus K_{p-1}} F_{p-1}(i)\rightarrow
\prod_{i\in K_p} F(i)
$$
used in the definition is obtained from
(\ref{eqn_long_exact_sequence_H^p}) as:
$$
\prod_{i\in \Ob_{p-1}(\P)\setminus K_{p-1}} F_{p-1}(i)
\hookrightarrow\prod_{i\in \Ob_{p-1}(\P)}
F_{p-1}(i)=\ker_{F_{p-1}}'\stackrel{\omega_p}\rightarrow \liminv
F_p\rightarrow \prod_{i\in K_p} F_p(i)
$$

For the rest of the section suppose we have a graded poset $\P$
with adequate covering family $\J$ and global covering family
$\K$. Consider the free, condensed and $\J$-determined functors
$c_\Z$, $F_1$, $F_2$,... given by Proposition \ref{cohomology with
adequate covering family}.

For each $p\geq 1$ we have, at the local level, a short exact
sequence
$$
0\rightarrow F_{p-1}(i_0)\stackrel{\lambda_{i_0}}\rightarrow
\prod_{i\in (i_0\downarrow \P)_p} F_{p-1}(i)
\stackrel{\pi_{i_0}}\rightarrow F_p(i_0)\rightarrow 0
$$
for $deg(i_0)\geq p$, and at the global level we have the short
exact sequence
$$
\xymatrix{0\ar[r]&\liminv F_{p-1}\ar[r]^{\iota_{p-1}}&\liminv
\ker_{F_{p-1}}'\ar[r]^{\omega_p}& \liminv F_p\ar[r]&H^p(\P;\Z)
\ar[r]&0 }.
$$

Recall that it was a set of integral equations involving $\J$ (the
adequateness property) that gave rise to the local short exact
sequences above. Next we study if the condition $H^p(\P;\Z)=0$ can
also be stated in terms of integral equations involving $\K$.

We begin defining adequate global covering families:

\begin{Defi}\label{defi global_adequate_covering_family}
Let $\P$ be a bounded above graded poset for which exists an
adequate covering family and let $\K$ be a global covering family.
We say that $\K$ is \emph{adequate} if
$$
\sum_{i\in \Ob_{p-1}(\P)} R^i_{p-1} =\sum_{i\in K_{p-1}}
R^i_{p-1}+\sum_{i\in K_p} R^i_p
$$
for $p\geq 1$ (the numbers $R^i_p$ depend on the local shape of
$\P$ and were defined in Definition \ref{defi_structural numbers
of graded poset}).
\end{Defi}

Consider the diagram
$$
\xymatrix{0\ar[r]&\liminv
F_{p-1}\ar[d]\ar[r]^{\iota_{p-1}}&\liminv
\ker_{F_{p-1}}'\ar[r]^{\omega_p}& \liminv
F_p\ar[r]\ar[d]&H^p(\P;\Z)
\ar[r]&0. \\
& \prod_{i\in K_{p-1}} F_{p-1}(i) & & \prod_{i\in K_p} F_p(i)}
$$

Recall that the six abelian groups involved in the diagram are
free abelian groups but, possibly, $H^p(\P;\Z)$.

Suppose first that $H^p(\P;\Z)=0$ for a fixed $p\geq 1$. Then, as
$\liminv \ker_{F_{p-1}}' \cong\prod_{i\in \Ob_{p-1}(\P)}
F_{p-1}(i)$ and $\dim F_{p-1}(i)=R^i_{p-1}$ (Proposition
\ref{cohomology with adequate covering family}), we have the
equality
$$
\dim\liminv F_{p-1}-\sum_{i\in \Ob_{p-1}(\P)} R^i_{p-1}
+\dim\liminv F_p=0.
$$
By definition of adequate covering family we also have the
equation
$$
\sum_{i\in K_{p-1}}R^i_{p-1}-\sum_{i\in \Ob_{p-1}(\P)} R^i_{p-1}
+\sum_{i\in K_p} R^i_p=0.
$$
From these two equations it is clear that
$$
\sum_{i\in K_{p-1}}R^i_{p-1}=\dim\liminv F_{p-1}\Leftrightarrow
\sum_{i\in K_p} R^i_p=\dim\liminv F_p.
$$
Thus, in case $H^p(\P;\Z)=0$ for $p\geq 1$, we obtain that the
condition
\begin{equation}\label{eqn dimliminvF_p=|K_p|}
\sum_{i\in K_p} R^i_p=\dim\liminv F_p
\end{equation}
holds for each $p\geq 0$ if and only if it holds for some $p\geq
1$.

Now we work in the opposite way: does Equation (\ref{eqn
dimliminvF_p=|K_p|}) for each $p\geq 1$ imply that $H^p(\P;\Z)=0$
for each $p\geq 1$? Recall that, by the definition of global
covering family, the map
$$
\liminv F_0\rightarrow \prod_{i\in K_0} F_0(i)
$$
is a pure monomorphism. Then Equation (\ref{eqn
dimliminvF_p=|K_p|}) for $p=0$ implies that this map is in fact an
isomorphism. Suppose now, inductively, that the map
$$
\liminv F_{p-1}\rightarrow \prod_{i\in K_{p-1}} F_{p-1}(i)
$$
is an isomorphism and that Equation (\ref{eqn dimliminvF_p=|K_p|})
holds for $p$. Because $\K$ is a global covering family the map
$$
\liminv F_p\rightarrow \prod_{i\in K_p} F_p(i)
$$
is a monomorphism. This implies that the quotient map
$$
\liminv \ker_{F_{p-1}}'/\liminv F_{p-1} \rightarrow \prod_{i\in
K_p} F_p(i)
$$
is injective too. By the induction hypothesis
$$
\liminv F_{p-1}\cong\prod_{i\in K_{p-1}} F_{p-1}(i)
$$
and thus
$$
\liminv \ker_{F_{p-1}}'/\liminv F_{p-1}= \prod_{i\in
\Ob_{p-1}(\P)\setminus K_{p-1}} F_{p-1}(i).
$$
As $\K$ is a global covering family the map
$$
\prod_{i\in \Ob_{p-1}(\P)\setminus K_{p-1}} F_{p-1}(i)\rightarrow
\prod_{i\in K_p} F(i)
$$
is pure with $\dim\liminv \ker_{F_{p-1}}'/\liminv
F_{p-1}=\sum_{i\in K_p} R^i_p=\dim\prod_{i\in K_p} F(i)$ ($\K$ is
adequate). This implies that it is an isomorphism, that
$$
\liminv \ker_{F_{p-1}}'/\liminv F_{p-1} \rightarrow \prod_{i\in
K_p} F_p(i)
$$
is an isomorphism too, that
$$\liminv \ker_{F_{p-1}}'\rightarrow \prod_{i\in K_p} F_p(i)$$
and
$$\liminv \ker_{F_{p-1}}'\rightarrow \liminv F_p$$
are epimorphisms and that
$$\liminv F_p\rightarrow \prod_{i\in K_p} F_p(i)$$
is an isomorphism. Moreover, $H^p(\P;\Z)=0$. We have proven

\begin{Thm}
Let $\P$ be a bounded above graded poset for which there exists an
adequate covering family and an adequate global covering family
$\K$. Suppose there exists $p_0\geq 1$ such that
$$
\dim\liminv F_{p_0}=\sum_{i\in K_{p_0}} R^i_{p_0}.
$$
Then
$$
\textit{$H^p(\P;\Z)=0$ for $p\geq 1\Leftrightarrow \dim\liminv
F_p=\sum_{i\in K_p} R^i_p$ for $p\geq 0$}.
$$
\end{Thm}

The beauty of this theorem is that it states the fact of $\P$
being acyclic in terms of integral equations. At first glance it
seems that the numbers $\dim\liminv F_p$ are totally unknown.
However, recall that $\dim\liminv F_0=\dim\liminv c_\Z$ is the
number of connected components of $\P$ (see also Lemma
\ref{lemma_liminvFp p+2_esqueleto}) and, if $\P$ is bounded below
with $p_0\in \Z$ minimal such that $\Ob_{p_0}(\P)\neq \emptyset$,
it turns out from the definition of global covering family that
$K_{p_0}=\Ob_{p_0}(\P)$ and $\dim\liminv F_{p_0}=\sum_{i\in
\Ob_{p_0}(\P)} R^i_{p_0}=\sum_{i\in K_{p_0}} R^i_{p_0}$. Thus, as
$R^i_0=1$ for each $i\in \Ob(\P)$, we obtain

\begin{Thm}\label{Thm_adequate_local y global characterization acyclicity}
Let $\P$ be a bounded graded poset for which there exist an
adequate covering family and an adequate global covering family
$\K$. Then $\P$ is acyclic if and only if $|K_0|$ equals the
number of connected components of $\P$. Moreover, in this case
$H^0(\P;\Z)=\Z^{|K_0|}$.
\end{Thm}



\section{Simplex-like posets.}\label{section_locally delta}
Recall (Section \ref{section simplicial complexes}) that
simplicial complexes can be viewed as a special kind of posets:
simplex-like posets. Another examples of simplex-like posets are
subdivision categories (see \cite{markus}) and the category
$\overline{s}dC$ used by Libman in \cite{libman} for the
normalizer decomposition for $p$-local finite groups. In this
section we shall see that if $\P$ is a simplex-like poset then it
is a graded poset, and $\P^{op}$ is a bounded above graded poset
that can be equipped with an adequate covering family. We begin
showing that $\P$ and $\P^{op}$ are graded:

\begin{Lem}\label{simplexlike_is_graded}
Let $\P$ be a simplex-like poset. Then $\P$ and $\P^{op}$ are
graded posets, $\P$ is bounded below and $\P^{op}$ is bounded
above.
\end{Lem}
\begin{proof}
That $\P$ and $\P^{op}$ are posets is immediate. To see that $\P$
is graded recall that, by definition of simplex-like poset, for
any $p\in \Ob(\P)$ the subcategory $(\P\downarrow p)$ is
isomorphic to the poset of all non-empty subsets of a finite set
$T$ (with inclusion as order relation). Define $deg(p)=|T|-1$.
Then $deg:\Ob(\P)\rightarrow \Z$ is a decreasing degree function
which assigns $0$ to minimal elements and $\P$ is graded. The same
function $deg$ is an increasing degree function for $\P^{op}$
which assigns $0$ to maximal elements. To see that $\P$ and
$\P^{op}$ are bounded just apply the definition.
\end{proof}

It is straightforward that the poset $\P$ is simplex-like $\P$ if
and only if for each $p\in \Ob(\P)$ the subcategory $(p\downarrow
\P^{op})=(\P\downarrow p)^{op}$ is isomorphic for some $n\geq 0$
to the poset $\Delta_{deg(p)}$ defined below:

\begin{Defi}\label{defi_Delta_n} For any integer $n\geq 0$
define the poset $\Delta_n$ as the subdivision category of the
poset $0\rightarrow 1\rightarrow 2...\rightarrow n$, i.e., it has
as objects the sequences of integers $[n_0,n_1,..,n_k]$ which are
increasing $n_0<n_1<..<n_k$, and with $0\leq n_j\leq n$. The
arrows are generated by the morphisms
$$d_l:[n_0,n_1,..,n_k]\rightarrow [n_0,n_1,..,\hat{n_l},..,n_k]$$
for $l=0,..,k$. It is a bounded above and below graded poset with
degree function $deg(n_0,n_1,..,n_k)=k$.
\end{Defi}
\begin{Ex} $\Delta_2$ looks like
$$
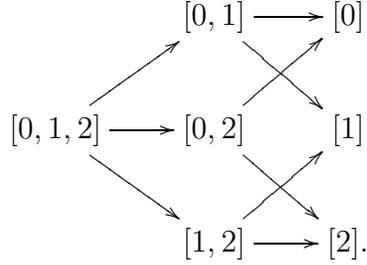
\xymatrix{ & [0,1]\ar[r]\ar[rd] & [0]\\
[0,1,2]\ar[r]\ar[ru]\ar[rd] & [0,2]\ar[ru]\ar[rd] & [1]\\
 & [1,2]\ar[ru]\ar[r] & [2].  }
$$
\end{Ex}

Next, we shall see that there is an adequate covering family for
$\Delta_n$. This is a first step to equip the whole $\P^{op}$ with
an adequate covering family.

\begin{Lem}\label{lemma_covering_for_Delta_n}
There is an adequate covering family $\J=\{J^{i}_p\}_{i\in
\Ob(\Delta_n)\textit{, $0\leq p\leq deg(i)$}}$ for $\Delta_n$.
Moreover, for $p\geq 0$ and $deg(i)=q\geq p$ we have
$$
R^{i}_p=\sum_{l=0}^p (-1)^{(p-l)}\big(\!\begin{smallmatrix} q+1
\\ l\end{smallmatrix}\!\big)
$$
\end{Lem}
\begin{proof}
For the object $i_0=[0,1,2,..,n]\in \Ob(\Delta_n)$ of degree $n$
define the subsets $J^{i_0}_p\subseteq (i_0\downarrow \Delta_n)_p$
for $0\leq p\leq n$ as the sequences of increasing integers
$[n_0,n_1,..,n_{p-1},n_p]$ which biggest value is $n$, i.e., such
that $n_p=n$. This definition coincides with the following
inductive one:
$$
J^{i_0}_0=\{[n]\}
$$
and
$$
J^{i_0}_p=\{[s,n]\textit{, $s\in (i_0\downarrow
\Delta_n)_{p-1}\setminus J^{i_0}_{p-1}$}\}
$$
for $p=1,2,..,n$.

For any other object $i\in \Ob(\Delta_n)$ different from
$i_0=[0,1,2,..,n]$ of degree $deg(i)=m<n$  we define the subsets
$J^i_p\subseteq (i\downarrow \Delta_n)_p$ for $0\leq p\leq m$
using the definition above for $\Delta_m$ instead of $\Delta_n$
and taking into account the natural isomorphism $\Delta_m\cong
(\Delta_n)^i$ (see the remark below). Purely combinatorial
arguments show that $\J$ is a covering family for $\Delta_n$.
\end{proof}
\begin{Ex}For $n=2$ have $J^{[0,1,2]}_{[0]}=\{[2]\}$,
$J^{[0,1,2]}_{[1]}=\{[0,2],[1,2]\}$,$J^{[0,1,2]}_{[2]}=\{[0,1,2]\}$:
$$
\xymatrix{ & [0,1]\ar[r]\ar[rd] & [0]\\
\mathbf{[0,1,2]}\ar[r]\ar[ru]\ar[rd] & \mathbf{[0,2]}\ar[ru]\ar[rd] & [1]\\
 & \mathbf{[1,2]}\ar[ru]\ar[r] & \mathbf{[2]}.  }
$$
Also, $J^{[0,1]}_{[0]}=\{[1]\}$, $J^{[0,2]}_{[0]}=\{[2]\}$ and
$J^{[1,2]}_{[0]}=\{[2]\}$
\end{Ex}
\begin{Rmk}\label{rmk_Rio_p=1}
For this covering family we have the following statement, which is
stronger than Definition
\ref{defi_covering_family}\ref{defi_covering_b}):
$J^i_q=J^{i_0}_q\cap (i\downarrow \Delta_n)_q$ for each $i_0$,
$0\leq p\leq deg(i_0)$, $i\in J^{i_0}_p$ and $0\leq q\leq deg(i)$.
Moreover, we have that
$$
R^{i}_p=1
$$
if $deg(i)=p$ by using the binomial expansion of $(1-1)^{p+1}=0$.
\end{Rmk}

\begin{Rmk}\label{rmk_Delta_m_isosdetermined + orderpreserverving}
Notice that any isomorphism of categories
$\varphi:\Delta_m\rightarrow \C$ is determined by the values
$\varphi([a])$ for $0\leq a\leq m$ ($[a]\in (\Delta_m)_0$). If
$\C=(\Delta_n)^i$ the natural isomorphism
$\varphi:\Delta_m\rightarrow (\Delta_n)^i$ used in the lemma above
is the only order preserving isomorphism, i.e., the only one such
that $a<b$ if and only if $\varphi([a])<\varphi([b])$.
\end{Rmk}

Now we reach the main result of this section:

\begin{Lem}\label{lem_coveringfamily_for_locallyDelta}
If $\P$ is a simplex-like poset then there is an adequate covering
family for the bounded above graded poset $\P^{op}$.
\end{Lem}
\begin{proof}
Above we have defined isomorphisms of categories $(p\downarrow
\P^{op})\simeq \Delta_n$ for each $p\in \Ob(\P)$, and we know that
$\Delta_n$ can be equipped with an adequate covering family. To
build an adequate covering family $\K=\{K^{i_0}_p\}_{i_0\in
\Ob(\P^{op})\textit{, $0\leq p\leq deg(i_0)$}}$ we have just to
choose appropriately the isomorphisms $(p\downarrow \P^{op})\simeq
\Delta_n$.

Consider the degree function $deg:\Ob(\P^{op})\rightarrow \Z$
defined in Lemma \ref{simplexlike_is_graded} and the set $T$ of
the maximal elements of $\P^{op}$, i.e., $T=\{p\in
\Ob(\P^{op})|deg(p)=0\}$. Choose a total order $<$ for $T$
(suppose $T$ is finite or use the Axiom of Choice \cite{jech}).
Then, given $p\in \Ob(\P^{op})$, consider the subset $(p\downarrow
\P^{op})_0\subseteq T$ and the restriction $((p\downarrow
\P^{op})_0,<)$ of the total order from $T$. There is a unique
isomorphism
$$
\varphi_p:(p\downarrow \P^{op})\simeq \Delta_{deg(p)}
$$
which induces an order preserving map
$$
((p\downarrow \P^{op})_0,<)\simeq
(\Delta_{deg(p)})_0=\{[0],[1],[2],...,[deg(p)]\}.
$$

Denote by $\J$ the covering family for $\Delta_{deg(p)}$ of Lemma
\ref{lemma_covering_for_Delta_n} and define, for $0\leq n\leq
deg(p)$,
$$
K^{p}_n=\varphi_{p}^{-1}(J^{\varphi(p)}_n).
$$

Then $\K$ fulfills condition \ref{defi_covering_a}) of Definition
\ref{defi_covering_family} because for $0\leq n<deg(p)$
$$
(p\downarrow
\P^{op})_n=\varphi_{p}^{-1}((\Delta_{deg(p)})_n)=\varphi_{p}^{-1}(\bigcup_{i\in
J^{\varphi_{p}(p)}_{n+1}} (i\downarrow \Delta_{deg(p)})_n)=
$$
$$=\bigcup_{\varphi_{p}^{-1}(i)\in \varphi_{p}^{-1}(K^{p}_{n+1})}
\varphi_{p}^{-1}((i\downarrow \Delta_{deg(p)})_n)=\bigcup_{i\in
K^{p}_{n+1}} (i\downarrow \P^{op})_n.
$$

To check condition \ref{defi_covering_b}) of Definition
\ref{defi_covering_family} take $i\in K^{p}_{n+1}$ for some $0\leq
n<deg(p)$ and call $\J'$ to the covering family for $\Delta_{n+1}$
of Lemma \ref{lemma_covering_for_Delta_n}. We want to see that
$K^i_n\subseteq K^{p}_n$. Recalling the natural inclusion
$(i\downarrow \P^{op})\subseteq (p\downarrow \P^{op})$ this is
equivalent to
$$
\varphi_i^{-1}(J'^{\varphi_i(i)}_n)\subseteq
\varphi_{p}^{-1}(J^{\varphi_{p}(p)}_n)
$$
and to
$$
\psi(J'^{\varphi_i(i)}_n)\subseteq J^{\varphi_{p}(p)}_n
$$
where $\psi=\varphi_{p}\circ \varphi_i^{-1}$. By construction
$\psi$ is order preserving (see Remark
\ref{rmk_Delta_m_isosdetermined + orderpreserverving}) and thus
this inclusion holds.
\end{proof}

Next, we list more familiar re-statements of results about
covering families applied to simplex-like posets:

\begin{Lem}\label{lema_adequate_covering_family_simplex-like posets}
Let $\P$ be a simplex-like poset and consider the bounded above
graded poset $\P^{op}$ (for which exists an adequate covering
family by Lemma {\rm\ref{lem_coveringfamily_for_locallyDelta}}).
Let $\K$ be a global covering family for $\P^{op}$. Then $\K$ is
\emph{adequate} if and only if
$$
|\Ob_{p-1}(\P)| =|K_{p-1}|+|K_p|
$$
for $p\geq 1$.
\end{Lem}
\begin{proof}
Apply Remark \ref{rmk_Rio_p=1} to Definition
\ref{defi_adequate_covering_family}.
\end{proof}
\begin{Lem}\label{lemma_analogue eulerchar_simplex-like posets}
Let $\P$ a simplex-like poset. Then
$$
\sum_i (-1)^i\dim H^i(\P;\Z)=\sum_i (-1)^i |\Ob_i(\P)|.
$$
\end{Lem}
\begin{proof}
Use Remark \ref{rmk_Rio_p=1}, and observe that $\P\simeq\P^{op}$
and that $\Ob(\P)=\Ob(\P^{op})$. Then apply Lemma
\ref{lemma_analogue eulerchar}.
\end{proof}

\chapter{Application: Webb's conjecture}\label{section_applications Webb_conjecture}
Denote by $\Ss_p(G)$ the Brown's complex of the finite group $G$
for the prime $p$, whose elements are the non-trivial
$p$-subgroups of $G$, and that was introduced by Brown
\cite{brown}. Webb conjectured that the orbit space $\Ss_p(G)/G$
(as topological space) is contractible. This conjecture was first
proven by Symonds in \cite{symonds}, generalized for blocks by
Barker \cite{barker1,barker2} and extended to arbitrary
(saturated) fusion system by Linckelmann \cite{markus}.

The works of Symonds and Linckelmann prove the contractibility of
the orbit space by showing that it is simply connected and
acyclic, and invoking Whitehead's Theorem. Both proofs of
acyclicity work on the subposet of normal chains (introduced by
Knorr and Robinson \cite{robinson} for groups). Symonds uses the
results from Thévenaz and Webb \cite{G-homot} that the subposet of
normal chains is $G$-equivalent to Brown's complex. Linckelmann
proves on his own that, also for fusion systems, the orbit space
and the orbit space on the normal chains has the same cohomology
\cite[Theorem 4.7]{markus}.

In this chapter we shall apply the results of Chapter
\ref{section_cohomology} to prove in an alternative way that the
orbit space on the normal chains is acyclic.

Let $(S,\Ff)$ be a saturated fusion system where $S$ is a
$p$-group. Consider its subdivision category $\Ss(\Ff)$ (see
\cite[2]{markus}) and the poset $[\Ss(\Ff)]$. An object in
$[\Ss(\Ff)]$ is an $\Ff$-isomorphism class of chains
$$
[Q_0<Q_1<...<Q_n]
$$
where the $Q_i$'s are subgroups of $S$. The subcategory
$([\Ss(\Ff)]\downarrow [Q_0<...<Q_n])$ has objects
$[Q_{i_0}<...<Q_{i_m}]$ with $0\leq m\leq n$ and $0\leq
i_1<i_2<...<i_m\leq n$ (see \cite[2]{markus} again). For example,
$([\Ss(\Ff)]\downarrow [Q_0<Q_1<Q_2])$ is
$$
\xymatrix{ [Q_0]\ar[r]\ar[rd] & [Q_0<Q_1]\ar[rd]\\
[Q_1]\ar[ru]\ar[rd]&[Q_0<Q_2]\ar[r]&[Q_0<Q_1<Q_2].\\
[Q_2]\ar[r]\ar[ru]&[Q_1<Q_2]\ar[ru] }
$$

Then it is clear that $[\Ss(\Ff)]$ is a simplex-like poset.
Following Linckelmann's notation denote by $\Ss_\lhd(\Ff)$ the
full subcategory of $\Ss(\Ff)$ which objects $Q_0<...<Q_n$ with
$Q_i\lhd Q_n$ for $i=0,...,n$. Also, denote by  $[\Ss_\lhd(\Ff)]$
the subdivision category of $\Ss_\lhd(\Ff)$, which is a sub-poset
of $[\Ss(\Ff)]$.

Our goal in this chapter is to prove that
$H^n([\Ss_\lhd(\Ff)];\Z)=0$ for $n\geq 1$ and
$H^0([\Ss_\lhd(\Ff)];\Z)=\Z$. It is straightforward that
$[\Ss_\lhd(\Ff)]$ is a simplex-like poset and thus, by Lemma
\ref{lem_coveringfamily_for_locallyDelta}, there exists an
adequate covering family for the bounded above graded poset
$[\Ss_\lhd(\Ff)]^{op}$. We shall build an adequate global covering
family for $[\Ss_\lhd(\Ff)]^{op}$ in order to apply Theorem
\ref{Thm_adequate_local y global characterization acyclicity}.

The definition of the global covering family is as follows, and it
is related with the pairing defined by Linckelmann in
\cite[Definition 4.7]{markus}. The notion of paired chains was
used by Knorr and Robinson in several forms throughout
\cite{robinson}.

\begin{Defi}\label{defi_global_covering_family for normal chains}
For the graded poset $[\Ss_\lhd(\Ff)]^{op}$ define the subsets
$\K=\{K_n\}_{n\geq 0}$ by
$$
K_n={\big\{}[Q_0<...<Q_n]\textrm{ $|$
}[Q_0<...<Q_n]=[Q'_0<...<Q'_n]\Rightarrow \cap_{i=0}^n
N_S(Q'_i)=Q'_n{\big\}}
$$
\end{Defi}

\begin{Lem}
The family $\K=\{K_n\}_{n\geq 0}$ defined in
$\ref{defi_global_covering_family for normal chains}$ is a global
covering family for $[\Ss_\lhd(\Ff)]^{op}$.
\end{Lem}

The proof is postponed to the next section. The argument uses the
properties of the (local) covering family of the simplex-like
category $[\Ss_\lhd(\Ff)]^{op}$ and that the chains $Q_0<...<Q_n$
can be ordered by $|Q_n|$. The fact that $\K$ is defined through a
pairing provides (see next section) a bijection
$\psi:\Ob_n([\Ss_\lhd(\Ff)]^{op})\setminus K_n\rightarrow
K_{n+1}$, which gives

\begin{Lem}
The global covering family $\K=\{K_n\}_{n\geq 0}$ defined in
$\ref{defi_global_covering_family for normal chains}$ is adequate
for $[\Ss_\lhd(\Ff)]^{op}$.
\end{Lem}
\begin{proof}
For any $n\geq 0$ we have
$$
\Ob_n([\Ss_\lhd(\Ff)]^{op})=K_n\cup
(\Ob_n([\Ss_\lhd(\Ff)]^{op})\setminus K_n).
$$
Then the bijection (Lemma \ref{lemma_bijection_K_n})
$$
\Ob_n([\Ss_\lhd(\Ff)]^{op})\setminus K_n\rightarrow K_{n+1}
$$
gives
$$
|\Ob_n([\Ss_\lhd(\Ff)]^{op})|=|K_n|+|K_{n+1}|
$$
as wished (see Lemma
\ref{lema_adequate_covering_family_simplex-like posets}).
\end{proof}


As $[\Ss(\Ff)]$ is connected then so is $[\Ss_\lhd(\Ff)]^{op}$.
Thus, $\dim \liminv c_Z=1$ over $[\Ss_\lhd(\Ff)]^{op}$. Also, by
elementary properties of $p$-groups it is clear that
$K_0=\{[S]\}$, i.e., the only subgroup of $S$ which equals its
normalizer in $S$ is $S$ itself, and $|K_0|=1$. Then Theorem
\ref{Thm_adequate_local y global characterization acyclicity}
gives
\begin{Thm}
Let $(S,\Ff)$ be a saturated fusion system. Then
$H^n([\Ss_\lhd(\Ff)]^{op};\Z)=0$ for $n\geq 1$ and
$H^0([\Ss_\lhd(\Ff)]^{op};\Z)=\Z$.
\end{Thm}

\section{\K\ is an adequate global covering
family.} In this section we prove that the family
$\K=\{K_n\}_{n\geq 0}$ defined in \ref{defi_global_covering_family
for normal chains} is an adequate global covering family for
$[\Ss_\lhd(\Ff)]^{op}$. We use terminology and results from
\cite[Appendix]{blo2}.

For any chain $Q_0<...<Q_n$ in $\Ss_\lhd(\Ff)$ define the
following subgroup of automorphisms of $Q_n$
$$
A_{Q_0<...<Q_n}=\{\alpha\in \Aut(Q_n)|\alpha(Q_i)=Q_i,
i=0,...,n\}.
$$

Then,
$$
N^{A_{Q_0<...<Q_n}}_S(Q_n)=\cap_{i=0}^n N_S(Q_i).
$$
If $[Q_0<...<Q_n]=[Q'_0<...<Q'_n]$ then there is $\varphi\in
\Iso_\Ff(Q_n,Q'_n)$ with $Q'_i=\varphi(Q_i)$ for $i=0,...,n$ and
$$
\varphi A_{Q_0<...<Q_n} \varphi^{-1} = A_{Q'_0<...<Q'_n}.
$$
By \cite[A.2(a)]{blo2} $Q_n$ is fully $A_{Q_0<...<Q_n}$-normalized
if and only if $|N^{A_{Q_0<...<Q_n}}_S(Q_n)|$ is maximum among
$|N^{A_{Q'_0<...<Q'_n}}_S(Q_n)|$ with
$[Q'_0<...<Q'_n]=[Q_0<...<Q_n]$, i.e., if and only if
$|\cap_{i=0}^n N_S(Q_i)|$ is maximum among $|\cap_{i=0}^n
N_S(Q'_i)|$ with $[Q'_0<...<Q'_n]=[Q_0<...<Q_n]$. Notice that in
the isomorphism class of chains $[Q_0<...<Q_n]$ always there is a
representantive $Q'_0<...<Q'_n$ which is fully
$A_{Q'_0<...<Q'_n}$-normalized, and that any two representantives
$Q'_0<...<Q'_n$ and $Q''_0<...<Q''_n$ of $[Q_0<...<Q_n]$ which are
fully $A_{Q'_0<...<Q'_n}$-normalized and fully
$A_{Q''_0<...<Q''_n}$-normalized respectively verify
$$
|\cap_{i=0}^n N_S(Q'_i)|=|\cap_{i=0}^n N_S(Q''_i)|.
$$

Thus, Definition \ref{defi_global_covering_family for normal
chains} is equivalent to
\begin{Defi}
For the graded poset $[\Ss_\lhd(\Ff)]^{op}$ define the subsets
$\K=\{K_n\}_{n\geq 0}$ by
$$
K_n=\{[Q'_0<...<Q'_n]\text{ $|$ $Q'_n$ fully
$A_{Q'_0<...<Q'_n}$-normalized and $\cap_{i=0}^{n}
N_S(Q'_i)=Q'_n$}\}
$$
\end{Defi}

\begin{Lem}\label{lemma_bijection_K_n}
For any $n\geq 0$ there is a bijection
$$
\psi:\Ob_n([\Ss_\lhd(\Ff)]^{op})\setminus K_n\rightarrow K_{n+1}
$$
\end{Lem}
\begin{proof}
Take $[Q_0<...<Q_n]\in \Ob_n([\Ss_\lhd(\Ff)]^{op})\setminus K_n$
and a representantive $Q'_0<...<Q'_n$ which is fully
$A_{Q'_0<...<Q'_n}$-normalized. Then $\cap_{i=0}^{n}
N_S(Q'_i)>Q'_n$. Define
$\psi([Q_0<...<Q_n])=[Q'_0<...<Q'_n<\cap_{i=0}^{n} N_S(Q'_i)]$.
The proof is divided in four steps:

\textbf{a) $\psi$ is well defined.} Take another representantive
$Q''_0<...<Q''_n$ which is fully $A_{Q''_0<...<Q''_n}$-normalized.
Then, by \cite[A.2(c)]{blo2}, there is a morphism
$$
\varphi\in
\Hom_\Ff(N_S^{A_{Q'_0<...<Q'_n}}(Q'_n),N^{A_{Q''_0<...<Q''_n}}_S(Q'_n))
$$
with $\varphi(Q'_i)=Q''_i$ for $i=0,..n$. As $|\cap_{i=0}^{n}
N_S(Q'_i)|=|\cap_{i=0}^{n} N_S(Q''_i)|$ then $\varphi$ is an
isomorphism onto $N^{A_{Q''_0<...<Q''_n}}_S(Q'_n)$ and thus
$$
[Q'_0<...<Q'_n<\cap_{i=0}^{n}
N_S(Q'_i)]=[Q''_0<...<Q''_n<\cap_{i=0}^{n} N_S(Q''_i)].
$$

\textbf{b) $\psi([Q_0<...<Q_n])$ belongs to $K_{n+1}$.} We have
$\psi([Q_0<...<Q_n])=[Q'_0<...<Q'_n<\cap_{i=0}^{n} N_S(Q'_i)]$
where $Q'_0<...<Q'_n$ is fully $A_{Q'_0<...<Q'_n}$-normalized.
Take any representantive $Q''_0<...<Q''_n<Q''_{n+1}$ in
$[Q'_0<...<Q'_n<\cap_{i=0}^{n} N_S(Q'_i)]$. If it were the case
that $Q''_{n+1}< \cap_{i=0}^{n+1} N_S(Q''_i)$ then we would have
$$
\cap_{i=0}^{n} N_S(Q'_i)\cong Q''_{n+1} < \cap_{i=0}^{n+1}
N_S(Q''_i)\leq  \cap_{i=0}^{n} N_S(Q''_i),
$$
which is in contradiction with $Q'_0<...<Q'_n$ being fully
$A_{Q'_0<...<Q'_n}$-normalized.

\textbf{c) $\psi$ is injective.} Suppose we have $[Q_0<...<Q_n]$
and $[R_0<..<R_n]$ with
$$
[Q'_0<...<Q'_n<\cap_{i=0}^{n}
N_S(Q'_i)]=[R'_0<...<R'_n<\cap_{i=0}^{n} N_S(R'_i)].
$$
Then
$[R_0<...<R_n]=[R'_0<...<R'_n]=[Q'_0<...<Q'_n]=[Q_0<...<Q_n]$.

\textbf{d) $\psi$ is surjective.} Take $[Q_0<...<Q_n<Q_{n+1}]$ in
$K_{n+1}$. We check that
$$\psi([Q_0<...<Q_n])=[Q_0<...<Q_n<Q_{n+1}]$$
Take a representantive $Q'_0<...<Q'_n$ in $[Q_0<...<Q_n]$ which is
fully $A_{Q'_0<...<Q'_n}$-normalized. Then $[Q_0<...<Q_n]\in
\Ob_n([\Ss_\lhd(\Ff)]^{op})\setminus K_n$ and
$\psi([Q_0<...<Q_n])=[Q'_0<...<Q'_n<\cap_{i=0}^{n} N_S(Q'_i)]$.
Then, by \cite[A.2(c)]{blo2}, there is $$\varphi\in
\Hom_\Ff(N_S^{A_{Q_0<...<Q_n}}(Q_n),N^{A_{Q'_0<...<Q'_n}}_S(Q_n))
$$
with $\varphi(Q_i)=Q'_i$ for $i=0,..n$. As
$Q_{n+1}=\cap_{i=0}^{i=n+1} N_S(Q_i)$ then $Q_{n+1}\leq
\cap_{i=0}^{n} N_S(Q_i)$ and $\varphi(Q_{n+1})\leq \cap_{i=0}^{n}
N_S(Q'_i)$. If it were the case that $\varphi(Q_{n+1})<
\cap_{i=0}^{n} N_S(Q'_i)$ then we would have
$$
\varphi(Q_{n+1})< N_{\cap_{i=0}^{n}
N_S(Q'_i)}(\varphi(Q_{n+1}))=\cap_{i=0}^{n+1}
N_S(\varphi(Q_i))=\varphi(\cap_{i=0}^{n+1}
N_S(Q_i))=\varphi(Q_{n+1}),
$$
a contradiction. Thus, $\varphi(Q_{n+1})=\cap_{i=0}^{n} N_S(Q'_i)$
and the proof is finished.
\end{proof}

\begin{Lem}
The family $\K=\{K_n\}_{n\geq 0}$ defined in
$\ref{defi_global_covering_family for normal chains}$ is a global
covering family for $[\Ss_\lhd(\Ff)]^{op}$.
\end{Lem}
\begin{proof}
For $n=0$ we have to prove that the map
$$
\liminv c_Z\rightarrow  \prod_{i\in K_0} F_0(i)=\prod_{i\in
\{[S]\}} F_0(i)=\Z
$$
is a pure monomorphism. In fact, as $[\Ss_\lhd(\Ff)]$ is
connected, this map is an isomorphism.

For $n\geq 1$ we have to prove that
$$
\liminv F_n\rightarrow \prod_{i\in K_n} F_n(i)
$$
is a monomorphism and the map
$$
\prod_{i\in \Ob_{n-1}(\P)\setminus K_{n-1}} F_{n-1}(i)\rightarrow
\prod_{i\in K_n} F_n(i)
$$
is pure. We begin proving the injectivity. Take $\psi\in \liminv
F_p$ such that $\psi(i)=0$ for each $i\in K_n$. If there is no
object of degree greater than $n$, i.e.
$\Ob_{\{...,n+2,n+1\}}([\Ss_\lhd(\Ff)]^{op})=\emptyset$ then
$K_n=\Ob_n([\Ss_\lhd(\Ff)]^{op})$ and we are done. If not, we
prove that $\psi(j)=0$ for each $j$ of degree $n+1$ by induction
on $|Q_{n+1}|$. This is enough to see that $\psi$ is zero as $F_n$
is $n$-condensed. We shall use the (local) covering family $\J$
defined in \ref{section_locally delta} for the simplex-like
category $[\Ss_\lhd(\Ff)]^{op}$.

The case base is $j=[Q_0<...<Q_{n+1}]$ with $|Q_{n+1}|$ maximal.
This implies that $J^j_n\subseteq K_n$. Then $\psi(j)$ goes to
zero by the monomorphism
$$
F_n(j)\rightarrow \prod_{i\in J^{j}_n} F_n(i),
$$
and thus $\psi(j)=0$. For the induction step consider
$j=[Q_0<...<Q_{n+1}]$ and $j'=[Q_0<...<\widehat{Q_l}<...<Q_{n+1}]
\in J^j_n$ with $0\leq l<n$. Then, either $j'\in K_n$ and
$\psi(j')=0$, or $j'\notin K_n$ and there is an arrow in
$[\Ss_\lhd(\Ff)]^{op}$
$$
j''=[Q'_0<...<\widehat{Q'_l}<...<Q'_{n+1}< \cap_{i=0,i\neq l}^{m}
N_S(Q'_i)]\rightarrow j'=[Q_0<...<\widehat{Q_l}<...<Q_{n+1}].
$$
In the latter case $\psi(j'')=0$ by the induction hypothesis, and
thus $\psi(j')=0$ too. As before, as the map $F_n(j)\rightarrow
\prod_{i\in J^{j}_n} F_n(i)$ is a monomorphism, $\psi(j)=0$.

Now we prove that the map
$$
\omega: \prod_{i\in \Ob_{n-1}(\P)\setminus K_{n-1}}
F_{n-1}(i)\rightarrow \prod_{i\in K_n} F_n(i)
$$
is pure. Take $y\in \prod_{i\in K_n} F_n(i)$, $n\geq 1$ and $x\in
\prod_{i\in \Ob_{n-1}(\P)\setminus K_{n-1}} F_{n-1}(i)$ with
$$
n\cdot y = \omega(x).
$$
We want to find $x'$ with $n\cdot x'=x$. We prove that $x_i$ is
divisible by $n$ for each $i=[Q_0<...<Q_{n-1}]\in
\Ob_{n-1}(\P)\setminus K_{n-1}$ by induction on $|Q_{n-1}|$.

The case base is $i=[Q_0<...<Q_{n-1}]$ with $|Q_{n-1}|$ maximal.
Consider the arrow in $[\Ss_\lhd(\Ff)]^{op}$
$$
j=[Q'_0<...<Q'_{n-1}< \cap_{i=0}^{n-1} N_S(Q'_i)]\rightarrow
[Q_0<...<Q_{n-1}].
$$
As $Q_{n-1}$ is maximal then $J^j_{n-1}\subseteq K_{n-1}$. Then
$n\cdot y_j=\omega(x)_j$ is the image of $(x_i,0,...,0)$ by the
map
$$
\ker_{F_{n-1}}'(j)=\prod_{l\in (j\downarrow
{[\Ss_{\lhd}(\Ff)]})_{n-1}} F_{n-1}(i)\stackrel{\pi_j}\rightarrow
F_n(j).
$$
As $F_n(j)\cong \prod_{l\in (j\downarrow {[\Ss_{\lhd}(\Ff)]}
)_{n-1}\setminus J^j_{n-1}} F_{n-1}(l)=F_{n-1}(i)$ by Remark
\ref{rmk_section} then $n$ divides $x_i$.

For the induction step consider $i=[Q_0<...<Q_{n-1}]\in
\Ob_{n-1}(\P)\setminus K_{n-1}$ and $j=[Q'_0<...<Q'_{n-1}<
\cap_{i=0}^{i=n-1} N_S(Q'_i)]$. As before, $n\cdot
y_j=\omega(x)_j$ is the image of $\tilde{x}=x|{(j\downarrow
{[\Ss_{\lhd}(\Ff)]})_{n-1}}$ by the map
$$
\ker_{F_{n-1}}'(j)=\prod_{l\in (j\downarrow
{[\Ss_{\lhd}(\Ff)]})_{n-1}} F_{n-1}(i)\stackrel{\pi_j}\rightarrow
F_n(j).
$$
By Remark \ref{rmk_section} again,
$$
n\cdot y_j=\pi_j(\tilde{x}-(\lambda_j\circ
s_j)(\tilde{x}))=\pi_j((x_i-(\lambda_j\circ
s_j)(x|{J^j_n}),0,...,0)).
$$
Now, by the induction hypothesis, for each $l\in J^j_n$ either
$l\in K_{n-1}$ and $x_l=0$, either $l\notin K_{n-1}$ and thus, by
the induction hypothesis, $n$ divides $x_l$. Then $n$ divides
$x|{J^j_n}$ and so, by the equation above and the isomorphism
$F_n(j)\cong F_{n-1}(i)$, $n$ divides $x_i$ too.
\end{proof}

\chapter{Application: homotopy colimit}\label{hocolim} In
this section we deal with the problem of when the natural map from
the homotopy colimit of a diagram of nerves of groups to the nerve
of the colimit of the groups is an isomorphism. By $\hocolim$ and
$\hocolim_*$ we denote the unpointed and pointed homotopy colimits
in the sense of Bousfield-Kan respectively. For any discrete group
$G_0\in \Grp$  we denote by $BG_0\in SSet$ the nerve of the
category with one object and with automorphism group $G_0$, which
is a classifying space for $G_0$, and by $\I G_0\subseteq \Z[G_0]$
the augmentation ideal of $G_0$. $\P$ always denote a small
category in this section.

More precisely, for a functor $G:\P\rightarrow \Grp$ with takes as
values discrete groups, and a cone
$$\tau: G\Rightarrow G_0,$$
we have and induced cone
$$B\tau:BG \Rightarrow BG_0$$
which gives a map
$$\hocolim BG\rightarrow BG_0.$$
We are interested in finding conditions such that this map induces
a weak homotopy equivalence
\begin{equation}\label{equation hocolimN=Ncolim}
\hocolim BG\cong BG_0
\end{equation}
in case $G_0=\limdir G$. We have the following preliminary result,
which proof is left to Section \ref{section_hocolim_thm}:

\begin{Thm}\label{hocolim_technic}
Let $G:\P\rightarrow \Grp$ be a functor and $\tau: G\Rightarrow
G_0$ be a cone, and call $F$ the homotopy fiber
$$F\rightarrow \hocolim BG\rightarrow BG_0.$$
Assume $\P$ is contractible and the cone $\tau:G\Rightarrow G_0$
is monic. Then
\begin{enumerate}
\item $\pi_1(F)=\ker(\limdir G\rightarrow G_0)$ \item
$\pi_0(F)=\coker(\limdir G\rightarrow G_0)$ \item
$H_j(F)=\limdir_{j-1} H$ for each $j\geq 2$,  where
$H:\P\rightarrow \Ab$ is a monic functor.
\end{enumerate}
\end{Thm}

Thus, if $G_0=\limdir G$, then $F$ is simply connected. Some
examples follow:

\begin{Ex}
\label{ex_hocolim_telescope} Consider a monic functor
$G:\P\rightarrow \Grp$ where $\P$ is the ``telescope category"
$\P$ with shape
$$ \xymatrix{a_0 \ar[r]^{f_1} & a_1 \ar[r]^{f_2} & a_2 \ar[r]^{f_3} &  a_3 \ar[r]^{f_4} & a_4 ...}$$
Then $\P$ is contractible and the cone $G\Rightarrow
G_0=\bigcup_{l\geq 0} G(a_l)$ is monic. Then by part Theorem
\ref{hocolim_technic} and Example \ref{ex_pushout_telescope
limdir_monic acyclic}
$$\hocolim BG\cong BG_0.$$
\end{Ex}

\begin{Ex}
\label{ex_hocolim_pushout} Consider a functor $G:\P\rightarrow
\Grp$ where $\P$ is the ``pushout category" $\P$ with shape
$$\xymatrix {a \ar[r]^{f}\ar[d]^{g} & b \\ c.}$$
Then $\P$ is contractible and, if $G(f)$ and $G(g)$ are injective,
then the cone $\tau:G\Rightarrow \limdir G$ is monic \cite{trees}.
Then, by Theorem \ref{hocolim_technic} and Example
\ref{ex_pushout_telescope limdir_monic acyclic} we have
$$\hocolim BG\cong BG_0.$$

Thus we obtain  the classical result of Whitehead that states that
$$\hocolim BG\cong B(G(b)*_{G(a)}G(c))$$
if both $G(f)$ and $G(g)$ are monomorphisms.

\begin{Rmk}
Assume $G(f)$ is a monomorphism. Then \cite{BK} $BG(f)$ is a
cofibration, the diagram $BG(b)\leftarrow BG(a)\rightarrow BG(c)$
is a cofibrant object for some structure of  closed model category
on $\SSet^\P$ and so
$$\hocolim BG=\limdir BG.$$
This has the following consequence
\begin{Claim}
If $G(f)$ and $G(g)$ are monomorphisms then
$$\hocolim BG\cong B(\limdir(G(b)\leftarrow G(a)\rightarrow G(c)))\cong \limdir BG.$$
\end{Claim}
Recall that the $n$-simplices in $\limdir BG$ are given by
$$(\limdir BG)_n=\limdir_{Set} BG_n.$$
This is a finite set if the groups $G(b)$ and $G(c)$ are, but the
$n$-simplices in $B(\limdir(G(b)\leftarrow G(a)\rightarrow G(c)))$
are in general infinite as $\limdir(G(b)\leftarrow G(a)\rightarrow
G(c))$ is in general infinite.

Thus the geometric realization of $\limdir BG$ gives a
dimension-wise finite $CW$-complex as classifying space for the
(possibly infinite) group $G(b)*_{G(a)}G(c)$.
\end{Rmk}

\end{Ex}
\begin{Ex}
Consider a monic functor $G:\P\rightarrow \Grp$ where $\P$ is a
graded poset which is a tree. Then $\P$ is contractible and the
cone $\tau:G\Rightarrow \limdir G$ is monic \cite{trees}. Then, by
Theorem \ref{hocolim_technic} and Corollary \ref{Ex_monic on tree
acyclic direct} we have
$$\hocolim BG\cong B\limdir G.$$
\end{Ex}
\begin{Ex}\label{ex_filtered_contractible}
Recall that a filtered category is a category $\P$ such that any
two objects $i$,$j\in \Ob(\P)$ can be joined
$$
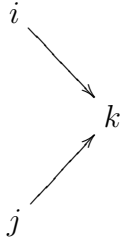
\xymatrix{i\ar[rd]&\\
&k\\
j\ar[ru]&}
$$
and such that any two parallel arrows
 $u,v:i\rightarrow j$ can be
co-equalized by $w:j\rightarrow k$ with $wu=wv$.

Any functor $F:\P\rightarrow \Ab$ with $\P$ filtered is
$\limdir$-acyclic as $\limdir:\Ab^\P\rightarrow \Ab$ is exact (see
\cite[2.6.15]{weibel}).

A finite poset which is filtered has a terminal object, and so
direct limits become trivial. The notion of pseudo-injectivity
gives an alternative condition to being filtered that also implies
$\limdir$-acyclicity but that does not impose the existence of a
terminal object when $\P$ is a finite poset.

Consider a contractible filtered category $\P$ and a monic functor
$G:\P\rightarrow \Grp$. Then
\begin{itemize}
\item The limiting cone $\tau:G\Rightarrow \limdir G$ is monic
(use that the forgetful functor $\Grp\rightarrow Set$ creates
filtered colimits, \cite[p.208]{MCL}). \item $\limdir_j H=0$ for
$j\geq 1$ because $\limdir:\Ab^\P\rightarrow \Ab$ is exact (see
\cite[2.6.15]{weibel}).
\end{itemize}
Thus, by the Theorem \ref{hocolim_technic}, we have that
$$\hocolim BG\cong B\limdir G.$$
\end{Ex}

\begin{Ex}\label{ex_locally finite group}
The last example applies to any locally finite group $G_0$: call
$\P$ to the poset category of its finite subgroups and
$G:\P\rightarrow \Grp$ to the monic functor which takes as values
the finite subgroups of $G_0$ and inclusion among them. Then
$G_0=\limdir_\P G$ and $\P$ is filtered and contractible (the
trivial group is an initial object). Thus:
$$\hocolim_{G\subseteq G_0,\textit{$G$ finite}} BG\cong BG_0$$
\end{Ex}

\begin{Ex}
Consider a group $G_0$ and the poset $\P$ (with the inclusion as
relation) of its normal finite $p$-subgroups for a fixed prime
$p$. If $H$ and $K$ are normal finite $p$-subgroups of $G_0$ then
$HK$ is a $p$-normal subgroup of $G_0$ too. This implies that $\P$
is directed. Moreover, $\P$ is contractible as $\{1\}\in \P$. Let
$G:\P\rightarrow \Grp$ the monic functor which takes each subgroup
of $\P$ to itself and inclusions to inclusions. Then by Example
\ref{ex_filtered_contractible}
$$\hocolim_{G\subseteq G_0,\textit{$G$ normal finite $p$-subgroup}} BG\cong B\limdir G.$$
If $G$ is finite then $\limdir G=O_p(G)$.
\end{Ex}

\begin{Ex}
Consider a \emph{finite} group $G_0$ and the poset $\P$ (with the
inclusion as relation) of the normal subgroups of $G_0$ which have
$p'$-quotient for a fixed prime $p$. Then $O^{p'}(G_0)$ is an
initial object of $\P$ and thus $\P$ is contractible. Moreover,
$G_0$ itself is a terminal object of $\P$ and thus $\P$ is
filtered. Let $G:\P\rightarrow \Grp$ the monic functor which takes
each subgroup of $\P$ to itself and inclusions to inclusions. Then
$\limdir G=G_0$ and, by Example \ref{ex_filtered_contractible},
$$\hocolim_{G\subseteq G_0,\textit{$G$ normal subgroup with $p'$ quotient}} BG\cong BG_0.$$
\end{Ex}

\section{Proof of the Theorem}\label{section_hocolim_thm}
\begin{Thm}
Let $G:\P\rightarrow \Grp$ be a functor and $\tau: G\Rightarrow
G_0$ be a cone, and call $F$ the homotopy fiber
$$F\rightarrow \hocolim BG\rightarrow BG_0.$$
Assume $\P$ is contractible and the cone $\tau:G\Rightarrow G_0$
is monic. Then
\begin{enumerate}
\item $\pi_1(F)=\ker(\limdir G\rightarrow G_0)$ \item
$\pi_0(F)=\coker(\limdir G\rightarrow G_0)$ \item
$H_j(F)=\limdir_{j-1} H$ for each $j\geq 2$,  where
$H:\P\rightarrow \Ab$ is a monic functor.
\end{enumerate}
\end{Thm}
\begin{proof}
Notice that because the spaces $BG_i$ are connected and we have a
pointed diagram $BG:\P\rightarrow SSets$, we can apply the Van
Kampen's spectral sequence \cite{stover} to obtain that
$\pi_1(\hocolim_* BG)=\limdir \pi_1(BG)=\limdir G $ and
$\pi_0(\hocolim_* BG)=0$.

The fibration
$$\P\rightarrow \hocolim BG\rightarrow \hocolim_* BG$$
gives, as $\pi_1(\P)=\pi_0(\P)=0$, $\pi_1(\hocolim BG)=\limdir G$
and $\pi_0(\hocolim BG)=0$. Then the fibration
$$
F\rightarrow \hocolim BG\rightarrow BG_0
$$
gives $\pi_1(F)=\ker(\lim G\rightarrow G_0)$ and
$\pi_0(F)=\coker(\lim G\rightarrow G_0)$.

It is readily checked that the usual construction of the homotopy
fiber $F_f$ of a map $f:A\rightarrow B$ is functorial on $A$ and
$B$. So the maps $B\tau_i:BG_i\rightarrow BG_0$ give a diagram
$F_.:\P\rightarrow SSet$ which take values the fibers $F_i$ of
$B\tau_i$:
$$F_i\rightarrow BG_i\rightarrow BG_0.$$
By \cite{cha}, the fiber $F$ is the homotopy colimit of the
diagram of the fibers $F_.:\P\rightarrow \SSet$:
$$F=\hocolim(F_.).$$

Consider the homology type first quadrant spectral sequence
(Bousfield-Kan) that converges to $H_*(F)=H_*(\hocolim F_i)$ and
such that $E^2_{p,q}=\limdir_p H_q(F_.)$. This spectral sequence
describes the homology of the fiber $F$. Consider again the
diagram $F_.$. It can be pointed taking $(1,c_1)\in F_i$ for each
$i\in \Ob(\P)$. Then we have natural transformations of functors
from $\P$ to $SSet_*$
$$F.\Rightarrow BG_.\Rightarrow BG_0.$$

For each $i\in \Ob(\P)$  the pointed fibration
$$F_i\rightarrow BG_i\rightarrow BG_0$$
gives the homotopy long exact sequence
$$.. 0\rightarrow \pi_1(F_i)\rightarrow \pi_1(BG_i)\rightarrow \pi_1(BG_0)
\rightarrow \pi_0(F_i)\rightarrow \pi_0(BG_i)\rightarrow
\pi_0(BG_0).$$ As the spaces $G_i$ and $\limdir G$ are discrete we
obtain that $\pi_j(F_i)=0$ for every $j\geq 2$.

Moreover, because $BG_i$ and $BG_0$ are connected we have:
$$0\rightarrow \pi_1(F_i)\rightarrow G_i \stackrel{\tau_i}\rightarrow \limdir G
\stackrel{p_i}\rightarrow \pi_0(F_i)\rightarrow 0.$$ This is a
exact sequence with three groups and a set. We can identify
$\pi_1(F_i)=\pi_1(F_i,(1,c_1))$ with $\ker(\tau_i)$ and write
\begin{equation}\label{4termseq}
0\rightarrow \ker(\tau_i)\rightarrow G_i
\stackrel{\tau_i}\rightarrow \limdir G \stackrel{p_i}\rightarrow
\pi_0(F_i)\rightarrow 0.
\end{equation}
If $\xi_i\in \pi_0(F_i)$ denotes the connected component of
$(1,c_1)\in F_i$ then the exactness at $\limdir G$ states that
$$p_i^{-1}(\xi_i)=\im(\tau_i).$$

Because the long homotopy exact sequence is natural
$$
\xymatrix{ 0\ar[r]&\ker(\tau_i)\ar[r]\ar[d]& G_i
\ar[r]^{\tau_i}\ar[d]&\limdir
G \ar@{=}[d]\ar[r]^{p_i}& \pi_0(F_i)\ar[d]\ar[r]& 0\\
0\ar[r]&\ker(\tau_j)\ar[r]& G_j \ar[r]^{\tau_j}&\limdir G
\ar[r]^{p_j}& \pi_0(F_j)\ar[r]& 0.  }
$$ then we have
in fact a exact sequence of functors
$$0\Rightarrow \ker(\tau_.)\Rightarrow G_. \stackrel{\tau_.}\Rightarrow \limdir G
\stackrel{p_.}\Rightarrow \pi_0(F_.)\Rightarrow 0$$


If the cone $\tau:G\Rightarrow \limdir G$ is monic then we have a
short exact sequence of functors
$$0\Rightarrow G_. \stackrel{\tau_.}\Rightarrow \limdir G
\stackrel{p_.}\Rightarrow \pi_0(F_.)\Rightarrow 0$$ and the spaces
$F_.$ are discrete. The Bousfield-Kan homology spectral sequence
reduces to
$$H_p(F)=\limdir_p H_0(F_.)$$
for $p\geq 0$.

Applying the functor free abelian group $\Z$ to the short exact
sequence above we obtain a sequence
$$\Z[G_.] \Rightarrow \Z[\limdir G]
\stackrel{\Z[p_.]}\Rightarrow \Z[\pi_0(F_.)]\Rightarrow 0$$ and
taking kernels we have a short exact sequence of functors
$$0\Rightarrow H \Rightarrow \Z[\limdir G]
\Rightarrow H_0(F_.)\Rightarrow 0$$ in $\Ab^\P$, where
$H(i)=\ker(\Z[\limdir G] \Rightarrow H_0(F_i))$, which is a
subgroup of the free abelian group $\Z[\limdir G]$ and so it is
free abelian too. Notice that the functor $H:\P\rightarrow \Ab$ is
monic: the arrow $H(i_1\rightarrow i_2)$ is the inclusion of the
subgroup $H(i_1)$ into the subgroup $H(i_2)$. Also it is clear
that the functor $H_0(F_\cdot)$ is epic: it is deduced from the
short exact sequence above where the middle functor $\Z[\limdir
G]$ is constant and takes the value $1_{\Z[\limdir G]}$ on
morphisms.

The long exact sequence of derived limits $\limdir_i$ for the
short exact sequence above gives
$$.. \rightarrow \limdir_1 H \rightarrow
\limdir_1  \Z[\limdir G] \rightarrow \limdir_1 H_0(F_.)
\rightarrow \limdir_0 H \rightarrow \limdir_0 \Z[\limdir G]
\rightarrow \limdir_0 H_o(F_.)\rightarrow 0$$ and recalling the
definition of homology of a simplicial set
$$.. \rightarrow \limdir_1 H \rightarrow
H_1(\P, \Z[\limdir G]) \rightarrow \limdir_1 H_0(F_.) \rightarrow
\limdir_0 H \rightarrow H_0(\P,\Z[\limdir G]) \rightarrow
\limdir_0 H_0(F_.)\rightarrow 0$$

Assume $\P$ is contractible. The we obtain
$$H_p(F)=\limdir_{p-1} H$$
for $p\geq 2$ by using the long exact sequence above.
\end{proof}

For convenience we give an explicit description of the functor
$H:\P\rightarrow \Ab$. Consider again the short exact sequence
with two groups and one set for any $i\in \Ob(\P)$
$$0\rightarrow G_i \stackrel{\tau_i}\rightarrow \limdir G
\stackrel{p_i}\rightarrow \pi_0(F_i)\rightarrow 0.$$

There is an action of $\limdir G=\pi_1(B\limdir G)$ on
$\pi_0(F_i)$ such that $p_i$ is $\limdir G$-equivariant when
$\limdir G$ is given the left action on itself, i.e., such that
$$
p_i(g\cdot g')=g\cdot p_i(g')
$$
for any $g,g'\in \limdir G$. Moreover, this action is natural in
the sense that for any arrow $i\rightarrow j$ in $\P$ we have
$$
\pi_0(i\rightarrow j)(g\cdot \eta)=g\cdot \pi_0(i\rightarrow
j)(\eta)
$$
for each $g\in \limdir G$ and $\eta\in \pi_0(F_i)$.

As $BG_i$ is path-connected the action of $\limdir G$ on
$\pi_0(F_i)$ is transitive. Fix, for each $\eta\in \pi_0(F_i)$, an
element $g_\eta\in \limdir G$ such that
$$g_\eta\cdot \xi_i=\eta$$
and $g_{\xi_i}=1$, where $\xi_i\in\pi_0(F_i)$ is the connected
component of $(1,c_1)$. Notice that $p_i(g_\eta)=\eta$. Recall
that the exactness of the short exact sequence above means that
$p_i^{-1}(\xi_i)=G_i\subseteq \limdir G$. Then
$$
p_i^{-1}(\eta)=g_\eta\cdot G_i
$$
and
$$
\limdir G=\bigcup_{\eta_\in\pi_0(F_i)} g_\eta \cdot G_i
$$
where the union is disjoint. In fact, $\pi_0(F_i)$ can be
identified with the set of cosets $G_i\backslash G$
\cite[p.73]{adem-milgram}, but we do not do it here. From this it
is straightforward that
$$
H(i)=\bigoplus_{\eta_\in\pi_0(F_i)} g_\eta \cdot IG_i,
$$
where $IG_i$ is the augmentation ideal of $G_i$, which is a free
abelian group with basis $\{g-1\}_{1\neq g\in G_i}$. Then the set
$$
B_i=\{\textit{$g_\eta\cdot g-g_\eta$, $1\neq g\in G_i$,
$\eta\in\pi_0(F_i)$}\}
$$
is a basis for the the free abelian group $H(i)$.

\section{Another example.}\label{section_poset}
A. Libman proposed the following example as one in which $\hocolim
BG\cong B\limdir G$ does not hold in spite the category $\P$ is
contractible and the cone $\tau$ is monic: Consider the product
category of the ``pushout category" $c\leftarrow a \rightarrow b$
with itself and denote it by $\P$. Then consider for any fixed
group $G_0$ the functor $G:\P\rightarrow \Grp$ which takes value
the trivial group $1$ on $(a,a)$ and value $G_0$ on the rest.

We compute the fiber $F$
$$F\rightarrow \hocolim BG\rightarrow B\limdir G=BG_0$$
by means of the tools developed in earlier sections:
\begin{Ex}\label{Ex_Libman}
Consider the contractible category $\P$ and the functor
$G:\P\rightarrow \Grp$ defined above. By Theorem
\ref{hocolim_technic} the fiber $F$
$$F\rightarrow \hocolim BG\rightarrow BG_0$$
is simply connected and $H_j(F)=\limdir_{j-1} H$ for $j\geq 2$,
where $H:\P\rightarrow \Ab$ is a monic functor. By dimensional
reasons $\limdir_{j} H=0$ for $j\geq 3$. The graded poset $\P$ has
shape:
$$
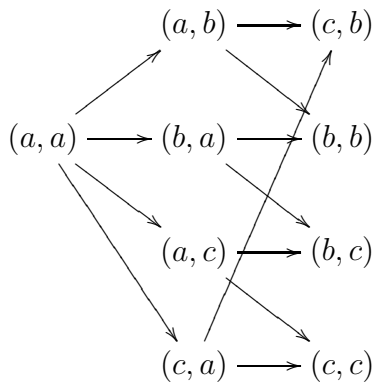
\xymatrix{
      & (a,b)\ar[r]\ar[rd] & (c,b) \\
(a,a)\ar[r]\ar[ru]\ar[rd]\ar[rdd] & (b,a)\ar[r]\ar[rd] & (b,b)\\
      & (a,c)\ar[r]\ar[rd] & (b,c)\\
      & (c,a)\ar[r]\ar[ruuu] & (c,c)\\
 }
$$
and the functor $G:\P\rightarrow \Grp$ takes values
$$
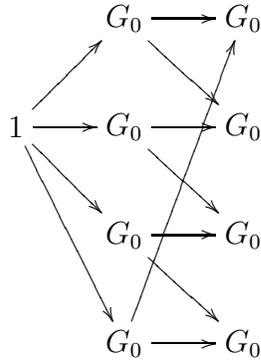
\xymatrix{
      & G_0\ar[r]\ar[rd] & G_0 \\
1\ar[r]\ar[ru]\ar[rd]\ar[rdd] & G_0\ar[r]\ar[rd] & G_0\\
      & G_0\ar[r]\ar[rd] & G_0\\
      & G_0\ar[r]\ar[ruuu] & G_0\\
 }
$$
The functor $H:\P\rightarrow \Ab$ is
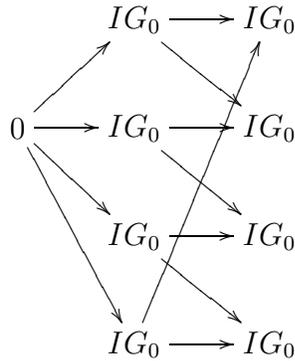
$$
\xymatrix{
      & IG_0\ar[r]\ar[rd] & IG_0 \\
0\ar[r]\ar[ru]\ar[rd]\ar[rdd] & IG_0\ar[r]\ar[rd] & IG_0\\
      & IG_0\ar[r]\ar[rd] & IG_0\\
      & IG_0\ar[r]\ar[ruuu] & IG_0\\
 }
$$
The functor $K_1$ (Section \ref{section_liminvII}) is isomorphic
to
$$
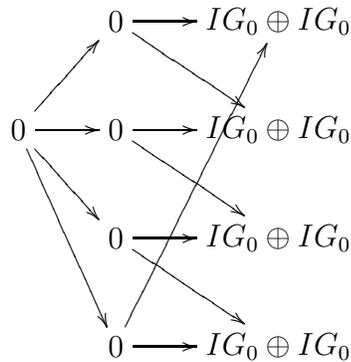
\xymatrix{
      & 0\ar[r]\ar[rd] & IG_0\oplus IG_0 \\
0\ar[r]\ar[ru]\ar[rd]\ar[rdd] &0\ar[r]\ar[rd] & IG_0\oplus IG_0\\
      & 0\ar[r]\ar[rd] & IG_0\oplus IG_0\\
      & 0\ar[r]\ar[ruuu] & IG_0\oplus IG_0\\
 }
$$
This implies that the functor $K_2$ is identically zero, and thus
$H_3(F)=\limdir_2 H=\ker\{\limdir K_2\rightarrow \limdir
K'_1\}=0$. The higher limit $\limdir_1 H$ corresponds to pairs
$(x_i,y_i)\in IG_0\oplus IG_0$ for $i=1,2,3,4$ with
$$x_i+y_i=0$$
and
$$x_i+x_{i+1}=0$$
for $i=1,2,3,4$. Then it is clear that $\limdir_1 H\cong IG_0$ and
thus we have a fibration
$$\bigvee_{\alpha\in G_0\setminus\{1\}} (S^2)_\alpha\rightarrow \hocolim BG\rightarrow BG_0$$
\end{Ex}



\end{document}